\documentclass{article}
\usepackage{amsmath}
\usepackage{amsthm,mathdots,cite}
\usepackage{amssymb}
\usepackage[ansinew]{inputenc}
\usepackage{url}
\usepackage[all,knot,arc,curve,color,frame]{xy}
\input xy
\xyoption{all}
\usepackage{extarrows}
\usepackage{tikz}
\usetikzlibrary{arrows}
\usepackage{comment}
\usepackage{eucal}
\usepackage{bm}
\usepackage{enumitem}
\usepackage{float}
\usepackage{chngcntr}
\counterwithin{figure}{section}

\theoremstyle{plain}
    \newtheorem{lemma}{Lemma}[section]
    
    \newtheorem{theorem}[lemma]{Theorem}
    \newtheorem{cor}[lemma]{Corollary}
    
    \newtheorem{prop}[lemma]{Proposition}
\theoremstyle{definition}
    
    \newtheorem{defi}[lemma]{Definition}
    \newtheorem{example}[lemma]{Example}
    \newtheorem{nota}[lemma]{Notation}
    \newtheorem{rem}[lemma]{Remark}
    
    \newtheorem{alg}[lemma]{Algorithm}

\DeclareMathOperator\fix{fix}
\DeclareMathOperator\im{im}

\DeclareMathOperator\hei{ht}

\DeclareMathOperator\dst{dst}

\DeclareMathOperator\length{length}
\DeclareMathOperator\lev{level}

\DeclareMathOperator\mrp{merge}

\newcommand{\al}{\alpha}
\newcommand{\bt}{\beta}
\newcommand{\gam}{\gamma}

\newcommand{\thet}{\theta}

\newcommand{\lam}{\lambda}

\newcommand{\vep}{\varepsilon}

\newcommand\tx{T(X)}
\newcommand\tn{T_n}
\newcommand\cn{\mathcal{C}_{n}}
\newcommand\cnine{\mathcal{C}_{9}}
\newcommand\xn{X_n}
\newcommand\lan{\langle}
\newcommand\ran{\rangle}

\newcommand\kjpq{K^j_{p,q}}

\newcommand\btm{B_t^m}
\newcommand\dtm{D_t^m}

\newcommand\ncn{N(\cn)}

\newcommand{\imp}{\Rightarrow}
\newcommand{\miff}{\Leftrightarrow}

\newcommand{\ga}{\Gamma}

\newcommand\ind{d^{-}}
\newcommand\od{d^{+}}
\newcommand\tr{\mathcal{T}}
\newcommand\pth{\mathcal{P}}
\newcommand\qth{\mathcal{Q}}
\newcommand\str{\mathcal{S}}

\newcommand\joi{\lor}

\newcommand\ceid{C_{id}(\vep)}
\newcommand\cenil{C_{nil}(\vep)}
\newcommand\ceidmu{C_{id}(\mu)}
\newcommand\cenilmu{C_{nil}(\mu)}
\newcommand\cnst{\cn^*}

\newcommand\lamp{\lambda^+}
\newcommand\tzppx{\tr^{0r}(\pth_x)}
\newcommand\toppx{\tr^{1r}(\pth_x)}
\newcommand\tzmpx{\tr^{0l}(\pth_x)}
\newcommand\tompx{\tr^{1l}(\pth_x)}

\newcommand\tppx{\tr^r(\pth_x)}
\newcommand\tmpx{\tr^l(\pth_x)}
\newcommand\tppy{\tr^r(\pth_y)}

\newcommand\tppz{\tr^r(\pth_z)}

\newcommand\tmpf{\tr^l(\pth_{\!f})}
\newcommand\tppf{\tr^r(\pth_{\!f})}
\newcommand\tmpfpr{\tr^l(\pth^{(f'\!\!,x)})}
\newcommand\tppfpr{\tr^r(\pth^{(f'\!\!,x)})}
\newcommand\tppa{\tr^r(\pth_a)}
\newcommand\tppe{\tr^r(\pth_e)}
\newcommand\cid{C_{id}}
\newcommand\cnl{C_{nil}}

\title{Centralizers in Catalan Monoids}

\author{Janusz Konieczny\footnote{Corresponding author.}\\
		{\small Department of Mathematics, University of Mary Washington}\\
	{\small Fredericksburg, Virginia 22401, jkoniecz@umw.edu}
	\and Emrah Korkmaz\\
{\small Department of Mathematics, \c{C}ukurova University}\\
{\small Adana, Sar\i\c{c}am 01330, ekorkmaz@cu.edu.tr}}
\date{}
\begin{document}
\maketitle

\begin{abstract}
\vskip 2mm
Fix a chain $\xn=\{1<\cdots<n\}$, where $n$ is a positive integer. 
The \emph{Catalan monoid} $\mathcal C_n$ consists of all decreasing and
order-preserving transformations of $X_n$, with function composition
as the operation.
For $\bt\in\cn$, its \emph{centralizer} is $C(\bt)=\{\al\in\cn:\al\bt=\bt\al\}$.

We characterize the elements of $C(\vep)$, $C(\mu)$, and $C(\bt)$, where $\vep\in\cn$ is an idempotent,
$\mu\in\cn$ is a nilpotent, and $\bt\in\cn$ is a general element. 
For an idempotent $\vep\in\cn$, we derive formulas for
the cardinality of $C(\vep)$, as well as formulas for the numbers of idempotents
and nilpotents of $C(\vep)$. We relate indexing sets arising in these
formulas to Dyck paths and Narayana numbers.

For a nilpotent $\mu\in\mathcal C_n$, we present algorithms that construct
$C(\mu)$, its subset of idempotents, and its subsemigroup of nilpotents,
thereby determining their cardinalities.
We also describe the nilpotents whose centralizers consist only of idempotents, of nilpotents together with the identity, 
or of nilpotents and idempotents. We extend these descriptions to arbitrary elements of $\cn$.

To obtain our results, we represent $\bt\in\cn$ as a sequence $\lan\tr_1,\ldots,\tr_m\ran$ of rooted trees.
The sets of vertices of the trees form a partition of $\xn$. Each tree is level-ordered, which means that it can be drawn
in such a way that moving from the root down, level by level, from left to right at each level, we pass through 
an increasing sequence of consecutive
elements of $\xn$.

\vskip 2mm

\noindent\emph{$2020$ Mathematics Subject Classification\/}. {20M20, 05A15, 05A19.}

\end{abstract}

\section{Introduction}\label{sint}
\setcounter{equation}{0}

For a semigroup $S$ and an element $a\in S$, the set $C(a)=\{x\in S: ax=xa\}$ is called the \emph{centralizer} of $a$ in $S$. It is clear that $C(a)$ is a subsemigroup of $S$. If $S$ has an identity $1$, then $C(1)=S$, and if $S$ has a zero $0$, then $C(0)=S$. Thus, centralizers can be regarded as generalizations of $S$. It is therefore of interest to find out which of the methods used to study $S$ can be extended to the centralizers of its elements, and how these centralizers differ from $S$ as semigroups. Centralizers also appear in other areas of research, for example,
in the study of automorphism groups of semigroups \cite{ArKo2009} and in the study of commuting graphs of semigroups
\cite{ArKiKo2011,ArBeKo2015,GiKu2016,Ku2021} and groups \cite{BaBuPeRo2003,Bu2006,IrJa2008}.

For an arbitrary set $X$, denote by $\tx$ the \emph{full transformation semigroup} on $X$, that is, the semigroup of all functions from $X$ to $X$ under the operation of function composition. We compose functions from left to right.
When $X$ is finite with $n$ elements, we take $\xn=\{1,\ldots,n\}$, ordered by $1<\cdots<n$, and denote $T(\xn)$ by $\tn$.
Numerous papers have been published on centralizers in transformation semigroups (see, for example, \cite{We1960,No1963,Hi1988,LiKo1995,KoLi1998,Ko1999,Ko2002,Ko2004a}). Green's relations in the centralizer of $\al$ have been determined when $\al$ is a permutation \cite{Ko1999}, an idempotent \cite{Ko2002}, or an injective nilpotent \cite{Ko2004a}. The centralizer of an arbitrary $\al\in\tx$ was described in \cite{ArKo2013}. For an infinite set $X$, the centralizers of idempotents in $\tx$ were investigated in \cite{ArKo2003,ArKo2004a,Sz2005}, and the cardinalities of
$C(\al)$ for certain types of $\al\in\tx$ were determined for a countable set $X$ in \cite{Ko2004b,Ko2006,Ko2009}.
Centralizers have also been studied in the symmetric inverse semigroup \cite{Li1996,Ko2013,Ko2018,Ko2020a}.
For $\al\in\tx$, the centralizer $C(\al)$ is the semigroup of endomorphisms of the functional digraph of $\al$ (see 
Definition~\ref{dfgr} and Lemma~\ref{lend}); we will use this description throughout the paper.

A transformation $\al\in\tn$ is called \emph{order-preserving} if $x<y$ implies $x\al\leq y\al$ for all $x,y\in\xn$, and
\emph{decreasing} if $x\al\leq x$ for all $x\in\xn$. We denote by $\cn$ the subsemigroup of $\tn$ consisting of all
order-preserving and decreasing transformations. The semigroup $\cn$ is a monoid with a zero, and its cardinality is
$\frac{1}{n+1}\binom{2n}{n}$, the $n^{\mathrm{th}}$
\emph{Catalan number} \cite{Hi1993b}. For this reason, $\cn$ is known in the literature as the \emph{Catalan monoid} \cite{So1996,MaSt2012}. Its combinatorial properties were studied in
\cite{AyAyKo2011,GaMa2009,Hi1993a,Hi1993b,LaUm2006,YaKo2020b}, and its algebraic properties, including presentations,
nilpotent and isolated subsemigroups, and ranks, were investigated in \cite{So1996,Hi1993a,KoAy2021,Ko2022,YaKo2020c,MaSt2012,MaGr2014}. The zero-divisor graph of $\cn$ was studied in \cite{To2021}.

In this paper, we study the centralizers of elements of the Catalan monoids. For $\bt\in\cn$, the centralizer $C(\bt)=\{\al\in\cn:\al\bt=\bt\al\}$ consists of those endomorphisms of the functional digraph $\ga(\bt)$ of $\bt$ that are also order-preserving and decreasing. 

In Section~\ref{scnasgraph}, we introduce the concept of a \emph{level-ordered}
tree (Definition~\ref{dlevord}) and show that each $\bt\in\cn$ can be represented
as a chain $\lan\tr_1<\cdots<\tr_m\ran$ of such trees (Proposition~\ref{paincn}). This representation is the main tool
we use to obtain our results.

Section~\ref{side} is devoted to the study of centralizers of idempotents of $\cn$. For an idempotent $\vep\in\cn$, we describe $C(\vep)$ 
(Theorem~\ref{tcide}) and derive formulas for $|C(\vep)|$ and for the numbers of idempotents and nilpotents in $C(\vep)$
(Theorems~\ref{tcidecard}, ~\ref{tceide}, and ~\ref{tcenil}).
The indexing set in the formula for $|C(\vep)|$ is in one-to-one correspondence with Dyck paths having a given number
of valleys (Proposition~\ref{pbtm}). Its cardinality is therefore given by Narayana numbers.

In Section~\ref{snil}, we describe $C(\mu)$ for a nilpotent $\mu\in\cn$ (Theorem~\ref{tcnilg}). 
Formulas for the cardinalities of $C(\mu)$, $\cid(\mu)$, and  $\cnl(\mu)$,
analogous to those for idempotents, do not seem feasible. Instead, we present algorithms for constructing $C(\mu)$, its subset $\cid(\mu)$ of idempotents, and its subsemigroup $\cnl(\mu)$ of nilpotents (Algorithms~\ref{acentmu}, \ref{acentnilmu}, and~\ref{acentidmu}), which
enable us to determine their cardinalities. 

Motivated by Examples~\ref{ecentmu}, \ref{ecentnilmu}, and~\ref{ecentidmu}, we also determine the nilpotents $\mu$ for which $C(\mu)$ consists only of idempotents (this happens only in a trivial case), of nilpotents together with the identity (Theorem~\ref{tiii}), or of nilpotents and idempotents
(Theorem~\ref{tiv}, which turned out to be the most challenging result of the paper).

In Section~\ref{sgen}, we describe $C(\bt)$ for an arbitrary $\bt\in\cn$ (Theorem~\ref{tcgen}), thereby generalizing Theorems~\ref{tcide} and~\ref{tcnilg}. We also extend Theorems~\ref{tiii} and~\ref{tiv}, obtained for nilpotents,
to non-nilpotent elements. More specifically, 
we show that, for $n\geq3$, if $C(\bt)$ contains no idempotents other than the identity and the zero, then $\bt$ is a nilpotent (Theorem~\ref{tgeniii});
and we determine the non-nilpotent elements $\bt\in\cn$ for which $C(\bt)$ consists of nilpotents and idempotents (Theorem~\ref{tgeniv}). Such elements occur only when $n\leq5$: four in $\mathcal{C}_4$ and one in $\mathcal{C}_5$ (see Figures~\ref{fig4} and~\ref{fig5}).

\section{Elements of $\cn$ as digraphs}
\label{scnasgraph}
\setcounter{equation}{0}
In this section, we introduce notation used throughout the paper and  show how to represent  each transformation $\al\in\cn$ by a chain of rooted
trees.

For an integer $m\geq1$, we set $[m]=\{1,\ldots,m\}$.
Let $\al\in\tn$. We denote the \emph{kernel}, \emph{image}, and \emph{fix} of~$\al$ by $\ker(\al)=\{(x,y)\in\xn\times\xn:x\al=y\al\}$,
$\im(\al)=\{x\al:x\in\xn\}$, and $\fix(\al)=\{x\in\xn:x\al=x\}$, respectively.
It is customary to write $\al\in\tn$ as
\begin{equation}\label{eqaltn}
\al=\begin{pmatrix}
A_1 & \ldots & A_m\\
a_1 & \ldots & a_m
\end{pmatrix},
\end{equation}
where $m\geq1$, $\im(\al)=\{a_1<\cdots<a_m\}$ and $A_p=a_p\al^{-1}$ for each $p\in[m]$. Note that $\{A_1,\ldots,A_m\}$ is the partition of $\xn$ induced by $\ker(\al)$.

Let $S$ be a semigroup. The set of idempotent elements of $S$ is denoted by $E(S)$, that is, $E(S)=\{\,e\in S: e^2=e\,\}$.
If $S$ has a zero, its set of nilpotent elements is denoted by $N(S)$, that is, 
$N(S)=\{\,a\in S: a^k=0\,\mbox{ for some integer $k\geq1$}\}$. A finite semigroup $S$ with zero is nilpotent, that is, $S^m=\{0\}$ for some positive integer $m$, if and only if
its only idempotent is zero \cite[Prop.~8.1.2]{GaMa2009}. 
It is well known that for every $\al\in\tn$, $\al$ is an idempotent if and only if $\fix(\al)=\im(\al)$.
As shown in \cite[Lemma~1.4]{LaUm2004}, an element $\al$ of $\cn$ is a nilpotent if and only if $\fix(\al)=\{1\}$.

A \emph{directed graph} (or a \emph{digraph}) is a pair $\Gamma=(V,E)$, where $V$ is a finite nonempty set
and~$E$ is a binary relation on $V$. Any element $x\in V$ is called a \emph{vertex} of $\Gamma$,
and any pair $(x,y)\in E$ is called an \emph{edge} of~$\Gamma$. We will write $x\to y$ if $(x,y)\in E$.
A \emph{walk} in $\ga$ is a sequence
$\pth=x_0\to\cdots\to x_k$, where $k\geq0$,
$x_0,\ldots,x_k\in V(\ga)$, and $x_i\to x_{i+1}$ is an edge
of $\ga$ for every $0\leq i<k$.
The walk starts at $x_0$ and ends at $x_k$; if $k=0$,
it consists of a single vertex.
The number $k$ is called
the \emph{length} of $\pth$, denoted by $\length(\pth)$. The walk $\pth$ is called a \emph{path} in $\ga$ from $x_0$ to $x_k$
if $x_0,\ldots,x_k$ are pairwise distinct.
A digraph $x_1\to\cdots\to x_k\to x_1$, where $k\geq1$ and $x_1,\ldots,x_k$ are pairwise distinct vertices of $\ga$ (which reduces to $x_1\to x_1$ if $k=1$) is called
a \emph{cycle} of length $k$ (or a \emph{$k$-cycle}) in $\ga$. We will denote a $k$-cycle by $(x_1\ldots\, x_k)$.

For $x\in V$, the \emph{in-degree} of $x$, denoted by $\ind(x)$, is the integer $\ind(x)=|\{y\in V:(y,x)\in E\}|$; and
the \emph{out-degree} of $x$, denoted by $\od(x)$, is the integer $\od(x)=|\{y\in V:(x,y)\in E\}|$.

A digraph is \emph{connected} if its underlying undirected graph is connected.
Every digraph is the disjoint union of its \emph{connected components},
which are its maximal connected subdigraphs.
If $\Gamma$ is a digraph, we write $V(\ga)$ for the set of vertices of $\ga$, and $E(\ga)$ for the set of edges of $\ga$.

Let $\ga=(V,E)$ be a digraph, and let $W$ be a nonempty subset of $V$. The \emph{subdigraph} of $\ga$ \emph{spanned by $W$} is the digraph $\ga_W=(W,F)$, where
$F=\{(x,y)\in W\times W: (x,y)\in E\}$.

Let $\ga_1=(V_1,E_1)$ and $\ga_2=(V_2,E_2)$ be digraphs.
A mapping $\phi:V_1\to V_2$
is called a \emph{homomorphism} from $\ga_1$ to $\ga_2$ if it preserves edges, that is, for all $x,y\in V_1$,
if $(x,y)\in E_1$, then $(x\phi,y\phi)\in E_2$ \cite{HeNe2004}.

We call a digraph $\Gamma$ a \emph{rooted tree} if it is connected, it has exactly one vertex $e$ 
(called the \emph{root} of the tree) such that $e\to e$ is an edge of $\ga$,
and for every vertex $x$ of $\Gamma$, $\od(x)=1$.
In this paper, by a ``tree'' we will always mean a ``rooted tree.''

Let $\tr$ be a tree with root $e$, and $x,y$ be vertices of $\tr$ such that there is a path (necessarily unique)
from $x$ to $y$. We call the length of this path the \emph{distance} from $x$ to $y$ and denote it by
$\dst(x,y)$. 
If $k=\dst(x,e)$, we say that $x$ is a \emph{$k$-level} vertex (or a vertex at \emph{level $k$}) in $\tr$. We denote by 
$\lev(x)$ the integer $k$ such that $x$ is at level~$k$ in $\tr$.
Note that the root $e$ is the only vertex at level $0$, and that every vertex $x\ne e$ such that $x\to e$ is at level $1$. 
The \emph{height} of $\tr$, denoted by $\hei(\tr)$, is the largest level of a vertex of $\tr$.
For $0\leq k\leq \hei(\tr)$, we denote by $V_k(\tr)$ the set of vertices of $\tr$ at level $k$. Note that $V_0(\tr)=\{e\}$.
If $k>\hei(\tr)$, we agree that $V_k(\tr)=\emptyset$.

Let $V$ be a finite nonempty set. Any transformation $\al\in T(V)$
can be represented by the digraph $\ga(\al)=(V,E)$, where for all $x,y\in V$, $(x,y)\in E$ if and only if $x\al=y$.

\begin{defi}\label{dfgr}
A digraph $\ga=(V,E)$ is called \emph{functional} if $\ga=\ga(\al)$ for some $\al\in T(V)$.
\end{defi}

The following lemma follows immediately from \cite{Hi1986,Hi1988}.

\begin{lemma}\label{lcfgr}
Let $\ga$ be a connected functional digraph. Then,
\begin{itemize} 
  \item [\rm(1)] there exists a unique cycle $(x_1\ldots\, x_k)$ in $\ga$;
  \item [\rm(2)] the unique cycle in $\ga$ is $(x_1)$ if and only if $\ga$ is a tree.
\end{itemize}
\end{lemma}

\begin{defi}\label{dint}
A subset $I$ of $\xn$ is called an \emph{interval} of $\xn$ if there are $i,p\in\xn$ such that $i+p-1\leq n$
and $I=\{i<i+1<\cdots<i+p-1\}$.
\end{defi}

\begin{defi}\label{dlevord}
Let $\tr$ be a tree such that $V(\tr)\subseteq\xn$, and let $h=\hei(\tr)$. We say that $\tr$ is \emph{level-ordered} if
$V(\tr)$ is an interval of $\xn$ and
for every $0\leq k\leq h$,
\begin{itemize} 
  \item [(a)] $V_k(\tr)$ is an interval of $\xn$;
  \item [(b)] for all $x,y\in V(\tr)$, if $k<h$, $x\in V_k(\tr)$, and $y\in V_{k+1}(\tr)$, then $x<y$;
  \item [(c)] for all $x,y,v,w\in V(\tr)$, if $x,y\in V_k(\tr)$, $x<y$, $x\to v$, and $y\to w$, then $v\leq w$.
\end{itemize}
\end{defi}

\begin{defi}\label{dordtree}
Let $\tr$ and $\str$ be trees such that $V(\tr),V(\str)\subseteq\xn$. We say that $\tr<\str$ if $x<y$ for all $x\in V(\tr)$ and $y\in V(\str)$.
\end{defi}
The binary relation $<$ from Definition~\ref{dordtree} is a strict partial order  (a binary relation that is irreflexive and transitive) on the set of trees whose vertices belong to $\xn$.
We write $\tr\leq\str$ if $\tr<\str$ or $\tr=\str$.

\begin{prop}\label{paincn}
Let $\al\in\tn$. Then, $\al\in\cn$ if and only if the following conditions are satisfied:
\begin{itemize} 
  \item [\rm(1)] for every 
  $x\in \xn$, the inverse image $x\al^{-1}$ is either empty or an interval of $\xn$;
  \item [\rm(2)] every connected component of $\ga(\al)$ is a level-ordered tree;
 \item [\rm(3)] for all distinct connected components $\tr$ and $\str$ of $\ga(\al)$, either $\tr<\str$ or $\str<\tr$.
\end{itemize}
\end{prop}
\begin{proof}
Suppose $\al\in\cn$.
To prove (1), let $x\in\xn$ be such that $x\al^{-1}\ne\emptyset$,
and let $x\al^{-1}=\{z_1<\cdots<z_p\}$. Suppose to the contrary that $z_{i+1}>z_i+1$ for some $i\in\{1,\ldots,p-1\}$. Then, since $\al$ is order-preserving, 
$x=z_i\al\leq(z_i+1)\al\leq z_{i+1}\al=x$. Thus, $(z_i+1)\al=x$, which is a contradiction since $z_i+1\notin x\al^{-1}$.
Hence, $x\al^{-1}$ is an interval of $\xn$.

We will prove (2) and (3) simultaneously. Let $\tr$ be a connected component of $\ga(\al)$. 
Since $\al$ is decreasing, every cycle of~$\tr$ is a loop: if $x_1\to\cdots\to x_r\to x_1$, then  $x_1\geq\cdots\geq x_r\geq x_1$,
which implies that the cycle is $(x_1)$. Thus, $\tr$ is a tree by Lemma~\ref{lcfgr}.

Let $h=\hei(\tr)$. We prove that $V_k(\tr)$ is an interval
of $\xn$ for every $0\leq k\leq h$.
The statement is immediate when $k=0$.
Assume $h\geq1$ and proceed by induction on $k$.
For $k=1$, $V_k(\tr)= V_1(\tr)=e\al^{-1}\setminus\{e\}$, where $e$ is the root of $\tr$. Thus, $V_1(\tr)$ is
an interval of $\xn$ since $e\al^{-1}$ is an interval by~(1), and $e$ is the smallest element of $e\al^{-1}$ because $\al$ is decreasing.
Let $1\leq k<h$ and suppose that
$V_k(\tr)=\{x_1<\cdots<x_s\}$ is an interval of $\xn$. By (1), $x_i\al^{-1}$ is empty or an interval of $\xn$, for every $i\in[s]$.
Since $k<h$, at least one of the sets $x_i\al^{-1}$ is
nonempty. If exactly one is nonempty, their union is an
interval by (1). Otherwise, consider $i,j\in[s]$ such that $i<j$, $x_i\al^{-1}$ and $x_j\al^{-1}$ are not empty, and
$x_r\al^{-1}=\emptyset$ for every $i<r<j$.
Let $z$ be the largest element of $x_i\al^{-1}$ and $w$ the smallest element of $x_j\al^{-1}$.
Suppose to the contrary that $w\ne z+1$. If $w\leq z$, then $x_j=w\al\leq z\al=x_i<x_j$, which is a contradiction.
Thus, $w>z+1$. Since $\al$ is order-preserving, 
$x_i=z\al\leq(z+1)\al\leq w\al=x_j$, and so $z+1\in x_r\al^{-1}$ for some $i\leq r\leq j$, which contradicts the choice of $i$, $j$, $z$, and $w$.
Thus, $w = z + 1$. It follows that $V_{k+1}(\tr)=x_1\al^{-1}\cup\cdots\cup x_s\al^{-1}$ is an interval of $\xn$. We have proved (a) of Definition~\ref{dlevord}.

Suppose $k<h$, $x\in V_k(\tr)$, and $y\in V_{k+1}(\tr)$.
Then, $y\al< y$ since $\al$ is decreasing and $y\al\in V_k(\tr)$. Thus, $x,y\al\in V_k(\tr)$ and $y\not\in V_k(\tr)$,
and so we cannot have $y\al<y\leq x$ since $V_k(\tr)$ is an interval of~$\xn$. Hence, $x<y$, which proves (b) of Definition~\ref{dlevord}.
Condition~(c) of Definition~\ref{dlevord} follows from the fact that $\al$ is order-preserving.

Let $\tr$ and $\str$ be distinct connected components of $\ga(\al)$. We already know that they are trees.
Let $e$ and $f$ be roots of $\tr$ and $\str$, respectively. Suppose $e<f$. 
Let $x\in V(\tr)$, $y\in V(\str)$, $k_1=\lev(x)$, $k_2=\lev(y)$, and $k=\max\{k_1,k_2\}$. Then, $x\al^k=e$ and $y\al^k=f$.
Suppose to the contrary that $x\geq y$.
Since $\al$ is order-preserving, $e=x\al^k\geq y\al^k=f$, which is a contradiction.
Hence, $x<y$, and so $\tr<\str$. If $f<e$, then $\str<\tr$ by symmetry. We have proved~(3).

Let $\tr_1,\ldots,\tr_m$ be the connected components (trees) of $\ga(\al)$ ordered in such a way
that $e_1<\cdots<e_m$, where $e_i$ is the root of $\tr_i$. 
By the foregoing argument, $\tr_1<\cdots<\tr_m$.
Since their vertex sets partition $\xn$, each $V(\tr_i)$
is an interval of $\xn$.
This concludes the proof of~(2).

Conversely, suppose (1)-(3) are satisfied. Let $x,y\in\xn$, with $x\in V(\tr)$ and $y\in V(\str)$, where $\tr$ and $\str$ are connected components
(trees) of $\ga(\al)$.
Suppose that $y=x\al$. Then, $\str=\tr$ by the definition of a connected component of a functional digraph.
If $x$ is the root of $\tr$, then $x=x\al=y$. Otherwise, $y=x\al<x$ by~(2) (see (b)~of Definition~\ref{dlevord}).
Thus, $\al$ is decreasing. 

Suppose $x<y$. Then, $\tr\leq\str$ by (3). Note that $x\al\in V(\tr)$ and $y\al\in V(\str)$. If $\tr<\str$, then $x\al<y\al$. Suppose $\tr=\str$.
Then, since $x<y$, we have $\lev(x)\leq\lev(y)$ by~(2), and so $\lev(x\al)\leq\lev(y\al)$. 
If $\lev(x\al)<\lev(y\al)$, then $x\al<y\al$ by~(2) (see (b)~of Definition~\ref{dlevord}). Suppose $\lev(x\al)=\lev(y\al) (= k)$.
If $k=0$, then $x\al=y\al=\mbox{the root of $\tr$}$. If $k>0$, then  $\lev(x)=\lev(y)$, and so $x\al\leq y\al$ by~(2)
(see (c)~of Definition~\ref{dlevord}). We have proved that $\al$ is order-preserving. Therefore, $\al\in\cn$.
\end{proof}

We will often use Proposition~\ref{paincn} in our proofs, not always explicitly.
For $\al\in\cn$, when we say \emph{trees of $\ga(\al)$}, we will mean the connected
components of $\ga(\al)$. In drawings, we omit arrowheads and the loop at the root. For each edge $x\to y$ with
$x\ne y$, we draw $x$ below $y$. The root is placed at the top, and vertices at the same level lie on a
horizontal line.

\begin{nota}\label{ntrcn}
Let $\al\in\cn$ with $\fix(\al)=\{e_1<\cdots<e_m\}$, where $e_1=1$. By Proposition~\ref{paincn}, the trees of $\ga(\al)$
can be ordered $\tr_1<\cdots<\tr_m$ (see Definition~\ref{dordtree}), where $e_i$ is the root of $\tr_i$ for each $i\in[m]$. We will write $\ga(\al)$ as the chain
\[
\ga(\al)=\{\tr_1<\cdots<\tr_m\}.
\]
\end{nota}

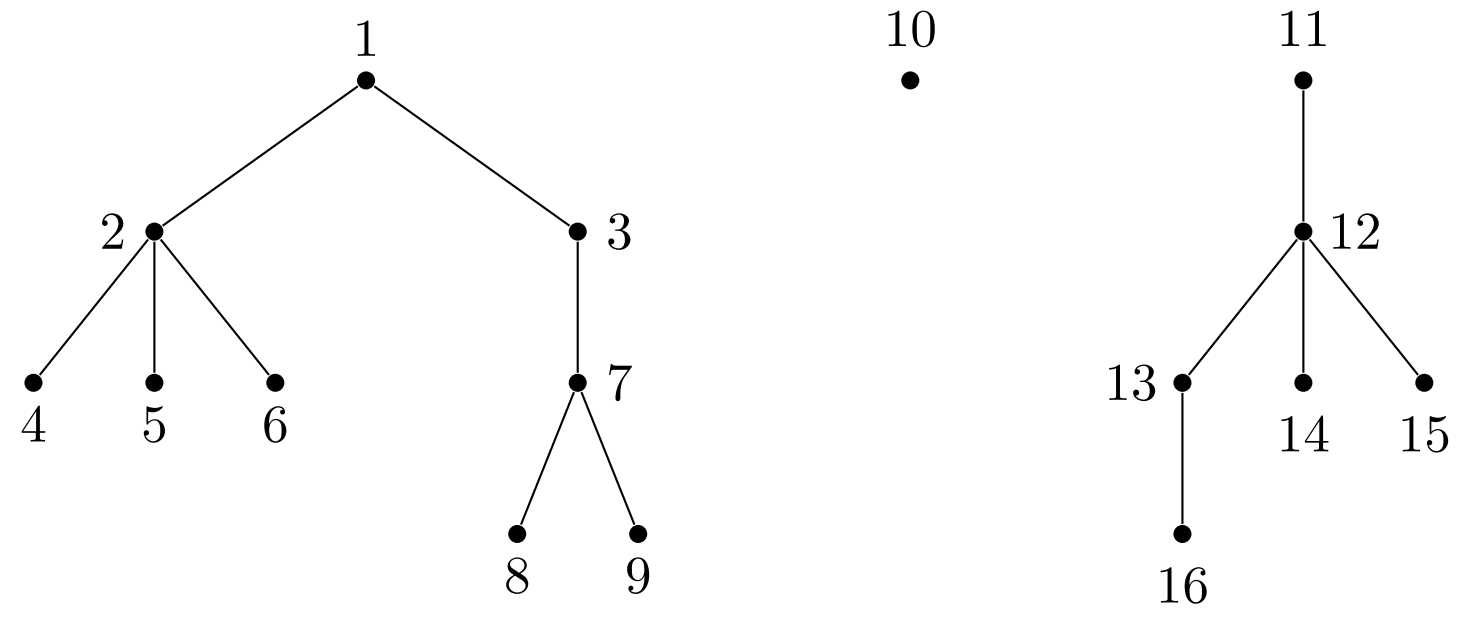
\begin{figure}[htbp]
	\centering
	\begin{tikzpicture}[
	level distance=1.0cm,
	every node/.style={circle, fill=black, inner sep=1.2pt}
	]
	\begin{scope}[level 1/.style={sibling distance=2.8cm},
	level 2/.style={sibling distance=0.8cm},
	level 3/.style={sibling distance=0.8cm}]
	\node[label=above:1] at (0,0) {}
	child { node[label=left:2] {}
		child { node[label=below:4] {} }
		child { node[label=below:5] {} }
		child { node[label=below:6] {} }
	}
	child { node[label=right:3] {}
		child { node[label=right:7] {}
			child { node[label=below:8] {} }
			child { node[label=below:9] {} }
		}
	};
	\end{scope}
	\node[label=above:10] at (3.6,0) {};
	\begin{scope}[level 1/.style={sibling distance=1.0cm},
	level 2/.style={sibling distance=0.8cm},
	level 3/.style={sibling distance=0.8cm}]
	\node[label=above:11] at (6.2,0) {}
	child { node[label=right:12] {}
		child { node[label=left:13] {}
			child { node[label=below:16] {} }
		}
		child { node[label=below:14] {} }
		child { node[label=below:15] {} }
	};
	\end{scope}
	\end{tikzpicture}
	\caption{The digraph of a transformation $\alpha \in \mathcal{C}_{16}$ with $\ga(\al)=\{\tr_1<\tr_2<\tr_3\}$.}
	\label{fig2}
\end{figure}

The digraph in Figure~\ref{fig2} represents
a transformation $\al\in\mathcal{C}_{16}$ with three connected components (trees). 
The middle component is a tree whose only vertex is its root (so $10\al=10$).

\begin{defi}\label{dconv}
Let $X=\{x_1<\cdots<x_k\}$, where $k\geq1$, be a chain of positive integers.
A sequence $\lan A_1,\ldots,A_m\ran$, where $1\leq m\leq k$, is called an \emph{ordered convex partition} of $X$
if $\{A_1,\ldots,A_m\}$ is a partition of $X$ and,
for every $1\leq i<m$, every element of $A_i$ is smaller
than every element of $A_{i+1}$.
\end{defi}
For example, $\lan\{1,2\},\{3,4,5\},\{6\}\ran$ is an ordered convex partition of $X_6$.

\begin{lemma}\label{lconv}
Let $X=\{x_1<\cdots<x_k\}$, where $k\geq1$, be a chain
of positive integers. For every $m\in[k]$, the number
of ordered convex partitions of $X$ into $m$ blocks
is $\binom{k-1}{m-1}$.
\end{lemma}
\begin{proof}
Let $m\in[k]$. There is a one-to-one correspondence between the ordered convex partitions $\lan A_1,\ldots,A_m\ran$ of $X$
and the chains $\{1=j_1<\cdots<j_m\}$ of elements of $[k]$ via the bijection 
\[
\lan A_1,\ldots,A_m\ran\to\{1=j_1<\cdots<j_m\},
\]
where for every $i\in[m]$, $x_{j_i}$ is the smallest element of $A_i$. 
Since $j_1=1$, the remaining $m-1$ indices can be chosen
from $\{2,\ldots,k\}$ in $\binom{k-1}{m-1}$ ways.
\end{proof}

For a finite set $A$ of positive integers with $A\ne\emptyset$, we will denote by 
$\min(A)$ and $\max(A)$ the smallest and largest, respectively, integer in $A$.
The following proposition follows immediately from Proposition~\ref{paincn}.
\begin{prop}\label{pcnel}
Let $\al\in \tn$ be written as in \eqref{eqaltn}. Then, $\al\in\cn$ if and only if 
$\lan A_1,\ldots,A_m\ran$ is an ordered convex partition
of $\xn$, and $a_i\leq\min(A_i)$ for every $i\in[m]$
(hence $a_1=1$).
\end{prop}

\section{Centralizers of idempotents of $\cn$}\label{side}
\setcounter{equation}{0}

In this section, we study the centralizers of idempotents in $\cn$.
For a given $\vep\in E(\cn)$,
we characterize the elements of $C(\vep)$ and present formulas for calculating
$|C(\vep)|$, $|\cid(\vep)|=|C(\vep)\cap E(\cn)|$ (the cardinality of the set of idempotents in $C(\vep)$), and
 $|\cnl(\vep)|=|C(\vep)\cap N(\cn)|$ (the cardinality of the semigroup of nilpotents in $C(\vep)$).

\subsection{Elements of $C(\vep)$}\label{side1}
Recall that a transformation $\vep\in\tn$ is idempotent
if and only if $\im(\vep)=\fix(\vep)$.
We use the following notation for an idempotent in $\cn$.

\begin{nota}\label{nide}
An idempotent $\vep\in\cn$ is completely determined by $\fix(\vep)=\{e_1<\cdots<e_m\}$. 
Let $\ga(\vep)=\{\tr_1<\cdots<\tr_m\}$ (see Notation~\ref{ntrcn}). Since $\vep\in E(\cn)$, each tree $\tr_i$
has height at most~$1$. For $i\in[m]$, set $E_i=e_i\vep^{-1}$. Then, $E_i=V(\tr_i)$ and, following \eqref{eqaltn}, we can write
\[
\vep=\begin{pmatrix}E_1  & \ldots & E_m \\ e_1  & \ldots & e_m\end{pmatrix}.
\]
By Proposition~\ref{pcnel},
$\lan E_1,\ldots,E_m\ran$ is an ordered convex partition of $\xn$. Since $\vep$ is decreasing and $e_i\in E_i$,
we have $e_i=\min(E_i)$ for every $i\in[m]$. Consequently, $E_i=\{e_i,e_i+1,\ldots,e_{i+1}-1\}$ if $i<m$, and
$E_m=\{e_m,e_m+1,\ldots,n\}$. (See Figure~\ref{fig3}.)
\end{nota}
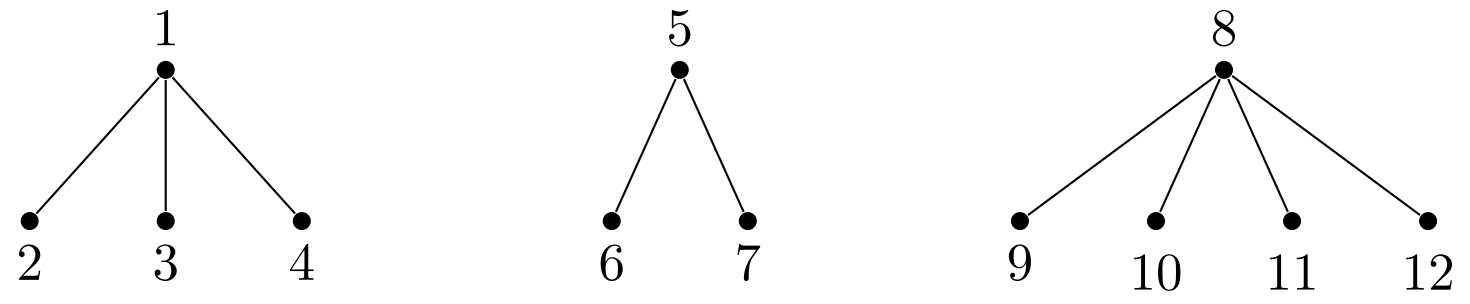
\begin{figure}[htbp]
	\centering
	\begin{tikzpicture}[
	level distance=1.0cm,
	level 1/.style={sibling distance=0.9cm},
	every node/.style={circle, fill=black, inner sep=1.2pt}
	]
	\node[label=above:1] at (0,0) {}
	child { node[label=below:2] {} }
	child { node[label=below:3] {} }
	child { node[label=below:4] {} };
	
	\node[label=above:5] at (3.4,0) {}
	child { node[label=below:6] {} }
	child { node[label=below:7] {} };
	
	\node[label=above:8] at (7,0) {}
	child { node[label=below:9] {} }
	child { node[label=below:10] {} }
	child { node[label=below:11] {} }
	child { node[label=below:12] {} };
	\end{tikzpicture}
	\caption{The digraph of an idempotent $\vep \in \mathcal{C}_{12}$.}
	\label{fig3}
\end{figure}

\begin{defi}\label{drelint}
Let $I=\{i<i+1<\cdots<i+p-1\}$ and $J=\{j<j+1<\cdots<j+q-1\}$ be intervals of $\xn$. We say that $I<J$ if $i+p-1<j$.
\end{defi}

The binary relation $<$ from Definition~\ref{drelint} is a strict partial order on the set of intervals of $\xn$.
We write $I\leq J$ if $I<J$ or $I=J$.

\begin{nota}\label{ncnij}
Let $I=\{i<i+1<\cdots<i+p-1\}$ and $J=\{j<j+1<\cdots<j+q-1\}$ be intervals of $\xn$, where $I\leq J$.
We  denote by $\cn(J,I)$ the set 
\[
\cn(J,I)=\{\al:J\to I: j\al=i\mbox{ and $\al$ is decreasing and order-preserving}\}.
\]
\end{nota}

The following theorem characterizes $C(\vep)$ for an arbitrary idempotent $\vep\in\cn$.

\begin{theorem}\label{tcide}
Let $\vep\in E(\cn)$ (written as in {\rm Notation~\ref{nide}}) and $\al\in\cn$. Then, $\al\in C(\vep)$ if and only if
 the following conditions are satisfied:
\begin{itemize} 
  \item [\rm(1)] $(\fix(\vep))\al\subseteq\fix(\vep)$;
  \item [\rm(2)] if $\{e_{i_1}<\cdots<e_{i_t}\}=(\fix(\vep))\al$, $\{j_1<\cdots<j_t\}$ is such that $j_r=\min\{j\in[m]:e_j\al=e_{i_r}\}$,
and $j_{t+1}=m+1$, then for all $j\in[m]$ and $r\in[t]$,
\begin{itemize} 
 \item [\rm(a)] $\al|_{E_j}\in\cn(E_j,E_{i_r})$ if $j=j_{r+1}-1$;
  \item [\rm(b)] $E_j\al=\{e_{i_r}\}$ if $j_r\leq j<j_{r+1}-1$.
\end{itemize}
\end{itemize} 
\end{theorem}

\begin{proof}
Let $\al\in C(\vep)$. Suppose $e_{i}\in \fix(\vep)$, where $i\in[m]$. Then, since $\vep\al=\al\vep$, we have
$(e_i\al)\vep=(e_i\vep)\al=e_i\al$,
and so $e_i\al\in\fix(\vep)$. This proves~(1).

Let $\{e_{i_1}<\cdots<e_{i_t}\}$, $\{j_1<\cdots<j_t\}$, and $j_{t+1}$ be as in~(2), and let $j\in[m]$ and $r\in[t]$. 
Since $\al$ is order-preserving, 
\begin{equation}\label{tcidee1}
e_j\al=e_{i_r}\miff j_r\leq j<j_{r+1}.
\end{equation}
Also, $i_1=j_1=1$, and $i_r\leq j_r$ for every $r\in[t]$,
since $e_{i_r}=e_{j_r}\al\leq e_{j_r}$.

Suppose that $j=j_{r+1}-1$ and $x\in E_j$, so  $j_r\leq j<j_{r+1}$ by the choice of $j_r$.
By \eqref{tcidee1}, $(x\al)\vep=(x\vep)\al = e_{j}\al = e_{i_r}$, which implies $x\al\in E_{i_r}$.
Thus, $E_j\al\subseteq E_{i_r}$, and so, since
$\al|_{E_j}$ is decreasing and order-preserving, $\al|_{E_j}\in\cn(E_j,E_{i_r})$
(see Notation~\ref{ncnij}). This proves~2(a).

Suppose that $j_r\leq j<j_{r+1}-1$ and $x\in E_j$.
Then, $e_{j_r}\leq x<e_{j+1}\leq e_{j_{r+1}-1}$. Since $\al$ is order-preserving,
$e_{i_r}=e_{j_r}\al\leq x\al\leq e_{j_{r+1}-1}\al=e_{i_r}$.
Hence, $E_j\al=\{e_{i_r}\}$, proving 2(b).

Conversely, suppose that (1) and (2) hold. Let $x\in\xn$, and choose $j\in[m]$ such that $x\in E_j$.
Since $j_1=1$ and $j_{t+1}=m+1$, there exists $r\in[t]$ such that $j_r\leq j<j_{r+1}$.
By \eqref{tcidee1} and (2), we have $e_j\al=e_{i_r}$ and $x\al\in E_{i_r}$. Therefore,
$x(\vep\al)=e_j\al=e_{i_r}=(x\al)\vep=x(\al\vep)$.
Therefore, $\vep\al=\al\vep$, that is, $\al\in C(\vep)$, as required.
\end{proof}

\subsection{Formula for calculating $|C(\vep)|$}\label{sidel2}
In this subsection, we present a formula for calculating the cardinality of $C(\vep)$, where $\vep\in E(\cn)$.  
We first establish the necessary counting lemmas and then
relate the indexing sets to Dyck paths and Narayana numbers.

We begin by calculating the cardinality of $\cn(J,I)$, where $I$ and $J$
are intervals of $\xn$ such that $I\leq J$.
The following lemma follows from the fact that $\cn(I,I)\cong\mathcal{C}_p$.
\begin{lemma}\label{lcnii}
Let $I=\{i<i+1<\cdots<i+p-1\}$ be an interval of $\xn$. Then, $\lvert\cn(I,I)\rvert=C_p=\frac1{p+1}\binom{2p}{p}$.
\end{lemma}

\begin{lemma}\label{lcnij}
Let $I=\{i<i+1<\cdots<i+p-1\}$ and $J=\{j<j+1<\cdots<j+q-1\}$ be intervals of $\xn$, where $i+p-1<j$. Then, 
$\lvert\cn(J,I)\rvert=\binom{p+q-2}{p-1}$.
\end{lemma}
\begin{proof} 
Let $t=\min\{p,q\}$ and fix $r\in[t]$. Then, every $\al\in\cn(J,I)$ of rank $r$ can be written as
\[
\al=\begin{pmatrix}A_1  & \ldots & A_r \\ a_1  & \ldots & a_r\end{pmatrix},
\]
where $\{i=a_1<\cdots<a_r\}$ is a chain of $I$ and $\lan A_1,\ldots,A_r\ran$ is an ordered convex partition of $J$. We can choose
$\{i=a_1<\cdots<a_r\}$ in $\binom{p-1}{r-1}$ ways and $\lan A_1,\ldots,A_r\ran$ in $\binom{q-1}{r-1}$ ways
(by Lemma~\ref{lconv}).
Thus, the number of $\al\in\cn(J,I)$ of rank $r$ is $\binom{p-1}{r-1}\binom{q-1}{r-1}$. Hence,
\[
\lvert\cn(J,I)\rvert
=\sum_{r=1}^{t}\binom{p-1}{p-r}\binom{q-1}{r-1}
=\binom{p+q-2}{p-1},
\]
where the last equality follows from Vandermonde's identity \cite[Problem~26, p.~167]{St1993}.
\end{proof}

Motivated by Lemmas~\ref{lcnii} and \ref{lcnij}, we define the integer $\kjpq$, where $j$, $p$, and $q$ are positive integers:
\begin{equation}\label{ekjpq}
\kjpq=\left\{
\begin{array}{ll}
\frac1{p+1}\binom{2p}{p} & \mbox{if $j=1$}\\[2mm]
\binom{p+q-2}{p-1} & \mbox{if $j>1$.}
\end{array}
\right.
\end{equation}

\begin{nota}\label{nbtm}
Let $m\geq1$ be an integer and $t\in[m]$. We denote by $\btm$ the set of all pairs 
\[
(\{i_1<\cdots<i_t\},\{j_1<\cdots<j_t\})
\]
such that $i_r,j_r\in[m]$ and $i_r\leq j_r$ for every $r\in[t]$, $i_1=1$, and $j_1=1$.
\end{nota}

\begin{theorem}\label{tcidecard}
Let $\vep\in E(\cn)$ be written as in {\rm Notation~\ref{nide}}. Then
\[
|C(\vep)|=\sum_{t=1}^m\sum_{(\{i_1<\cdots<i_t\},\{j_1<\cdots<j_t\})\in\btm}\prod_{r=1}^t K^{j_{r+1}-i_r}_{p_r,q_r},
\]
where $\btm$ is defined in {\rm Notation~\ref{nbtm}}, $j_{t+1}=m+1$, $p_r=|E_{i_r}|$, $q_r=|E_{j_{r+1}-1}|$, and $K^{j_{r+1}-i_r}_{p_r,q_r}$ is defined in \eqref{ekjpq}.
\end{theorem}

\begin{proof}
Fix $t \in [m]$ and a pair $(\{i_1 < \cdots < i_t\}, \{j_1 < \cdots < j_t\}) \in B^m_t$. Let $A$ be the set of all $\al\in\cn$ such that
for every $j\in[m]$, $\al|_{E_j}\in \cn(E_j,E_{i_r})$ if $j = j_{r+1} - 1$ for some $r\in[t]$,
and $\im(\al|_{E_j})=\{e_{i_r}\}$  if $j_r\leq j<j_{r+1} - 1$ for some $r\in[t]$. Then, by
Theorem~\ref{tcide}, $A$ is the set of all $\al\in C(\vep)$ such that $(\fix(\vep))\al=\{e_{i_1}<\cdots<e_{i_t}\}$
and, for every $\al\in A$, $\{j_1<\cdots<j_t\}$ is such that $j_r=\min\{j\in[m]:e_j\al=e_{i_r}\}$ for every $r\in[t]$.

For $j\in[m]$, let $A_j=\{\al|_{E_j}:\al\in A\}$. 
If $j=j_{r+1}-1$ and $j_{r+1}-i_r = 1$, for some $r\in[t]$, then $j=i_r$, and so
\[
|A_j| = |\cn(E_j,E_{i_r})| = |\cn(E_{i_r}, E_{i_r})| = \frac{1}{p_r + 1} \binom{2p_r}{p_r}=K^{j_{r+1}-i_r}_{p_r,q_r}
\]
by Lemma~\ref{lcnii} and~\eqref{ekjpq}.
If $j=j_{r+1}-1$ and $j_{r+1}-i_{r}>1$, for some $r\in[t]$, then $i_r<j_{r+1}-1$, and so
\[
|A_j| =|\cn(E_j,E_{i_r})|=|\cn(E_{j_{r+1}-1},E_{i_r})| = \binom{p_r+q_r-2}{p_r-1}=K^{j_{r+1}-i_r}_{p_r,q_r}
\]
by Lemma~\ref{lcnij} and~\eqref{ekjpq}. For any other value of $j$, $|A_j|=1$. 

The maps on the intervals $E_j$ can be chosen independently.
Indeed, choose maps $\al_j:E_j\to\xn$ satisfying the local
conditions defining $A$, and define $\al:\xn\to\xn$ by
$\al|_{E_j}=\al_j$ for every $j\in[m]$.
Each $\al_j$ is order-preserving and decreasing.
For each $r\in[t]$, the maps $\al_j$ with
$j_r\leq j<j_{r+1}-1$ are constant with value $e_{i_r}$,
whereas $\al_{j_{r+1}-1}$ takes values in $E_{i_r}$.
Since $e_{i_r}=\min E_{i_r}$, the map $\al$ is
order-preserving on
$E_{j_r}\cup\cdots\cup E_{j_{r+1}-1}$.
Moreover, $E_{i_1}<\cdots<E_{i_t}$, so $\al$ is
order-preserving on all of $\xn$.
It is also decreasing, and hence $\al\in\cn$.
By construction, $\al$ satisfies the conditions defining $A$,
so $\al\in A$.
Since the intervals $E_j$ partition $\xn$, this extension
is unique. Thus, the map $A\to A_1\times\cdots\times A_m$
is a bijection. Consequently,
\begin{equation}\label{ecidecard}
|A|=\prod_{j=1}^m |A_j|
=\prod_{r=1}^t K^{j_{r+1}-i_r}_{p_r,q_r}.
\end{equation}

Summing the rightmost expression in \eqref{ecidecard} over all pairs in $B^m_t$, we obtain the cardinality of the set of 
$\al\in C(\vep)$ such that $|(\fix(\vep))\al|=t$. Finally, 
summing over all $t \in [m]$, we obtain
\[
|C(\vep)|=\sum_{t=1}^{m} \sum_{(\{i_1<\cdots<i_t\},\{j_1<\cdots<j_t\})\in B_t^m} \prod_{r=1}^{t} K^{j_{r+1}-i_r}_{p_r,q_r},
\]
as required.
\end{proof}

\begin{example}\label{exm2}
We will apply Theorem~\ref{tcidecard} to an idempotent $\vep\in\cn$ with $m=|\fix(\vep)|=2$.
Let $\vep=\begin{pmatrix}E_1& E_2 \\ e_1 & e_2\end{pmatrix}$. Suppose $t=1$. Then, $\btm=B_1^2=\{(\{1\},\{1\})\}$ and $j_2=3$. Let $r\in[t]$, so $r=1$.
Then, $p_1=|E_1|=e_2-1$, $q_1=|E_2|=n-e_2+1$, and $j_2-i_1=3-1=2$. Hence,
\[
K_{p_1,q_1}^{j_2-i_1}=K_{p_1,q_1}^2=\binom{p_1+q_1-2}{p_1-1}=\binom{e_2-1+n-e_2+1-2}{e_2-1-1}=\binom{n-2}{e_2-2}.
\]
Suppose $t=2$. Then, $\btm=B_2^2=\{(\{1<2\},\{1<2\})\}$ and $j_3=3$. Let $r\in[2]$. If $r=1$, then
$p_1=|E_1|=e_2-1$, $q_1=|E_1|=e_2-1$, $j_2-i_1=2-1=1$,
and hence $K_{p_1,q_1}^{j_2-i_1}=C_{e_2-1}$.
If $r=2$, then
$p_2=|E_2|=n-e_2+1$, $q_2=|E_2|=n-e_2+1$, $j_3-i_2=3-2=1$, and so $K_{p_2,q_2}^{j_3-i_2}=C_{n-e_2+1}$.

By the foregoing calculations and Theorem~\ref{tcidecard}, we obtain
\[
|C(\vep)|=\binom{n-2}{e_2-2}+C_{e_2-1}C_{n-e_2+1}.
\]
\end{example}

We now determine $|\btm|$ using  \emph{lattice paths} \cite[p.~199]{St2024}. We will consider the lattice paths with steps $(0,1)$ and $(1,0)$. There are several definitions
of \emph{Dyck paths} \cite[Cor.~6.2.3(iv)(v)]{St2024}. Here, for an integer $m\geq1$, by a Dyck path we will mean a lattice path $L$ from $(0,0)$ to $(m,m)$ such that no point $(a,b)$ of~$L$
lies below the line $y=x$. A \emph{valley} of a Dyck path $L$ is a lattice point $(a,b)$ of $L$ such that
$(a,b+1)$ lies on $L$ and $(a,b-1)$ does not. (An analogous notion is that of a \emph{peak}.)
Under this convention, $(0,0)$ counts as a valley.
Thus, the number of valleys equals the number of maximal
vertical runs, and hence the number of peaks.
It is well known that the number of Dyck paths with $t$ valleys (or peaks) is equal to the \emph{Narayana number} $N(m,t)=\frac{1}{m} \binom{m}{t} \binom{m}{t-1}$ \cite[Ex.~6.36(a)]{St2024}.

\begin{figure}[H]

	\begin{tikzpicture}[thick]
	
	\tikzset{
		grid lines/.style={help lines, step=1cm, gray!30},
		diagonal/.style={gray, dashed},
		path/.style={blue, thick},
		valley/.style={red, fill=red, inner sep=2pt, opacity=0.8}
	}
	
	\begin{scope}[xshift=0cm]
	\draw[grid lines] (0,0) grid (3,3);
	\draw[diagonal] (0,0) -- (3,3);
	\draw[path]
	(0,0) -- (0,1) -- (0,2) -- (1,2) -- (1,3) -- (2,3) -- (3,3);
	\fill[valley] (0,0) circle (3pt) node[below right, black, font=\small] {Valley};
	\fill[valley] (1,2) circle (3pt) node[below right, black, font=\small] {Valley};
	\node at (1.5,-0.6) {(a) \((\{1,2\}, \{1,3\})\)};
	\end{scope}
	
	\begin{scope}[xshift=6cm]
	\draw[grid lines] (0,0) grid (3,3);
	\draw[diagonal] (0,0) -- (3,3);
	\draw[path]
	(0,0) -- (0,1) -- (1,1) -- (1,2) -- (1,3) -- (2,3) -- (3,3);
	\fill[valley] (0,0) circle (3pt) node[below right, black, font=\small] {Valley};
	\fill[valley] (1,1) circle (3pt) node[below right, black, font=\small] {Valley};
	\node at (1.5,-0.6) {(b) \((\{1,2\}, \{1,2\})\)};
	\end{scope}
	
	\begin{scope}[xshift=12cm]
	\draw[grid lines] (0,0) grid (3,3);
	\draw[diagonal] (0,0) -- (3,3);
	\draw[path]
	(0,0) -- (0,1) -- (0,2) -- (1,2) -- (2,2) -- (2,3) -- (3,3);
	\fill[valley] (0,0) circle (3pt) node[below right, black, font=\small] {Valley};
	\fill[valley] (2,2) circle (3pt) node[below right, black, font=\small] {Valley};
	\node at (1.5,-0.6) {(c) \((\{1,3\}, \{1,3\})\)};
	\end{scope}
	
	\end{tikzpicture}
	\caption{{Dyck paths from $(0,0)$ to $(3,3)$ with exactly two valleys, counting the origin.}}
	\label{fig1}
\end{figure}
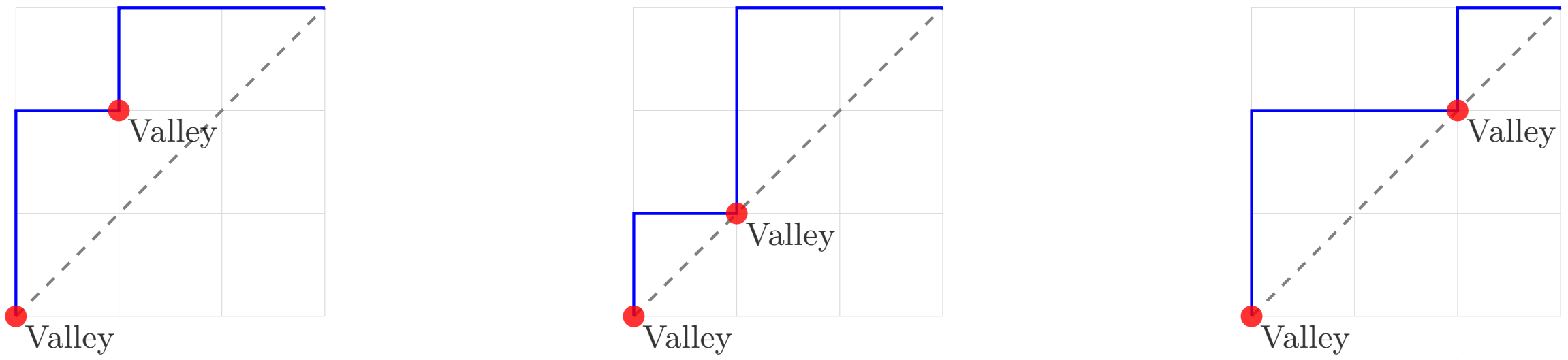

\begin{prop}\label{pbtm}
Let $m\geq 1$ be an integer, $t\in[m]$, and $\btm$ be the set of all pairs
\[
(\{i_1<\cdots<i_t\},\{j_1<\cdots<j_t\})
\]
such that $i_r,j_r\in[m]$ and $i_r\leq j_r$ for every $r\in[t]$, $i_1=1$, and $j_1=1$. Then
\[
|\btm|=N(m,t)=\frac{1}{m}\binom{m}{t} \binom{m}{t-1}.
\]
\end{prop}
\begin{proof}
Define a map from $\btm$ to the set of Dyck paths with $t$ valleys by
\[
(\{1=i_1<\cdots<i_t\leq m\}, \{1=j_1<\cdots<j_t\leq m\})\to L,
\]
where $L$ is the unique Dyck path whose valleys are
$(0,0),(i_2-1,j_2-1),\ldots,(i_t-1,j_t-1)$.
The inequalities $i_r\leq j_r$ ensure that the path
stays on or above the diagonal.
Conversely, listing the valleys in their order along
the path and adding $1$ to each coordinate recovers
the two index sets. Thus, the map is a bijection.
The result follows from the Narayana formula. (See Figure~\ref{fig1}.)
\end{proof}

\subsection{Formula for calculating $|\ceid|$}\label{sidel3}
In this subsection, we give a formula for the number of
idempotents in $C(\vep)$, where $\vep\in E(\cn)$.

\begin{lemma}\label{lceide}
Let $\vep\in E(\cn)$ be written as in {\rm Notation~\ref{nide}}, and fix $\rho\in\ceid$.
Set $p_i=|E_i|$ for every $i\in[m]$. For $i,j\in[m]$ with $e_j\rho=e_i$ and an integer
$s$ with $0\leq s\leq p_i-1$, let $A_{i,j}^s$ be the set of restrictions $\xi|_{E_j}$, where
$\xi\in\ceid$ satisfies $\xi|_{\fix(\vep)}=\rho|_{\fix(\vep)}$ and $|(E_i\setminus\{e_i\})\cap\fix(\xi)|=s$.
When $e_i\rho=e_i$, write $A_i^s=A_{i,i}^s$.
Then:
\begin{itemize}
\item [\rm(1)] If $j>i$, then $|A_{i,j}^0|=1$
and $|A_{i,j}^s|=0$ for $s>0$.

\item [\rm(2)] If $i<m$ and $e_{i+1}\rho=e_i$,
then $|A_i^0|=1$ and $|A_i^s|=0$ for $s>0$.

\item [\rm(3)] If $e_i\rho=e_i$ and either $i=m$
or $i<m$ and $e_{i+1}\rho=e_{i+1}$, then
$|A_i^s|=\binom{p_i-1}{s}$.
\end{itemize}
\end{lemma}
\begin{proof} 
Suppose $j>i$, and let $\xi\in\ceid$ agree with $\rho$ on $\fix(\vep)$. Then $e_j\xi=e_i$, and so
$e_i\xi=(e_j\xi)\xi=e_j(\xi\xi)=e_j\xi=e_i$. Since $\xi$ is order-preserving,
$e_i=e_i\xi\leq y\xi\leq e_j\xi=e_i$ for every $y\in E_i$, and so  $E_i\xi=\{e_i\}$.
By Theorem~\ref{tcide}, $E_j\xi\subseteq E_i$. Since every point in $\im(\xi)$ is fixed by $\xi$, $E_i\xi=\{e_i\}$ implies
$E_j\xi=\{e_i\}$. Thus, the only possible restriction of $\xi\in\ceid$ to $E_j$ is the constant map with value $e_i$.
This proves~(1).

For every $x\in E_i$, we have $e_i=e_i\xi\leq x\xi\leq e_{i+1}\xi=e_i$. Hence, $\xi|_{E_i}$ is the constant map with value~$e_i$,  which proves~(2).

Suppose the hypotheses of~(3) hold. For each  subset $S$ of $E_i\setminus\{e_i\}$ with $|S|=s$, there is a unique idempotent
$\zeta\in\cn(E_i,E_i)$ with $\fix(\zeta)=\{e_i\}\cup S$. Replace $\rho|_{E_i}$ by $\zeta$, leaving $\rho$
unchanged outside $E_i$. Then, since $e_{i+1}\rho=e_{i+1}$ if $i<m$, $\zeta$ is an idempotent in $\cn$ that commutes with $\vep$.
Every such choice of~$\zeta$ gives an element of $A_i^s$, and every element of $A_i^s$ arises uniquely in this way.
There are $\binom{p_i-1}{s}$ choices of $\zeta$, proving~(3).
\end{proof}

\begin{theorem}\label{tceide}
Let $\vep\in E(\cn)$ be written as in {\rm Notation~\ref{nide}}. Then,
\[
|\ceid|=\sum_{t=1}^{m}\sum_{1=i_1<\cdots<i_t\leq m}2^{p_{i_{r_1}}+\cdots+p_{i_{r_u}}-u},
\]
where $i_1=1$, $p_{i_r}=|E_{i_r}|$ for every $r\in[t]$, and $\{r_1<\cdots<r_u\}=\{r\in[t]:i_{r+1}-i_r=1\}$, where $i_{t+1}=m+1$.
(The sum in the exponent is understood to be $0$ when $u=0$.)
\end{theorem}
\begin{proof} 
Fix $1=i_1<\cdots<i_t\leq m$. Select $\rho\in C(\vep)$ defined by $x\rho=e_{i_r}$ for every $x\in E_j$ with 
$r\in[t]$ and $\ i_r\leq j<i_{r+1}$,
and note that $\fix(\rho)=\{e_{i_1}<\cdots<e_{i_t}\}$. 
Let $A$ be the set of all $\xi\in\ceid$ such that $\fix(\vep)\cap\fix(\xi)=\{e_{i_1}<\cdots<e_{i_t}\}$.
If $u=0$, parts~(1) and~(2) of Lemma~\ref{lceide} force $\xi=\rho$ for every $\xi\in A$.
Hence, $A=\{\rho\}$ and $|A|=1$.

Now suppose $u\geq1$. For each $j\in[u]$, let $A_{i_{r_j}}=\{\xi|_{E_{i_{r_j}}}:\xi\in A\}$.
Put $q=i_{r_j}$. If $q<m$, then
$i_{r_j+1}=q+1$, so both $e_q$ and $e_{q+1}$
are fixed by $\xi$.
If $q=m$, no condition on a subsequent root is needed.
In either case, part~(3) of Lemma~\ref{lceide} applies. Hence,
\[
|A_{i_{r_{\!j}}}|=\sum_{s=0}^{p_{i_{r_{\!j}}}-1}|A_{i_{r_{\!j}}}^s|=\sum_{s=0}^{p_{i_{r_{\!j}}}-1}\binom{p_{i_{r_{\!j}}}-1}{s}=2^{p_{i_{r_{\!j}}}-1}.
\]
By (1) and (2) of Lemma~\ref{lceide},
the function $f:A\to A_{i_{r_1}}\times\cdots\times A_{i_{r_u}}$ defined by $\xi f=(\xi|_{E_{i_{r_1}}},\ldots,\xi|_{E_{i_{r_u}}})$ is a bijection. Thus,
\begin{equation}\label{eceide}
|A|=|A_{i_{r_1}}||A_{i_{r_2}}|\cdots|A_{i_{r_u}}|=2^{p_{i_{r_1}}-1}2^{p_{i_{r_2}}-1}\cdots2^{p_{i_{r_u}}-1}=2^{p_{i_{r_1}}+\cdots+p_{i_{r_u}}-u}.
\end{equation}

We also have $|A|=2^{p_{i_{r_1}}+\cdots+p_{i_{r_u}}-u}$
when $u=0$ since then both sides are equal to $1$.
Summing the rightmost expression in \eqref{eceide} over all chains $1=i_1<\cdots<i_t\leq m$, we obtain the cardinality of the set of $\xi\in\ceid$ 
such that $|\fix(\vep)\cap\fix(\xi)|=t$. Finally, 
summing over all $t\in[m]$, we obtain
\[
|\ceid|=\sum_{t=1}^{m}\sum_{1=i_1<\cdots<i_t\leq m}2^{p_{i_{r_1}}+\cdots+p_{i_{r_u}}-u},
\]
as required.
\end{proof}

\subsection{Formula for calculating $|\cnl(\vep)|$}\label{sidel4}
In this subsection, we give a formula for the number of
nilpotents in $C(\vep)$, where $\vep\in E(\cn)$. Unlike the set $\ceid$, which  need not be closed under composition,
$\cenil$ is a subsemigroup of $C(\vep)$. Indeed, the zero of $\cn$ is the constant transformation $\theta$ with image $\{1\}$. 
Recall that $\al\in\cn$ is nilpotent if and only if
$\fix(\al)=\{1\}$, equivalently, $x\al<x$ for every $x>1$. If $\al,\bt\in\ncn$, then $x(\al\bt)\leq x\al<x$ for every $x>1$.
Thus, $\ncn$ is a subsemigroup of $\cn$.

We will need the following subset $\dtm$ of $\btm$ (see Notation~\ref{nbtm}).

\begin{nota}\label{ndtm}
Let $m\geq1$ be an integer and $t\in[m]$. We denote by $\dtm$ the set of all pairs
\[
(\{i_1<\cdots<i_t\},\{j_1<\cdots<j_t\})
\]
such that $i_r,j_r\in[m]$, $i_1=j_1=1$, $i_t<m$, and $i_r<j_r$ for every $2\leq r\leq t$.
\end{nota}

\begin{theorem}\label{tcenil}
Let $\vep\in E(\cn)$ be written as in {\rm Notation~\ref{nide}}. Then,
\begin{itemize} 
  \item [\rm(1)] if $m\geq2$, then 
\[
|\cenil| =\sum_{t=1}^{m-1} \sum_{(\{i_1<\cdots<i_t\},\{j_1<\cdots<j_t\})\in\dtm}\prod_{r=1}^{t}\binom{p_r+q_r-2}{p_r-1},
\]
where $\dtm$ is defined in {\rm Notation~\ref{ndtm}}, $p_r=|E_{i_r}|$, $q_r=|E_{j_{r+1}-1}|$, and $j_{t+1}=m+1$;
  \item [\rm(2)] if $m=1$, then $|\cenil|=\frac1{n}\binom{2n-2}{n-1}$.
\end{itemize}
\end{theorem}
\begin{proof}
 For (1), follow the proof of Theorem~\ref{tcidecard} with the following modifications. First, sum over all pairs in $\dtm$, not $\btm$,
since if $\al\in N(\cn)$, then $\fix(\al)=\{1\}$, so $i_r<j_r$ for every $2\leq r\leq t$. 
Moreover, $e_m\al=e_{i_t}<e_m$, so $i_t<m$. Consequently, $t\leq m-1$.

Conversely, suppose that the pair associated with $\al\in C(\vep)$ belongs to $\dtm$. Then $e_j\al<e_j$ for every $2\leq j\leq m$.
In particular, $e_2\al=1$, so order preservation gives $x\al=1$ for every $x\in E_1$. If $x\in E_j$ with $j\geq2$, then
$(x\al)\vep=(x\vep)\al=e_j\al<e_j$. Hence, $x\al$ belongs to an interval preceding $E_j$, and therefore $x\al<x$.
Thus, $\fix(\al)=\{1\}$, proving that $\al$ is nilpotent.

Second, sum over $t\in[m-1]$, not $[m]$, since, by the definition of $\dtm$, we have $t<m$ if $m\geq2$. 
For $r<t$, we have $j_{r+1}>i_{r+1}\geq i_r+1$. For $r=t$, we have $j_{t+1}=m+1>i_t+1$. Thus, for every $r\in[t]$,
\[
K^{j_{r+1}-i_r}_{p_r,q_r}
=\binom{p_r+q_r-2}{p_r-1}
\]
by~\eqref{ekjpq}, and statement~(1) follows.

To prove (2), suppose $m=1$. Then, $\vep$ is the zero in $\cn$, and so $|\cenil|=|N(\cn)|=\frac1{n}\binom{2n-2}{n-1}$
by \cite[Prop.~2.3]{LaUm2006}.
\end{proof}

Let $m\geq2$ and $t\in[m-1]$. Then, there is a natural bijection from $B^{m-1}_t$ to $\dtm$:
\[
(\{i_1<i_2<\cdots<i_t\}, \{j_1<j_2<\cdots<j_t\})\to(\{i_1<i_2<\cdots<i_t\}, \{j_1<j_2+1<\cdots<j_t+1\}).
\]
Hence, $|\dtm|=|B_t^{m-1}|=N(m-1,t)=\frac1{m-1}\binom{m-1}{t}\binom{m-1}{t-1}$ (see Prop.~\ref{pbtm}).

\begin{example}\label{excidecard}
We will apply Theorems~\ref{tcidecard}, ~\ref{tceide}, and~\ref{tcenil} to an idempotent $\vep\in\mathcal{C}_{12}$ with $\fix(\vep)=\{1,5,8\}$, that is,
\[
\vep=\begin{pmatrix}\{1,2,3,4\}&\{5,6,7\}& \{8,9,10,11,12\}\\ 1 & 5 & 8\end{pmatrix}
\]
(see Figure~\ref{fig3}).
Then $m=3$, $|E_1|=4$, $|E_2|=3$, and $|E_3|=5$.

We first calculate $|C(\vep)|$. Suppose $t=1$. Then, $B_1^3=\{(\{1\},\{1\})\}$
and $j_2=m+1=4$, $p_1 = |E_1| = 4$, and $q_1 = |E_{j_2-1}| = |E_{4-1}| = |E_3| = 5$. 
We have $j_2-i_1=4-1=3$. Thus, $K_{4,5}^3 = \binom{4+5-2}{4-1} = \binom{7}{3} = 35$. The contribution for $t=1$ is $35$.

Suppose $t=2$. Then, $B_2^3=\{(\{1<2\},\{1<2\}),(\{1<2\},\{1<3\}),(\{1<3\},\{1<3\})\}$ and $j_3=m+1=4$.
Consider the pair $(\{1<2\}, \{1<2\})$. If $r=1$, then $p_1 = |E_1| = 4$, $q_1 = |E_{2-1}| = |E_1| = 4$, and $j_2-i_1= 2-1 = 1$. 
So, $K_{4,4}^1 = \frac{1}{4+1}\binom{2 \cdot 4}{4} = \frac{1}{5} \cdot 70 = 14$. 
If $r=2$, then $p_2 = |E_2| = 3$, $q_2 = |E_{4-1}| = |E_3| = 5$, and $j_3-i_2 = 4-2 = 2$. So, $K_{3,5}^2 = \binom{3+5-2}{3-1} = \binom{6}{2} = 15$. 
The product for this pair is $14 \cdot 15 = 210$. 
By similar calculations, the product for $(\{1<2\}, \{1<3\})$ is $150$, and the product for $(\{1<3\}, \{1<3\})$ is $420$. 
The contribution to the sum for $t=2$ is $210+150+420=780$.

 Suppose $t=3$. Then, $B_3^3 = \{(\{1<2<3\}, \{1<2<3\})\}$ and $j_4=m+1=4$. If $r=1$, then $K_{4,4}^1 = 14$.
 If $r=2$, then $K_{3,3}^1 = \frac{1}{3+1}\binom{2 \cdot 3}{3} = \frac{1}{4} \cdot 20 = 5$. If $r=3$, then $K_{5,5}^1 = 42$. The product for this pair is $14 \cdot 5 \cdot 42 = 2940$. The contribution for $t=3$ is $2940$.

Summing all contributions, we obtain $|C(\vep)| = 35 + 780 + 2940 = 3755$.

We now calculate $|\ceid|$. If $t=1$, then the set $\{r_1, \dots, r_u\}$ is empty ($u=0$), so the contribution to the sum is $2^{0-0} = 2^0 = 1$. 
Suppose $t=2$. The chains of $[3]$ of size $t=2$ with $i_1=1$ are $\{1<2\}$ and $\{1<3\}$. Also, $i_3=m+1=4$.
Consider the chain $\{i_1<i_2\} = \{1<2\}$. Then, $i_2 - i_1 = 2-1 = 1$
and $i_3-i_2 = 4-2=2$. So, we have $\{r_1\} = \{1\}$, which means $u=1$.
Thus, the contribution of $\{1<2\}$ is $2^{p_{i_{r_1}} - u} = 2^{p_1 - 1} = 2^{4-1} = 2^3 = 8$.
Similarly, we check that the contribution of $\{1<3\}$ is $16$.
The total for $t=2$ is $8 + 16 = 24$. 

Suppose $t=3$.
The only chain of $[3]$ of size $t=3$ with $i_1=1$ is $\{1<2<3\}$ and $i_4=m+1=4$.
In this case, $u=3$, and $\{1<2<3\}$ contributes $512$ to the count of $|\ceid|$.

Summing the contributions of all cases, we obtain $|\ceid| = 1 + 24 + 512 = 537$.

Finally, we calculate $|\cenil|$. Suppose $t=1$. Then, $D_1^3=\{(\{1\},\{1\})\}$, and so, as in the calculation of $|C(\vep)|$,
the contribution for $t=1$ is $35$.

Suppose $t=2$. Then, $D_2^3=\{(\{1<2\},\{1<3\})\}$ and $j_3=m+1=4$.
If $r=1$, then $p_1 = |E_1| = 4$, $q_1 = |E_{3-1}| = |E_2| = 3$, and $j_2-i_1 = 3-1 = 2$. So, $K_{4,3}^2=\binom{4+3-2}{4-1}=\binom{5}{3}=10$. 
If $r=2$, then $p_2 = |E_2| = 3$, $q_2 = |E_{4-1}| = |E_3| = 5$, and so $\binom{3+5-2}{3-1} = \binom{6}{2} = 15$. The product for this pair,
and the contribution for $t=2$, is $10 \cdot 15 = 150$. 

Note that $t=3$ is impossible. Summing the contributions of all cases, we obtain $|\cenil| =35+150 =185$.

\end{example}

\section{Centralizers of nilpotents of $\cn$}\label{snil}
\setcounter{equation}{0}
Recall that the monoid $\cn$ has a zero, namely the constant transformation $\theta$ with $\im(\theta)=\{1\}$,
and an identity, that is the transformation $\vep_n$ that fixes every element of $\xn$.
It is easy to see that $\mu\in\cn$ is a nilpotent if and only if $\fix(\mu)=\{1\}$, and consequently the subset $N(\cn)$ of nilpotents
of $\cn$ is a subsemigroup of $\cn$.

In this section, we study the centralizers of nilpotents in $\cn$. For a given $\mu\in N(\cn)$,
we characterize the elements of $C(\mu)$ and present algorithms that construct
$C(\mu)$, $\cid(\mu)=C(\mu)\cap E(\cn)$ (the set of idempotents of $C(\mu)$), and
 $\cnl(\mu)=C(\mu)\cap N(\cn)$ (the semigroup of nilpotents of $C(\mu)$). We also describe the nilpotents $\mu\in\cn$
for which the centralizer $C(\mu)$ consists only of idempotents, of nilpotents together with the identity, 
or of nilpotents and idempotents; that is,
for which $C(\mu)=\cid(\mu)$, $C(\mu)=\cnl(\mu)\cup\{\vep_n\}$, and $C(\mu)=\cnl(\mu)\cup\cid(\mu)$.

Proposition~\ref{paincn} implies the following result.

\begin{cor}\label{cnil}
Let $\mu\in\cn$. Then, $\mu$ is a nilpotent if and only if $\ga(\mu)$ is a level-ordered tree with root~$1$, and 
for every $x\in\xn$, the inverse image $x\mu^{-1}$ is either empty or an interval of $\xn$.
\end{cor}

The leftmost tree in Figure~\ref{fig2} represents a nilpotent $\mu\in\mathcal{C}_9$.

\subsection{Elements of $C(\mu)$}\label{snil1}
In this subsection, we characterize the elements of $C(\mu)=\{\al\in\cn:\al\mu=\mu\al\}$ for an arbitrary nilpotent $\mu\in\cn$. 
The following result was proved in \cite[Prop.~3.1]{ArKo2013}. 

\begin{lemma}\label{lend}
Let $\al,\bt\in\tn$. Then $\al\in C(\bt)$ if and only if $\al$ is a homomorphism from $\ga(\bt)$ to $\ga(\bt)$.
\end{lemma}

We also need to extend $\cn(J,I)$ (see Notation~\ref{ncnij}).

\begin{nota}\label{ncnstij}
Let $I=\{i<i+1<\cdots<i+p-1\}$ and $J=\{j<j+1<\cdots<j+q-1\}$ be intervals of $\xn$, where $I\leq J$.
We  denote by $\cnst(J,I)$ the set 
\[
\cnst(J,I)=\{\al:J\to I:\mbox{$\al$ is decreasing and order-preserving}\}.
\]
\end{nota}

The difference between $\cn(J,I)$ and $\cnst(J,I)$ is that in the latter, we do not require $j\al=i$.
The following lemma follows from the fact that $\cnst(I,I)=\cn(I,I)\cong\mathcal{C}_p$.
\begin{lemma}\label{lcnstii}
Let $I=\{i<i+1<\cdots<i+p-1\}$ be an interval of $\xn$. Then, $\lvert\cnst(I,I)\rvert=C_p=\frac1{p+1}\binom{2p}{p}$.
\end{lemma}

\begin{lemma}\label{lcnstij}
Let $I=\{i<i+1<\cdots<i+p-1\}$ and $J=\{j<j+1<\cdots<j+q-1\}$ be intervals of $\xn$, where $I<J$. Then, $\lvert\cnst(J,I)\rvert=\binom{p+q-1}{p-1}$.
\end{lemma}
\begin{proof} 
Let $t=\min\{p,q\}$ and fix $m\in[t]$. Then, every $\al\in\cnst(J,I)$ of rank $m$ can be written as
\[
\al=\begin{pmatrix}A_1  & \ldots & A_m \\ a_1  & \ldots & a_m\end{pmatrix},
\]
where $\{a_1<\cdots<a_m\}$ is a chain of $I$ and $\lan A_1,\ldots,A_m\ran$ is an ordered convex partition of $J$. We can choose
$\{a_1<\cdots<a_m\}$ in $\binom{p}{m}$ ways, and $\lan A_1,\ldots,A_m\ran$ in $\binom{q-1}{m-1}$ ways.
Thus, the number of $\al\in\cnst(J,I)$ of rank $m$ is $\binom{p}{m}\binom{q-1}{m-1}$. Hence,
\[
\lvert\cnst(J,I)\rvert=\sum_{m=1}^{t}\binom{p}{m}\binom{q-1}{m-1}=\sum_{m=1}^{t}\binom{p}{(p-1)-(m-1)}\binom{q-1}{m-1}=\binom{p+q-1}{p-1},
\]
where the last equality is a version of Vandermonde's identity \cite[Problem  26, p.\ 167]{St1993}.
\end{proof}

\begin{theorem}\label{tcnilg}
Let $\mu\in \ncn$, $\tr=\ga(\mu)$, and $\al\in\cn$.
Then, $\al\in C(\mu)$ if and only if the following conditions are satisfied:

\begin{itemize} 
    \item [\rm(1)] $(\im(\mu))\al\subseteq \im(\mu)$;
    \item [\rm(2)] for all $x,y\in\im(\mu)$, if $x=y\al$, then $\al|_{J}\in\cnst(J,I)$, where $I=x\mu^{-1}$ and $J=y\mu^{-1}$ (see {\rm Notation~\ref{ncnstij}}).
\end{itemize} 
\end{theorem}
\begin{proof}
Let $\al\in C(\mu)$. Suppose $y \in \im(\mu)$. Then, $y=x\mu$ for some $x\in\xn$. Thus, $y\al=(x\mu)\al=(x\al)\mu\in\im(\mu)$,
and so~(1) holds.

To prove~(2), let $x,y\in\im(\mu)$ with $x=y\al$. Suppose $z\in J$. Then, $z\to y$ in $\tr$, and so $z\al\to y\al=x$
by Lemma~\ref{lend}. Hence, $(z\al)\mu=x$, and so $z\al\in I$. We have proved that
$\al$ maps $J$ to $I$, that is, $\al|_J: J \to I$. Since $\al$ is decreasing and order-preserving, $\al|_J \in \cnst(J,I)$.
	
Conversely, suppose that conditions (1) and (2) hold. Let $z \in\xn$. Set $y=z\mu$ and $x=y\al$. Then, $x\in\im(\mu)$ by (1), and so
$z\al\in x\mu^{-1}$ by (2). Hence, $z(\al\mu)=x=y\al=z(\mu\al)$,
which shows that $\al\in C(\mu)$.
\end{proof}

\subsection{Algorithm for constructing $C(\mu)$}\label{snil2}

Let $\mu\in N(\cn)$. It does not seem feasible to find a formula for $|C(\mu)|$ of the kind we have in Theorem~\ref{tcidecard}
for  $|C(\vep)|$. Instead, we will present an algorithm that constructs the centralizer $C(\mu)$, using the representation
of $\mu$ as $\ga(\mu)$ (Algorithm~\ref{acentmu}). We first
need two preliminary algorithms.

The first algorithm, given an interval $I=\{i,i+1,\ldots,i+p-1\}$ of $\xn$
as input (see Definition~\ref{dint}), returns the set $\cnst(I,I)$ as output (see Notation~\ref{ncnstij}).
(Recall that $\cnst(I,I)=\cn(I,I)$.) The algorithm successively constructs $\cn^1(I,I),\ldots,\cn^p(I,I)$,
where $\cn^m(I,I)$ is the set of elements of $\cnst(I,I)$ of rank $m$.

\begin{alg}\label{acnii}
	Let $I=\{i,i+1,\ldots,i+p-1\}$ be an interval of $\xn$.
	\begin{enumerate} 
		\item Set $m=1$, $\cn^m(I,I)=\emptyset$, and go to step 2.
		\item For all chains $\{k_1<\cdots<k_m\}$ and $\{l_1<\cdots<l_m<l_{m+1}\}$, where $k_r,l_r\in[p]$ and $k_r\leq l_r$ for every $2\leq r\leq m$,
		$k_1=l_1=1$, and $l_{m+1}=p+1$, do the following:
		\begin{itemize}
			\item [(i)]	set $a_r=i+k_r-1$ and $A_r=\{i+l_r-1,\ldots,i+l_{r+1}-2\}$ for every $r\in[m]$;
			\item [(ii)] add $\al=\begin{pmatrix} A_1&\ldots&A_m\\a_1&\ldots&a_m\end{pmatrix}$ to $\cn^m(I,I)$.
		\end{itemize}
		Go to step 3.
		\item Increase $m$ by $1$.
		If $m\leq p$, set $\cn^m(I,I)=\emptyset$, and go to step 2. If $m=p+1$, RETURN $\cnst(I,I)=\cn^1(I,I)\cup\ldots\cup\cn^p(I,I)$ and STOP.
	\end{enumerate}
\end{alg}

The second algorithm, given intervals $I=\{i,i+1,\ldots,i+p-1\}$ and $J=\{j,j+1,\ldots,j+q-1\}$ of $\xn$ such that $I<J$
as input (see Definition~\ref{drelint}), returns the set $\cnst(J,I)$ as output (see Notation~\ref{ncnstij}).
The algorithm successively constructs $\cn^1(J,I),\ldots,\cn^t(J,I)$,
where $t=\min\{p,q\}$, and $\cn^m(J,I)$ is the set of elements of $\cnst(J,I)$ of rank $m$.

\begin{alg}\label{acnij}
	Let $I=\{i,i+1,\ldots,i+p-1\}$ and $J=\{j,j+1,\ldots,j+q-1\}$ be intervals of $\xn$ such that $I<J$.
		\begin{enumerate} 
		\item Set $t=\min\{p,q\}$, $m=1$, $\cn^m(J,I)=\emptyset$, and go to step 2.
		\item For all chains $\{k_1<\cdots<k_m\}$ and $\{l_1<\cdots<l_m<l_{m+1}\}$, where $k_r\in[p]$ and $l_r\in[q]$ for every $r\in[m]$,
		$k_1\geq1$, $l_1=1$, and $l_{m+1}=q+1$, do the following:
		\begin{itemize}
			\item [(i)]	set $a_r=i+k_r-1$ and $A_r=\{j+l_r-1,\ldots,j+l_{r+1}-2\}$ for every $r\in[m]$;
			\item [(ii)] add $\al=\begin{pmatrix} A_1&\ldots&A_m\\a_1&\ldots&a_m\end{pmatrix}$ to $\cn^m(J,I)$.
		\end{itemize}
		Go to step 3.
		\item Increase $m$ by $1$.
		If $m\leq t$, set $\cn^m(J,I)=\emptyset$, and go to step 2. If $m=t+1$, RETURN $\cnst(J,I)=\cn^1(J,I)\cup\ldots\cup\cn^t(J,I)$ and STOP.
		\end{enumerate}
\end{alg}

Algorithm~\ref{acnij} is justified by the proof of Lemma~\ref{lcnstij}.

\begin{nota}\label{nmuk}
Let $\mu\in N(\cn)$, $\tr=\ga(\mu)$, $h=\hei(\tr)$, and $0\leq k\leq h$. We denote by $\mu_k$ the restriction of $\mu$ to
$A_k=V_0(\tr)\cup\ldots\cup V_k(\tr)$. Let $A\subseteq\xn$. We denote by $\cn(A)$ the set of all order-preserving and decreasing transformations from $A$ to $A$. Note that $\mu_k\in\cn(A_k)$ and $\mu_h=\mu$.
\end{nota}

\begin{defi}\label{djoin}
	Let $\al:A\to\xn$ and $\bt:B\to\xn$ be functions such that $A$ and $B$ are nonempty subsets of $\xn$ with $A\cap B=\emptyset$.
	The \emph{join} of $\al$ and $\bt$ is the function $\al\joi\bt:A\cup B\to\xn$ defined by
\[
x(\al\joi\bt)= \left\{
	\begin{array}{ll}
		x\al & \mbox{if $x\in A$,}\\
		x\bt & \mbox{if $x\in B$.}
\end{array}
\right.
\]
\end{defi}

The following algorithm, given a tree $\tr=\ga(\mu)$, where $\mu\in N(\cn)$, as input (see Corollary~\ref{cnil}),
returns the centralizer $C(\mu)$ as output. It constructs inductively the centralizers 
\[
C(\mu_0),C(\mu_1),\ldots,C(\mu_h)=C(\mu),
\]
where $h=\hei(\tr)$ (see Notation~\ref{nmuk}).

\begin{alg}\label{acentmu}
	Let $\mu\in N(\cn)$, $\tr=\ga(\mu)$, and $h=\hei(\tr)$.
	\begin{enumerate} 
		\item Set $C(\mu_0)=\{\al_1\}$, where $1\al_1=1$.
		If $h=0$, RETURN $C(\mu)=C(\mu_0)$ and STOP. If $h\geq1$, set $k=0$ ($k$ will range from $0$ to $h$) and go to step 2.
		\item At this point, $k<h$, and we have calculated $C(\mu_k)=\{\al_1,\ldots,\al_t\}$.  Find the set $W_k=\{y_1,\ldots,y_s\}$ of all elements
		$y\in V_k(\tr)$ such that $y\mu^{-1}\ne\emptyset$. Set $C(\mu_{k+1})=\emptyset$, $i=1$ ($i$ will range from $1$ to $t$), and go to step 3.
		\item For every $j\in [s]$, set $x_j=y_j\al_i$. If $x_j\mu^{-1}=\emptyset$ for some $j\in[s]$, go to step~5.
		Otherwise, go to step~4.
		\item Calculate $\cnst(y_j\mu^{-1},x_j\mu^{-1})=\{\bt_1^j,\ldots,\bt_{v_j}^j\}$ for each $j\in[s]$ (see Algorithms~\ref{acnii} and~\ref{acnij}).
		For every $(j_1,\ldots,j_s)\in[v_1]\times\ldots\times[v_s]$
		add $\al_i\joi\bt^1_{j_1}\joi\ldots\joi\bt^s_{j_s}$ to $C(\mu_{k+1})$ if it is order-preserving, and go to step 5.
		\item Increase $i$ by 1. If $i\leq t$, go to step 3. If $i=t+1$, go to step 6.
		\item Increase $k$ by $1$. If $k<h$, go to step 2. If $k=h$, RETURN $C(\mu)=C(\mu_h)$ and STOP.
	\end{enumerate}
\end{alg}
Using Theorem~\ref{tcnilg} and induction on $k$, it is straightforward to show that Algorithm~\ref{acentmu} constructs 
all elements of $C(\mu)$ (without repetitions). 

\begin{rem}\label{rcentmu}
In step 4 of Algorithm~\ref{acentmu}, $\al_i\joi\bt^1_{j_1}\joi\ldots\joi\bt^s_{j_s}$,
which maps $A_{k+1}$ to $A_{k+1}$, is decreasing by the construction
of $\al_i\joi\bt^1_{j_1}\joi\ldots\joi\bt^s_{j_s}$, so this does not have to be checked. Checking if $\al_i\joi\bt^1_{j_1}\joi\ldots\joi\bt^s_{j_s}$ is order-preserving
is straightforward. Suppose $A_{k+1}=\{1,\ldots,m\}$. Then, 
$\al_i\joi\bt^1_{j_1}\joi\ldots\joi\bt^s_{j_s}=\begin{pmatrix}1&\ldots&m\\a_1&\ldots&a_m\end{pmatrix}$
is order-preserving if and only if $a_1\leq\ldots\leq a_m$.
\end{rem}

\begin{figure}[htbp]
	\centering
	\begin{tikzpicture}[
	level distance=1.0cm,
	level 1/.style={sibling distance=3.0cm},
	level 2/.style={sibling distance=1.6cm},
	level 3/.style={sibling distance=1.2cm},
	every node/.style={circle, fill=black, inner sep=1.2pt}
	]
	\node[label=above:1] {}
	child { node[label=left:2] {}
		child { node[label=below:4] {} }
		child { node[label=right:5] {}
			child { node[label=below:7] {} }
		}
	}
	child { node[label=right:3] {}
		child { node[label=right:6] {}
			child { node[label=below:8] {} }
			child { node[label=below:9] {} }
		}
	};
	\end{tikzpicture}
	\caption{The tree $\ga(\mu)$ for the nilpotent $\mu$ in Example~\ref{ecentmu}.}
	\label{fcentmu}
\end{figure}

\begin{example}\label{ecentmu}
	Consider  $\mu=\left( \begin{smallmatrix}
	1 & 2 & 3 & 4 & 5 & 6 & 7 & 8 & 9 \\
	1 & 1 & 1 & 2 & 2 & 3 & 5 & 6 & 6 
	\end{smallmatrix}\right) \in N(\mathcal{C}_{9})$. The tree
$\tr=\ga(\mu)$ is presented in Figure~\ref{fcentmu}. We will construct $C(\mu)$ using Algorithm~\ref{acentmu}. 
	The tree $\tr$ has height $h=3$, and its levels are $V_0(\tr) = \{1\}$, $V_1(\tr) = \{2, 3\}$, $V_2(\tr) = \{4, 5, 6\}$, and $V_3(\tr) = \{7, 8, 9\}$.
\begin{itemize}
	\item[(1)] $C(\mu_0)=\{\al_1\}$, where $1\al_1=1$, by step 1 of Algorithm~\ref{acentmu}.
	
	\item[(2)]  $C(\mu_1)=\cnine^*(1\mu^{-1},1\mu^{-1})=\{\al_1=111,\al_2=112,\al_3=113,\al_4=122,\al_5=123\}$ (see steps 2--5 of Algorithm~\ref{acentmu}).
	
	\item[(3)] We will present the construction of $C(\mu_2)$ in some detail. At the beginning, $k=1$, $C(\mu_1)$ is as in (2), $t=|C(\mu_1)|=5$, $C(\mu_2)=\emptyset$,
	 and $W_1=\{y_1=2,y_2=3\}$ (see step 2 of Algorithm~\ref{acentmu}). We have to go through $\al_1,\ldots,\al_5$ and either reject $\al_i$ completely or add some of its extensions
	 to $C(\mu_2)$ (see steps 3 and 4 of Algorithm~\ref{acentmu}).
	 
	 Consider $\al_1=111\in C(\mu_1)$. Then, $x_j=y_j\al_1=1$ for $j=1,2$. We calculate 

\begin{align*}
\cnine^*(y_1\mu^{-1},x_1\mu^{-1})&=\cnine^*(2\mu^{-1},1\mu^{-1})=\cnine^*(\{4<5\},\{1<2<3\})\\
&=\mbox{{\LARGE $\{$}}\bt^1_1=\begin{pmatrix} 4 & 5\\1 & 1\end{pmatrix},\bt^1_2=\begin{pmatrix} 4 & 5\\1 & 2\end{pmatrix},
\bt^1_3=\begin{pmatrix} 4 & 5\\1 & 3\end{pmatrix},\bt^1_4=\begin{pmatrix} 4 & 5\\2 & 2\end{pmatrix},
\bt^1_5=\begin{pmatrix} 4 & 5\\2 & 3\end{pmatrix},\\
&\hskip 7.5mm\bt^1_6=\begin{pmatrix} 4 & 5\\3 & 3\end{pmatrix}\mbox{{\LARGE $\}$}},\\
\cnine^*(y_2\mu^{-1},x_2\mu^{-1})&=\cnine^*(3\mu^{-1},1\mu^{-1})=\cnine^*(\{6\},\{1<2<3\})\\
&=\mbox{{\LARGE $\{$}}\bt^2_1=\begin{pmatrix} 6\\1 \end{pmatrix},\bt^2_2=\begin{pmatrix} 6\\2 \end{pmatrix},\bt^2_3=\begin{pmatrix} 6\\3\end{pmatrix}\mbox{{\LARGE $\}$}}.
\end{align*}

There are $18$ choices for $(j_1,j_2)\in[6]\times[3]$, so we have $18$ extensions of $\al_1$:
$\al_1\joi\bt^1_1\joi\bt^2_1=111111$, $\al_1\joi\bt^1_1\joi\bt^2_2=111112$, $\al_1\joi\bt^1_1\joi\bt^2_3=111113$, 
$\al_1\joi\bt^1_2\joi\bt^2_1=111121$, $\al_1\joi\bt^1_2\joi\bt^2_2=111122$, $\al_1\joi\bt^1_2\joi\bt^2_3=111123$,
$\al_1\joi\bt^1_3\joi\bt^2_1=111131$, $\al_1\joi\bt^1_3\joi\bt^2_2=111132$, $\al_1\joi\bt^1_3\joi\bt^2_3=111133$,
$\al_1\joi\bt^1_4\joi\bt^2_1=111221$, $\al_1\joi\bt^1_4\joi\bt^2_2=111222$, $\al_1\joi\bt^1_4\joi\bt^2_3=111223$,
$\al_1\joi\bt^1_5\joi\bt^2_1=111231$, $\al_1\joi\bt^1_5\joi\bt^2_2=111232$, $\al_1\joi\bt^1_5\joi\bt^2_3=111233$,
$\al_1\joi\bt^1_6\joi\bt^2_1=111331$, $\al_1\joi\bt^1_6\joi\bt^2_2=111332$, and $\al_1\joi\bt^1_6\joi\bt^2_3=111333$.
Rejecting eight that are not order-preserving, for example, $\al_1\joi\bt^1_2\joi\bt^2_1=111121$, we add the remaining ten to
$C(\mu_2)$: $111111$, $111112$, $111113$,  $111122$, $111123$, $111133$, $111222$, $111223$,  $111233$, and $111333$.

Repeating the above calculations for $\al_2,\al_3,\al_4,\al_5$, we obtain	
	\[
	\begin{aligned}
	C(\mu_2)=\{\,&
	111111,\,111112,\,111113,\,111122,\, 111123,\, 111133,\,111222,\\
	&111223,\,111233,\,111333,\,112224,\,112225,\,112234,\,112235,\\
	&112334,\,112335,\,113336,\, 122444,\,122445,\, 122455,\, 123446,\,123456\}.
	\end{aligned}
	\]
	
	\item[(4)] Finally, we calculate
\[
\begin{aligned}
C(\mu)=C(\mu_3)=\{&
111111111,\,111111112,\,111111113,\,111111122,\,111111123,\,111111133,\\
&111111222,\,111111223,\,111111233,\,111111333,\,111112244,\,111112245,\\
&111112255,\,111112344,\,111112345,\,111112355,\,111113366,\,111122444,\\
&111122445,\,111122455,\,111122555,\,111123466,\,111123566,\,111133666,\\
&111222444,\,111222445,\,111222455,\,111222555,\,111223466,\,111223566,\\
&111233666,\,111333666,\,112225577,\,112235677,\,112335677,\,113336688,\\
&113336689,\,122455777,\,123456788,\,123456789\},
\end{aligned}
\]
so $|C(\mu)|=40$.
\end{itemize}
\end{example}

\subsection{Algorithm for constructing $\cnl(\mu)$}\label{snil3}
Algorithm~\ref{acentmu} can be modified easily to construct $\cenilmu=C(\mu)\cap N(\cn)$.

\begin{nota}\label{ncnnilii}
	Let $I=\{i<i+1<\cdots<i+p-1\}$ and $J=\{j<j+1<\cdots<j+q-1\}$ be intervals of $\xn$ with $I\leq J$. 
	We  denote by $\cn^{nil}(I,I)$ the  following subset of $\cn(I,I)$ (see Notation~\ref{ncnij}): 
	\[
	\mathcal{C}_n^{nil}(I,I)=\{\al\in\cn(I,I): (i+r)\al\ne i+r\mbox{ for every $r\in\{1,\ldots,p-1$\}}\}, 
	\]
and extend this definition to $\mathcal{C}_n^{nil}(J,I)$ by
\[
\mathcal{C}_n^{nil}(J,I)=
\begin{cases}
\mathcal{C}_n^{nil}(I,I) & \text{if } I=J,\\[2mm]
\mathcal{C}_n^*(J,I) & \text{if } I<J.
\end{cases}
\]
\end{nota}

Recall that $\lvert\cn(I,I)\rvert=C_p=\frac1{p+1}\binom{2p}{p}$ (Lemma~\ref{lcnii}).

\begin{lemma}\label{lcnnilii}
	Let $I=\{i<i+1<\cdots<i+p-1\}$ be an interval of $\xn$. Then, $\lvert\cn^{nil}(I,I)\rvert=\frac{1}{p}\binom{2p-2}{p-1}$.
\end{lemma}
\begin{proof}
If $p=1$, then the result is clearly true. Suppose $p\geq2$.
Consider the interval  $I'=\{i<i+1<\cdots<i+p-2\}$ of $\xn$. There is a natural bijection from $\cn(I',I')$ to $\cn^{nil}(I,I)$:
\[
\begin{pmatrix} i & i+1 & i+2 & i+3 & \ldots & i+p-2 \\ i & a_{i+1} & a_{i+2} & a_{i+3} & \ldots & a_{i+p-2}\end{pmatrix}\rightarrow 
\begin{pmatrix} i & i+1 & i+2 & i+3 & \ldots & i+p-1 \\ i & i & a_{i+1} & a_{i+2} & \ldots & a_{i+p-2}\end{pmatrix}.
\]
Hence, $\lvert\cn^{nil}(I,I)\rvert=\lvert \cn(I',I')\rvert=C_{p-1}=\frac{1}{p}\binom{2p-2}{p-1}$.
\end{proof}

The following algorithm, given an interval $I$ of $\xn$
as input, returns the set $\cn^{nil}(I,I)$ as output.

\begin{alg}\label{acnnilii}
	Let $I=\{i,i+1,\ldots,i+p-1\}$ be an interval of $\xn$.
\begin{enumerate} 
	\item Set $m=1$, $\cn^m(I,I)=\emptyset$, and go to step 2.
	\item For all chains $\{k_1<\cdots<k_m\}$ and $\{l_1<\cdots<l_m<l_{m+1}\}$, where $k_r,l_r\in[p]$ and $k_r<l_r$ for every $2\leq r\leq m$,
	$k_1=l_1=1$, and $l_{m+1}=p+1$, do the following:
	\begin{itemize}
		\item [(i)]	set $a_r=i+k_r-1$ and $A_r=\{i+l_r-1,\ldots,i+l_{r+1}-2\}$ for every $r\in[m]$;
		\item [(ii)] add $\al=\begin{pmatrix} A_1&\ldots&A_m\\a_1&\ldots&a_m\end{pmatrix}$ to $\cn^m(I,I)$.
	\end{itemize}
	Go to step 3.
	\item Increase $m$ by $1$.
	If $m\leq p$, set $\cn^m(I,I)=\emptyset$, and go to step 2. If $m=p+1$, RETURN $\cn^{nil}(I,I)=\cn^1(I,I)\cup\ldots\cup\cn^p(I,I)$ and STOP.
\end{enumerate}
\end{alg}
We note that, in Algorithm~\ref{acnnilii}, if $p\geq2$, then $\cn^p(I,I)=\emptyset$.
The following algorithm, given a tree $\tr=\ga(\mu)$, where $\mu\in N(\cn)$, as input,
returns $C_{nil}(\mu)$ as output.

\begin{alg}\label{acentnilmu}
	Let $\mu\in N(\cn)$, $\tr=\ga(\mu)$, and $h=\hei(\tr)$.
	\begin{itemize} 
		\item [1(a)] Set $C_{nil}(\mu_0)=\{\al\}$, where $1\al=1$.
		If $h\geq1$, go to step 1(b). If $h=0$, RETURN $C_{nil}(\mu)=C_{nil}(\mu_0)$ and STOP. 
		\item [1(b)] Set $C_{nil}(\mu_1)=\cn^{nil}(1\mu^{-1},1\mu^{-1})$ (see Algorithm~\ref{acnnilii}).
		If $h\geq2$, set $k=1$ and go to step 2. If $h=1$, RETURN $C_{nil}(\mu)=C_{nil}(\mu_1)$ and STOP.

	\end{itemize}
Steps 2--6 are exactly as in Algorithm~\ref{acentmu} except that we replace each $C(\mu_{k+1})$ with $C_{nil}(\mu_{k+1})$ and each $\cnst(y_j\mu^{-1},x_j\mu^{-1})$ with $\cn^{nil}(y_j\mu^{-1},x_j\mu^{-1})$.
\end{alg}
It follows from Notation~\ref{ncnnilii}, step 1(b), and induction on $k$ that for every $\al$ constructed
by Algorithm~\ref{acentnilmu}, the only point fixed by $\al$ is $1$, that is, $\al\in C_{nil}(\mu)$.

\begin{example}\label{ecentnilmu}
	Consider  $\mu=\left( \begin{smallmatrix}
	1 & 2 & 3 & 4 & 5 & 6 & 7 & 8 & 9 \\
	1 & 1 & 1 & 2 & 2 & 3 & 5 & 6 & 6 
	\end{smallmatrix}\right) \in N(\mathcal{C}_{9})$ as in Example~\ref{ecentmu} (see Figure~\ref{fcentmu}).
We can calculate $C_{nil}(\mu)$ using Algorithm~\ref{acentnilmu} (without calculating the whole $C(\mu)$):
\[
\begin{aligned}
C_{nil}(\mu_0)=\{&1\}, \ C_{nil}(\mu_1)=\{111,\,112\},\\	
C_{nil}(\mu_2)=\{
&111111,\,111112,\,111113,\,111122,\, 111123,\, 111133,\,111222,\,111223,\\
&111233,\,111333,\,112224,\,112225,\,112234,\,112235,\,112334,112335\},\\
C_{nil}(\mu)=C_{nil}(\mu_3)=\{
&111111111,\,111111112,\,111111113,\,111111122,\,111111123,\,111111133,\\
&111111222,\,111111223,\,111111233,\,111111333,\,111112244,\,111112245,\\
&111112255,\,111112344,\,111112345,\,111112355,\,111113366,\,111122444,\\
&111122445,\,111122455,\,111122555,\,111123466,\,111123566,\,111133666,\\
&111222444,\,111222445,\,111222455,\,111222555,\,111223466,\,111223566,\\
&111233666,\,111333666,\,112225577,\,112235677,\,112335677\},
\end{aligned}
\]
so $|C_{nil}(\mu)|=35$.
\end{example}

\subsection{Algorithm for constructing $\cid(\mu)$}\label{snil4}

Algorithm~\ref{acentmu} can also be modified to construct $\ceidmu=C(\mu)\cap E(\cn)$, but the modification is more complicated
than the one presented in Algorithm~\ref{acentnilmu} for calculating $\cenilmu$. Note that, unlike $\cenilmu$, $\ceidmu$ is a subset,
not necessarily a subsemigroup, of $C(\mu)$. 

\begin{rem}\label{rcnijc}
Let $I=\{i_1<\cdots<i_p\}$ and $J=\{j_1<\cdots<j_q\}$, where $p,q\geq1$, be chains of $\xn$. We can extend Notation~\ref{ncnstij}
from intervals to chains in an obvious way:
\[
\cnst(J,I)=\{\al:J\to I:\mbox{$\al$ is decreasing and order-preserving}\}.
\]
Lemmas~\ref{lcnstii} and~\ref{lcnstij} are true for chains ($I<J$ means that $i_p<j_1$), and Algorithms~\ref{acnii} and~\ref{acnij}
work with this extended notation.
\end{rem}

\begin{nota}\label{nvepia}
Let $I$ be an interval of $\xn$, and let $A$ be any subset of $I$ that contains the smallest element of $I$. We denote by $\vep(I,A)$
the unique idempotent $\vep\in\cn(I)$ such that $\im(\vep)=A$.

\end{nota}
For a set $A$, we denote by $\mathcal{P}(A)$ the power set of $A$.
The algorithm below, given a tree $\tr=\ga(\mu)$, where $\mu\in N(\cn)$, as input,
returns $\ceidmu$ as output.

\begin{alg}\label{acentidmu}
	Let $\mu\in N(\cn)$, $\tr=\ga(\mu)$, and $h=\hei(\tr)$.
	\begin{enumerate} 
		\item Set $C_{id}(\mu_0)=\{\al_1\}$, where $1\al_1=1$.
		If $h=0$, RETURN $C_{id}(\mu)=C_{id}(\mu_0)$ and STOP. If $h\geq1$, set $k=0$ ($k$ will range from $0$ to $h$) and go to step 2.
		\item At this point, $k<h$, and we have calculated $C_{id}(\mu_k)=\{\al_1,\ldots,\al_t\}$.  Find the set $W_k=\{y_1,\ldots,y_s\}$ of all elements
		$y\in V_k(\tr)$ such that $y\mu^{-1}\ne\emptyset$. Set $C_{id}(\mu_{k+1})=\emptyset$, $i=1$ ($i$ will range from $1$ to $t$), and go to step 3.
           \item For every $j\in [s]$, set $x_j=y_j\al_i$. If $x_j\mu^{-1}=\emptyset$ for some $j\in[s]$, go to step~11.
		Otherwise, go to step~4.
            \item Find the set $W^i_k=\{z_1,\ldots,z_l\}$ of all elements $z\in W_k$ such that $z\al_i=z$. If $W^i_k=\emptyset$, go to step~9. Otherwise, go to step~5.
            \item Find the set $\mathcal{P}_k$ of all elements $(A_1,\ldots,A_l)\in\mathcal{P}(z_1\mu^{-1})\times\cdots\times\mathcal{P}(z_l\mu^{-1})$
            such that for every $r\in[l]$, if $z$ is the smallest element of $z_r\mu^{-1}$, then $z\in A_r$. Go to step~6.
            \item If $\mathcal{P}_k=\emptyset$ (which will not happen when the jump is from step 5), go to step 11. Otherwise, select any $(A_1,\ldots,A_l)\in \mathcal{P}_k$ and go to step 7.
            \item For every $j\in[s]$, calculate 
            \[ \{\bt_1^j,\ldots,\bt_{v_j}^j\}= \left\{\begin{array}{lll}
            \cnst(y_j\mu^{-1},x_j\mu^{-1}\cap\,\im(\al_i))&\mbox{ if $x_j\notin W_k$ (so $y_j\notin W^i_k$),}\\
            \cnst(y_j\mu^{-1},A_r)&\mbox{ if $y_j\ne x_j=z_r\in W^i_k$ for some $r\in[l]$,}\\
            \{\vep(y_j\mu^{-1},A_r)\}&\mbox{ if $y_j=x_j=z_r\in W^i_k$ for some $r\in[l]$.}
            \end{array} \right.
            \]
            Go to step 8.
            \item For every $(j_1,\ldots,j_s)\in[v_1]\times\ldots\times[v_s]$,
		add $\al_i\joi\bt^1_{j_1}\joi\ldots\joi\bt^s_{j_s}$ to $C_{id}(\mu_{k+1})$ if it is order-preserving, remove $(A_1,\ldots,A_l)$ from $\mathcal{P}_k$, and go to step 6.
		\item For every $j\in[s]$, calculate $\cnst(y_j\mu^{-1},x_j\mu^{-1}\cap\im(\al_i))=\{\bt_1^j,\ldots,\bt_{v_j}^j\}$, and go to step 10.
		\item For every $(j_1,\ldots,j_s)\in[v_1]\times\ldots\times[v_s]$,
		add $\al_i\joi\bt^1_{j_1}\joi\ldots\joi\bt^s_{j_s}$ to $C_{id}(\mu_{k+1})$ if it is order-preserving, and go to step 11.
		\item Increase $i$ by 1. If $i\leq t$, go to step 3. If $i=t+1$, go to step 12.
		\item Increase $k$ by $1$. If $k<h$, go to step 2. If $k=h$, RETURN $C_{id}(\mu)=C_{id}(\mu_h)$ and STOP.
	\end{enumerate}
\end{alg}
Using Theorem~\ref{tcnilg}, the fact that $\al\in\cn$ is an idempotent if and only if $\im(\al)=\fix(\al)$, and induction on $k$, it is straightforward to show that Algorithm~\ref{acentidmu} constructs 
all elements of $\ceidmu$ (without repetitions).

\begin{example}\label{ecentidmu}
	Consider again  $\mu=\left( \begin{smallmatrix}
	1 & 2 & 3 & 4 & 5 & 6 & 7 & 8 & 9 \\
	1 & 1 & 1 & 2 & 2 & 3 & 5 & 6 & 6 
	\end{smallmatrix}\right) \in N(\mathcal{C}_{9})$ as in Example~\ref{ecentmu} (see Figure~\ref{fcentmu}).
The tree in Figure~\ref{fcentidmu} presents the calculations of $C_{id}(\mu_0)$, $C_{id}(\mu_1)$, $C_{id}(\mu_2)$, and $C_{id}(\mu_3)=C_{id}(\mu)$ (main loop, $k=0,1,2,3$, of Algorithm~\ref{acentidmu}). 
The $k$-level vertices of the tree are the elements of $C_{id}(\mu_k)$. 

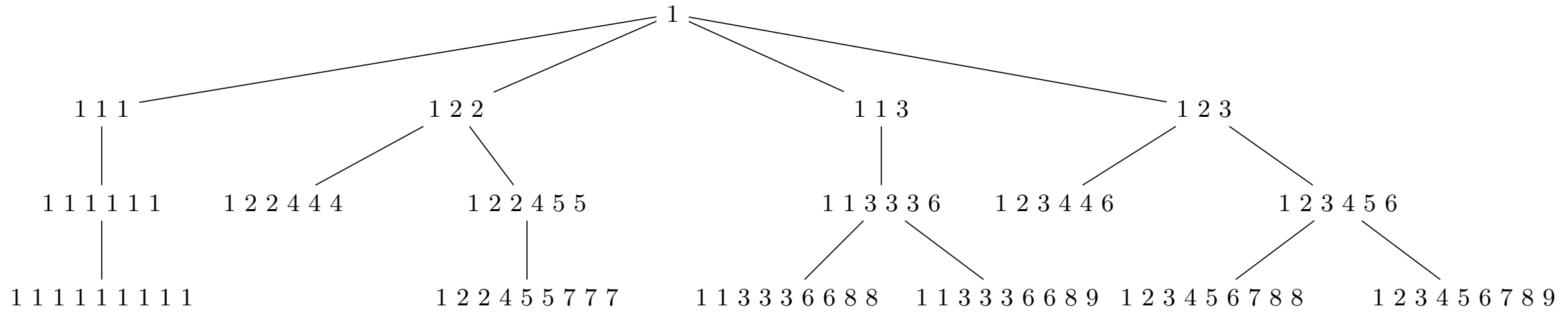
\begin{figure}[htbp]
	\centering
	\resizebox{\textwidth}{!}{
		\begin{tikzpicture}[
		every node/.style={font=\small, align=center},
		level distance=12mm,
		level 1/.style={sibling distance=45mm},
		level 2/.style={sibling distance=38mm}, 
		level 3/.style={sibling distance=32mm}, 
		edge from parent/.style={draw, -}
		]
		
		\node {1}
		child [xshift=-5mm]{node {1 1 1}
			child {node {1 1 1 1 1 1}
				child {node {1 1 1 1 1 1 1 1 1}}
			}
		}
		child [xshift=-5mm]{node {1 2 2}
			child[xshift=-3mm] {node {1 2 2 4 4 4}}
			child[xshift=-10mm] {node {1 2 2 4 5 5} 
				child {node {1 2 2 4 5 5 7 7 7}}
			}
		}
		child[xshift=4mm] {node {1 1 3} 
			child {node {1 1 3 3 3 6}
				child[xshift=4mm] {node {1 1 3 3 3 6 6 8 8}}
				child {node {1 1 3 3 3 6 6 8 9}}
			}
		}
		child {node {1 2 3}
			child {node {1 2 3 4 4 6}}
			child[xshift=-2mm] {node {1 2 3 4 5 6}
				child {node {1 2 3 4 5 6 7 8 8}}
				child {node {1 2 3 4 5 6 7 8 9}}
			}
		};
		
		\end{tikzpicture}
	}
	\caption{Tree presenting calculation of $C_{id}(\mu)$.}
	\label{fcentidmu}
\end{figure}

The elements $122444$ and $123446$ of $C_{id}(\mu_2)$ were ``rejected'' by step~3 of Algorithm~\ref{acentidmu}.
Note that $|C_{id}(\mu)|=6$.
\end{example}

\subsection{Nilpotents $\mu$ for which $C(\mu)$ consists of idempotents and nilpotents}\label{snil5}

For $\mu=\left( \begin{smallmatrix}
	1 & 2 & 3 & 4 & 5 & 6 & 7 & 8 & 9 \\
	1 & 1 & 1 & 2 & 2 & 3 & 5 & 6 & 6 
	\end{smallmatrix}\right) \in N(\mathcal{C}_{9})$,
we have $C(\mu)=C_{nil}(\mu)\cup C_{id}(\mu)$ (see Examples~\ref{ecentmu}, \ref{ecentnilmu}, and~\ref{ecentidmu}). This is not true in general. For example, 
for $\mu=111223\in N(\mathcal{C}_6)$, $\al=122445$ is in $C(\mu)$
but it is neither nilpotent nor idempotent.

In this subsection, we describe the nilpotents $\mu\in\cn$ for which $C(\mu)$ consists of idempotents,
nilpotents, or idempotents and nilpotents. More precisely, let $\vep_n$ and $\thet_n$ be the identity and zero of 
$\cn$, respectively, and note that $\vep_n,\thet_n\in C(\gam)$ for every $\gam\in\cn$. We describe the nilpotents $\mu\in\cn$ such that:
\begin{align}
  C(\mu)&=C_{id}(\mu)\label{eqa};\\
 C(\mu)&\mbox{ contains no idempotents except $\vep_n$ and $\thet_n$}\label{eqb};\\
 C(\mu)&=C_{nil}(\mu)\cup\{\vep_n\}\label{eqc};\\
 C(\mu)&=C_{nil}(\mu)\cup C_{id}(\mu)\label{eqd}.
\end{align}
Equality \eqref{eqa} holds only in a trivial case: for every $\mu\in N(\cn)$, $C(\mu)=C_{id}(\mu)$ if and only if $n\leq2$. 
Note that for $n\leq2$, $C(\mu)=C_{id}(\mu)=\{\vep_n,\thet_n\}$.
It turns out that for every $\mu\in N(\cn)$,  \eqref{eqb} (no nontrivial idempotents) and  \eqref{eqc} (only nilpotents and the identity) are equivalent. This is shown in Theorem~\ref{tiii},
which also describes the nilpotents $\mu\in\cn$  for which  \eqref{eqb} and  \eqref{eqc} hold.
A description of $\mu\in N(\cn)$ for which  \eqref{eqd} holds
is much more difficult. It is given in Theorem~\ref{tiv}, which is the main result of this subsection.

\begin{defi}\label{dloopless}
A \emph{loopless tree} with root $x$ is the digraph obtained from a tree with root $x$ by removing the loop $(x,x)$.
\end{defi}

\begin{nota}\label{ntrsub}
Let $\tr$ be a tree with $V(\tr)=\xn$ and root $1$, and let $x$ be a vertex of $\tr$. We denote by
$\tr_x$ the subdigraph of $\tr$ spanned by $W$, where $W$ is the set of vertices $z$ of $\tr$ such that there exists 
a path in $\tr$ from $z$ to $x$.

The subdigraph $\tr_x$ is a tree with root $x$ if $x=1$, and 
it is a loopless tree with root~$x$ if $x\ne1$.
The notions of a vertex of $\tr_x$ at level $k$, the height $\hei(\tr_x)$ of $\tr_x$, and
 the set $V_k(\tr_x)$ of vertices of $\tr_x$ at level $k$, where $0\leq k\leq\hei(\tr_x)$,
carry over to $\tr_x$ from a tree in the obvious way.
\end{nota}

\begin{defi}\label{ddophom}
Let $\ga_1$ and $\ga_2$ be digraphs with $V(\ga_1),V(\ga_2)\subseteq\xn$.
A \emph{decreasing and order-preserving homomorphism} from $\ga_1$ to $\ga_2$
is any function $\lam:V(\ga_1)\to V(\ga_2)$ such that $\lam$ is a digraph homomorphism, it is decreasing
($x\lam\leq x$ for every $x\in V(\ga_1)$), and it is order-preserving (for all $x,y\in V(\ga_1)$, if $x\leq y$, then $x\lam\leq y\lam$).
\end{defi}

\begin{nota}\label{nlam1}
Let $\mu\in N(\cn)$ with $\tr=\ga(\mu)$.
Suppose that $r,t\in\xn$ with $r\leq t$,
\[
\pth=r_0\to\cdots\to r_k=r
\] is a path in $\tr$,
and $\str$ is a loopless subtree of $\tr$ with root $t$
such that $\hei(\str)\leq k$ and
$r_{k-\dst(v,t)}\leq v$ for every $v\in V(\str)$.
We denote by $\lam(\str,\pth)$ the function
$\lam:V(\str)\to V(\pth)$ defined by
$v\lam=r_{k-\dst(v,t)}$ for every $v\in V(\str)$.
\end{nota}

\begin{nota}\label{nlam2}
Let $\mu\in N(\cn)$ with $\tr=\ga(\mu)$. Suppose that $r,t\in\xn$ with $\lev(r)=\lev(t)$ in $\tr$ and $r\leq t$,
$\pth=r_0\to\cdots\to r_k=r$ and $\qth=t_0\to\cdots\to t_l=t$ are paths in $\tr$ such that 
$l\leq k$ and $r_{k-l}\leq t_0$.
We denote by $\lam(\qth,\pth)$ the function $\lam:V(\qth)\to V(\pth)$ such that $t_i\lam=r_{i+k-l}$ for every $i\in\{0,\ldots,l\}$.
\end{nota}

\begin{nota}\label{nlam3}
	Let $\mu\in N(\cn)$ with $\tr=\ga(\mu)$. Suppose that $r,t\in\xn$ with  $\lev(r)=\lev(t)=1$ in $\tr$ and $t\leq r$,
	$\pth=r_0\to\cdots\to r_k=r$ is a path in $\tr_r$, and $\str$ is a loopless subtree of $\tr$ with root $t$ such that 
	$\hei(\str)\leq k+1$.	
	We denote by $\lamp(\str,\pth)$ the function $\lamp:V(\str)\to V(\pth)\cup\{1\}$ such that $t\lamp=1$ and for every $v\in V(\str)$
	with $v\ne t$,
	$v\lamp=r_i$, where $i\in\{0,\ldots,k\}$ is such that $\lev(r_i)=\lev(v)-1$ in $\tr$.
\end{nota}

The following lemma follows immediately from the definitions of $\lam(\str,\pth)$, $\lam(\qth,\pth)$, and $\lamp(\str,\pth)$, and
the fact that $\tr=\ga(\mu)$, where $\mu\in N(\cn)$, is level-ordered (see Definition~\ref{dlevord} and Corollary~\ref{cnil}).

\begin{lemma}\label{lllp}
The functions $\lam(\str,\pth)$, $\lam(\qth,\pth)$, and $\lamp(\str,\pth)$ from {\rm Notations~\ref{nlam1}}, {\rm \ref{nlam2}},
and~{\rm\ref{nlam3}}
are decreasing and order-preserving homomorphisms.
\end{lemma}

\begin{defi}\label{dttrl}
Let $\mu\in N(\cn)$ with $\tr=\ga(\mu)$. Let  $x\in V(\tr)$ with $x\ne1$,  $\pth_x=w_0\to\cdots\to w_k=x$ be a path in $\tr_x$ with $k\geq1$
and $\ind(w_0)=0$, and let $v\in V(\tr_x)$. We say that $v$ is \emph{to the right} [\emph{left}] of $\pth_x$ if $v$ lies on $\pth_x$
or there exists a path $v_0\to\cdots\to v_m=w_i$ in $\tr$ such that $m,i\geq1$ and $v_{m-1}>w_{i-1}$ [$v_{m-1}<w_{i-1}$].

The concepts of ``to the right'' and ``to the left'' are well defined by Proposition~\ref{paincn}(1) (see Definition~\ref{dlevord}).
\end{defi}

\begin{nota}\label{ntrx}
Let $\tr=\ga(\mu)$, where $\mu\in\cn$ is a nilpotent,  $x\in V(\tr)$ with $x\ne1$, and $\pth_x=w_0\to\cdots\to w_k=x$ be a path in $\tr_x$ with $k\geq1$
and $\ind(w_0)=0$.
We introduce the following notation:
\begin{itemize} 
    \item[(a)] $\tzppx$ is the subdigraph of $\tr_x$ spanned by $W$, where $W$ is the set of vertices $v$ of $\tr_x$
such that $v=x$ or $v$ is to the right of $\pth_x$ and the path from $v$ to $x$ passes 
through~$w_{k-1}$;
    \item[(b)] $\toppx$ is the subdigraph of $\tr_x$ spanned by $W$, where $W$ is the set of vertices $v$ of $\tr_x$ such that $v=x$
or the path from $v$ to $x$ passes through a vertex $w\in V_1(\tr_x)$ with $w>w_{k-1}$;
    \item[(c)] $\tzmpx$ is the subdigraph of $\tr_x$ spanned by $W$, where $W$ is the set of vertices $v$ of $\tr_x$
such that $v=x$ or $v$ is to the left of $\pth_x$ and the path from $v$ to $x$ passes 
through~$w_{k-1}$;
    \item[(d)] $\tompx$ is the subdigraph of $\tr_x$ spanned by $W$, where $W$ is the set of vertices $v$ of $\tr_x$ such that $v=x$
or the path from $v$ to $x$ passes through a vertex $w\in V_1(\tr_x)$ with  $w<w_{k-1}$;
   \item[(e)] $\tppx$ is the union of $\tzppx$ and $\toppx$;
   \item[(f)]  $\tmpx$ is the union of $\tzmpx$ and $\tompx$.
 \end{itemize}
If $V(\tr_x)=\{x\}$,  we set  $\tzppx=\toppx=\tzmpx=\tompx=\tppx=\tmpx=\tr_x$.

The subdigraphs $\tzppx$, $\toppx$, $\tzmpx$, $\tompx$, $\tppx$, and  $\tmpx$ are loopless trees with root~$x$. 
If $k\geq1$, then every vertex of $\tppx$ [$\tmpx$] is to the right [left] of $\pth_x$.
\end{nota}

The following lemma follows immediately from Lemma~\ref{lend} and the definition of a path in a tree.

\begin{lemma}\label{lpath}
Let $\mu\in N(\cn)$ with $\tr=\ga(\mu)$, and let $\al\in C(\mu)$. Then, for every path
$x_0\to\cdots\to x_k$ in $\tr$, $x_0\al\to\cdots\to x_k\al$ is a walk in $\tr$. Moreover, if $k=0$
or $k\geq1$ and $x_{k-1}\al\ne1$, then the walk $x_0\al\to\cdots\to x_k\al$ is a path in $\tr$.
\end{lemma}

\begin{lemma}\label{lhei1}
Let $\tr=\ga(\mu)$, where $\mu\in N(\cn)$, and let $a,b\in V(\tr)$ with $a\ne1$. Suppose $\al\in C(\mu)$
is such that $b\al=a$. Then, $\hei(\tr_a)\geq\hei(\tr_b)$.
\end{lemma}
\begin{proof}
Let $m=\hei(\tr_b)$. Select any vertex $t\in V(\tr_b)$ such that $t=t_0\to\cdots\to t_m=b$ ($m\geq0$)
is a path from $t$ to $b$ in $\tr_b$. By Lemma~\ref{lpath}, $t_0\al\to\cdots\to t_m\al=b\al=a$ is a path in $\tr_a$,
which implies $\hei(\tr_a)\geq m=\hei(\tr_b)$.
\end{proof}

The following theorem shows that  \eqref{eqb} and  \eqref{eqc} are equivalent and describes the nilpotents $\mu\in\cn$
for which  \eqref{eqb} and  \eqref{eqc} are true.

\begin{theorem}\label{tiii}
Let $n\geq3$ and $\mu\in N(\cn)$ with $\tr=\ga(\mu)$.  Then, the following conditions
are equivalent:
\begin{itemize} 
  \item [\rm(1)] $C(\mu)$ contains no idempotents except $\vep_n$ and $\thet_n$;
  \item [\rm(2)] $C(\mu)=C_{nil}(\mu)\cup\{\vep_n\}$;
  \item [\rm(3)] $|V_1(\tr)|=1$ and for all $x,y,z\in\xn$ with $x\ne1$, if $y,z\in V_1(\tr_x)$ and $y<z$, then $\hei(\tr_y)<\hei(\tr_z)$.
 \end{itemize}
\end{theorem}
\begin{proof}
($(1)\Rightarrow (2)$)
Suppose to the contrary that (1) holds and there is $\al\in C(\mu)$ such that $\al$ is not a nilpotent and $\al\ne\vep_n$.
Then, there exist $x,y\in\xn$ such that $x\ne1$, $x\al=x$, and $y\al\ne y$. Since $\cn$ is finite, $\al^k$ is an idempotent
for some $k\geq1$. Set $\vep=\al^k$ and note that $\vep\in C(\mu)$. By \cite[Prop.~2.3]{Ko2022},
$\fix(\al^k)=\fix(\al)$, which implies $x\in\fix(\vep)$ and $y\notin\fix(\vep)$. Hence, $\vep\ne\vep_n$ and $\vep\ne\thet_n$,
which is a contradiction.

($(2)\Rightarrow (3)$) Assume (2) holds. Suppose to the contrary  that $|V_1(\tr)|\geq2$.
Let $x=\min(V_1(\tr))$, $y=\max(V_1(\tr))$, and note that $x<y$. Fix $e\in V_1(\tr)$ such that
$\hei(\tr_e)\geq\hei(\tr_t)$ for every $t\in V_1(\tr)$, and a path $\pth_e$ of maximum length in $\tr_e$.
Define $\al:\xn\to\xn$ by:
\vskip -5mm
\begin{align*}
\al|_{V(\tr_s)}&=\mbox{$\lamp(\tr_s,\pth_e)$ for every $s\in V_1(\tr)$ with $s<e$;}\\[-1mm]
\al|_{V(\tppe)}&=\mbox{$\lam(\tppe,\pth_e)$;}\\[-1mm]
\al|_{V(\tr_s)}&=\mbox{$\lam(\tr_s,\pth_e)$ for every $s\in V_1(\tr)$ with $s>e$;}\\[-1mm]
\end{align*}
\noindent and $p\al=p$ for every other $p\in\xn$ (see Notations~\ref{nlam1}, \ref{nlam3}, and~\ref{ntrx}). 
Then, $\al\in\cn$ (by Lemma~\ref{lllp}) and $\al\in C(\mu)$ (by Lemma~\ref{lend}). 
Moreover, $\al$ is not a nilpotent since $e\ne1$ and $e\al=e$ (by the definition of $\al|_{V(\tppe)}$).
If $e=x$, then $\al\ne\vep_n$ since $x<y$ and $y\al=x$ (by the definition of $\al|_{V(\tr_s)}$ for $s\in V_1(\tr)$ with $s>e$)).
If $e>x$, then $\al\ne\vep_n$ since $x\ne1$  and $x\al=1$ (by the definition of $\al|_{V(\tr_s)}$ for $s\in V_1(\tr)$ with $s<e$).
This contradicts~(2). We have proved that $|V_1(\tr)|\leq1$, and so $|V_1(\tr)|=1$ since $n\geq3$.

 Now, suppose to the contrary that
there are $x,y,z\in\xn$ such that $x\ne1$, $y,z\in V_1(\tr_x)$, $y<z$, and $\hei(\tr_y)\geq\hei(\tr_z)$.
Let $a$ be the largest element in $V_1(\tr_x)$ such that $a<z$ and $\hei(\tr_{a})\geq\hei(\tr_z)$.
Fix a path $\pth_{a}$ of maximum length in $\tr_{a}$. 
Define $\vep:\xn\to\xn$ by:
\vskip -5mm
\begin{align*}
\vep|_{V(\tppa)}&=\mbox{$\lam(\tppa,\pth_a)$ (possible by the choice of $\pth_a$)};\\[-1mm]
\vep|_{V(\tr_d)}&=\mbox{$\lam(\tr_d,\pth_a)$ for every $a<d\leq z$  (possible since $\hei(\tr_a)\geq\hei(\tr_d)$  for every $a<d\leq z$)};\\[-1mm]
\end{align*}
\vskip -5mm
\noindent and  $p\vep=p$ for every other $p\in\xn$  (see Notations~\ref{nlam1} and~\ref{ntrx}). Then, $\vep\in\cn$
 (by Lemma~\ref{lllp}), $\vep\in C(\mu)$ (by Lemma~\ref{lend}), and $\vep$ is an idempotent by its definition.
Further, $\vep\ne\vep_n$ (since $z\vep=a<z$) and $\vep\ne\thet_n$ (since $a\vep=a\ne1$), which contradicts~(2).

($(3)\Rightarrow (1)$) Assume (3) holds. Then, in particular, $V_1(\tr)=\{2\}$.
Suppose to the contrary that there is an idempotent $\vep\in C(\mu)$ such that $\vep\ne\vep_n,\thet_n$.
Since $\vep\ne\vep_n$, there exist $y,z\in\xn$ such that $y<z$, $z\vep=y$, and $y\vep=y$. We may assume that $z$ is the smallest element
of $\xn$ for which this happens. 

Suppose to the contrary that $y=1$. Then, $z=2$ by the minimality of $z$, and so $2\vep=1$. Since $\vep\ne\theta_n$,
some $t>1$ is fixed by $\vep$. The unique path from $t$ to $1$ passes through $2$,
and so, by $t\vep=t$ and Lemma~\ref{lpath}, $2\vep=2$, contradicting $2\vep=1$.
Hence, $y>1$.

Let $y_1=y\mu$ and $z=z_0\to\cdots\to z_l=1$ be the unique path in $\tr$ from $z$ to $1$.
Then, $l\geq2$ (since $z>y>1$), $z_{l-1}\in V_1(\tr)$,  $z_1\vep=y_1$ (by Lemma~\ref{lend}), and $z_{l-1}\vep\leq z_1\vep=y_1$ (since $z_{l-1}\leq z_1$).
Suppose to the contrary that  $y\in V_1(\tr)$. Then, $y_1=y\mu=1$ and $z_{l-1}\ne y$ (since $y\vep=y>1$ and $z_{l-1}\vep\leq y_1=1$).
Thus, $|V_1(\tr)|\geq2$, which is a contradiction since $|V_1(\tr)|=1$.

Hence, $y\notin V_1(\tr)$ and $y>1$, and so $y_1>1$. Since $z_1\vep=y_1$, we must have $z_1=y_1$ by the minimality of $z$. Let $x=z_1=y_1$.
Then, $x\ne1$, $y,z\in V_1(\tr_x)$, $y<z$, and $z\vep=y$. By Lemma~\ref{lhei1}, $\hei(\tr_y)\geq\hei(\tr_z)$. On the other hand, by (3),
$\hei(\tr_y)<\hei(\tr_z)$. We have obtained a contradiction, and so (1) holds.
\end{proof}

We now turn to the problem of describing the nilpotents $\mu\in\cn$ for which  \eqref{eqd} is true, that is, the nilpotents $\mu\in\cn$ for which  
$C(\mu)=C_{nil}(\mu)\cup C_{id}(\mu)$. It will be more convenient (and the description will be clearer) to describe, equivalently,
$\mu\in N(\cn)$ for which $C(\mu)\ne C_{nil}(\mu)\cup C_{id}(\mu)$, that is, the nilpotents $\mu\in\cn$ such that
there exists at least one $\al\in C(\mu)$ that is neither nilpotent nor idempotent.

It is easy to see that if $n\leq3$, then every element of $\cn$ is a nilpotent or an idempotent. From now on,
we will assume that $n\geq4$.

\begin{lemma}\label{l1oru}
	Let $n\geq4$ and $\mu\in N(\cn)$ with $\tr=\ga(\mu)$. Then, $C(\mu)\ne C_{nil}(\mu)\cup C_{id}(\mu)$
if and only if there exists $\al\in C(\mu)$ that satisfies at least one of the following conditions:
\begin{itemize} 
  \item[\rm(1)] $w\al=v$, $v\al=1$, and $u\al=u$ for some $w,v,u\in\xn$ with $u>w>v>1$;
  \item[\rm(2)] $w\al=v$, $v\al=u$, and $u\al=u$ for some $w,v,u\in\xn$ with $w>v>u>1$.
\end{itemize}
\end{lemma}
\begin{proof} 
Suppose that $C(\mu)\ne C_{nil}(\mu)\cup C_{id}(\mu)$. Then, there exists $\al\in C(\mu)$ that is neither nilpotent nor idempotent.
Since $\al$ is not an idempotent and $\al$ is order-preserving, there are $w,v,v_1\in\xn$ such that
$w>v>v_1$, $w\al=v$ and $v\al=v_1$. Let $w$ be the smallest element for which such $w,v,v_1$ exist.
Then, $v_1\al=v_1$. Suppose $v_1=1$. Since $\al$ is not a nilpotent, there is $u\in\xn$ such that $u\ne1$ and $u\al=u$.
Suppose to the contrary that $u \leq w$. Then, since $\al$ is order-preserving,  $u = u\al \leq w\al = v$. This implies
$u \leq v$, and consequently $u = u\al \leq v\al = 1$, which contradicts $u \ne 1$. Therefore, $u > w$,
which yields $u > w > v > 1$.
Thus, (1) is satisfied. If $v_1>1$, setting $u=v_1$, we see that (2) is satisfied .

Conversely, if $\al\in C(\mu)$ satisfies (1) or (2), then $\al$ is not a nilpotent (since it fixes $u\ne1$) and is not an idempotent (since $w\al=v$ and $v\al\ne v$). Hence, $C(\mu)\ne C_{nil}(\mu)\cup C_{id}(\mu)$.
\end{proof}

The following lemma
describes $\mu\in N(\cn)$ for which there exists $\al\in C(\mu)$ such that (1) of Lemma~\ref{l1oru} holds.
\begin{lemma}\label{lwv1u}
	Let $n\geq4$ and $\mu\in N(\cn)$ with $\tr=\ga(\mu)$. Then, there exists $\al\in C(\mu)$ with
	$w\al=v$, $v\al=1$, and $u\al=u$ for some $w,v,u\in\xn$ such that $u>w>v>1$
	if and only if there are $x,y,z\in V_1(\tr)$ such that $x<y<z$,  and the following conditions hold:
\begin{itemize} 
  \item[\rm(a)] $\hei(\tr_y)\leq\hei(\tr_x)$;
  \item[\rm(b)] there exists $e\in V_1(\tr)$ such that $e\geq z$ and for every $s\in V_1(\tr)$, $\hei(\tr_s)\leq\hei(\tr_e)+1$
if $s<y$, and $\hei(\tr_s)\leq\hei(\tr_e)$ if $s\geq z$.
\end{itemize}
\end{lemma}
\begin{proof} 
$(\Rightarrow)$
 Suppose that there exists $\al\in C(\mu)$ that satisfies the given conditions. Let 
$w=w_0\to\cdots\to  w_{p}\to 1$, $v=v_0\to\cdots\to v_q\to 1$, and $u=u_0\to\cdots\to u_r\to 1$ be
the unique paths from $w$ to $1$, $v$ to $1$, and $u$ to $1$. Set $x=v_q$, $y=w_p$, and $z=u_r$.
Then,  $z\al=z$ (by $u_0\al=u\al=u=u_0$ and Lemma~\ref{lpath}),
$x\al=1$ (since $x\leq v$, and so $x\al\leq v\al=1$), and $x<z$ (since otherwise $1=x\al\geq z\al=z>1$).
Suppose to the contrary that $p\geq1$. 
Then, $w > z$ (since $\lev(w)\geq 2 > 1 = \lev(z)$), and so $v = w\al \geq z\al = z$.
Hence $1 = v\al \geq z\al = z > 1$.
This is a contradiction, and so $p=0$, which implies $w=w_0=w_p=y$. Hence, $v=x$ since $v=w\al\leq w=y$, and so $\lev(v)$ cannot be greater than
$\lev(y)=1$.
We have $x<y$ (since otherwise $1<v=w\al=y\al\leq x\al=1$) and $y<z$ (since otherwise $x=v=w\al=y\al\geq z\al=z>x$).
Hence, $x<y<z$.

We now have $x,y,z\in V_1(\tr)$ such that $x<y<z$, $x\al=1$, $y\al=x$, and $z\al=z$. We may assume that 
$y$ is the smallest element of $V_1(\tr)$ for which such $x,y,z$ exist.
Condition (a) holds by $y\al=x$ and Lemma~\ref{lhei1}.

To prove (b), let $e\in V_1(\tr)$ be such that  $e\geq z$ and $\hei(\tr_e)\geq\hei(\tr_t)$ for every $t\in V_1(\tr)$ with $t\geq z$.
Let $s\in V_1(\tr)$ with $s<y$. If $s\leq x$, then $s\al\leq x\al=1$, so $s\al=1$.
Suppose $x<s<y$. Then, $1=x\al\leq s\al\leq y\al=x$, so $s\al=1$ or $s\al=x$. But the latter is impossible by
the minimality of $y$. We have proved that $s\al=1$.
If $\hei(\tr_s)=0$, then clearly $\hei(\tr_s)\leq\hei(\tr_e)+1$. Suppose $\hei(\tr_s)\geq1$. Then, we can select
$s_1\in V_1(\tr_s)$ such that  $\hei(\tr_s)=\hei(\tr_{s_1})+1$. Let $s_2=s_1\al$.
Since $\lev(s_1)=2>1=\lev(z)$, we have $s_1>z$, and so $s_2=s_1\al\geq z\al=z$. Moreover, since $s_1\mu=s$ and $\al\in C(\mu)$, we have $s_2\mu=(s_1\al)\mu=(s_1\mu)\al=s\al=1$. Since $s_2\geq z>1$, it follows that $s_2\in V_1(\tr)$.
Hence, by the choice of $e$ and Lemma~\ref{lhei1},
$\hei(\tr_s)=\hei(\tr_{s_1})+1\leq\hei(\tr_{s_2})+1\leq\hei(\tr_e)+1$.

$(\Leftarrow)$ Conversely, suppose there exist $x, y, z \in V_1(\tr)$ such that $x < y < z$,
and conditions (a) and (b) hold.
Fix paths $\pth_x$ and $\pth_e$ of maximum length in $\tr_x$ and $\tr_e$, respectively. 
Define $\al:\xn\to\xn$ by:

\vskip -5mm
\begin{align*}
\al|_{V(\tr_s)}&=\mbox{$\lamp(\tr_s,\pth_e)$ for every $s\in V_1(\tr)$ with $s<y$ (possible by (b));}\\[-1mm]
\al|_{V(\tr_y)}&=\mbox{$\lam(\tr_y,\pth_x)$ (possible by (a));}\\[-1mm]
\al|_{V(\tppe)}&=\mbox{$\lam(\tppe,\pth_e)$;}\\[-1mm]
\al|_{V(\tr_s)}&=\mbox{$\lam(\tr_s,\pth_e)$ for every $s\in V_1(\tr)$ with $s>e$ (possible by $e\geq z$ and (b));}\\[-1mm]
\end{align*}
\vskip -5mm
\noindent and $p\al=p$ for every other $p\in\xn$ (see Notations~\ref{nlam1}, \ref{nlam3}, and \ref{ntrx}). 
Then, $\al\in\cn$ (by Lemma~\ref{lllp}) and $\al\in C(\mu)$ (by Lemma~\ref{lend}). 

Let $w=y$, $v=x$, and $u=z$. Then, $u>w>v>1$ (since $z>y>x>1$),
$w\al=y\al=x=v$  (by the definition of $\al|_{V(\tr_y)}$), $v\al=x\al=1$ (by the definition of $\al|_{V(\tr_s)}$ for $s\in V_1(\tr)$ with $s<y$),
and $u\al=z\al=z=u$ (since $\al$ fixes $z$).
\end{proof}

Let $\mu\in N(\cn)$ and suppose there is $\al\in C(\mu)$ such that
$w\al=v$, $v\al=u$, and $u\al=u$ for some $w,v,u\in\xn$ with $w>v>u>1$ (see (b) of Lemma~\ref{l1oru}).
By Corollary~\ref{cnil} and Definition~\ref{dlevord}, $\lev(w)\geq\lev(v)\geq\lev(u)>0$. Therefore, there are four cases
to consider: 
\begin{align}
\lev(w)&=\lev(v)=\lev(u)>0;\label{ecase1}\\
\lev(w)&=\lev(v)>\lev(u)>0;\label{ecase2}\\
\lev(w)&>\lev(v)=\lev(u)>0;\label{ecase3}\\
\lev(w)&>\lev(v)>\lev(u)>0.\label{ecase4}
\end{align}

The following lemma shows that \eqref{ecase4} is impossible. 

\begin{lemma}\label{lth1}
Let $n\geq4$ and $\mu\in N(\cn)$ with $\tr=\ga(\mu)$. Suppose $\al\in\cn$ with
$w\al=v$, $v\al=u$, and $u\al=u$ for some $w,v,u\in\xn$ such that $\lev(w)>\lev(v)>\lev(u)>0$.
Then, $\al\notin C(\mu)$.
\end{lemma}

\begin{proof}
Suppose to the contrary that $\al\in C(\mu)$. Consider the unique paths from $w$, $v$, and $u$ to $1$ in $\tr$:
$w=w_0\to\cdots\to w_k=1$, $v=v_0\to\cdots\to v_l=1$, and $u=u_0\to\cdots\to u_m=1$. 
Then,  $k>l>m\geq1$ (since $\lev(w)>\lev(v)>\lev(u)>0$), $u_{m-1}\al=u_{m-1}$
(by $u_0\al=u\al=u=u_0$ and Lemma~\ref{lpath}),
$v_m\al=u_m=1$ (by $v_0\al=v\al=u=u_0$ and Lemma~\ref{lpath}),
$\lev(w_{l-1})\geq2$ (since $k>l$ and $\lev(w_{k-1})=1$), and
$w_{l-1}\al=v_{l-1}$ (by $w_0\al=w\al=v=v_0$ and Lemma~\ref{lpath}).
Thus, $w_{l-1}>u_{m-1}$, and so $v_{l-1}=w_{l-1}\al\geq u_{m-1}\al=u_{m-1}>1$.
Hence, $v_{l-1}\al\geq u_{m-1}\al=u_{m-1}>1$.
On the other hand, $v_{l-1}\leq v_m$, and so $v_{l-1}\al\leq v_m\al=1$.
We obtained a contradiction, and so $\al\notin C(\mu)$.
\end{proof}

We now describe $\mu\in N(\cn)$ for which there is $\al\in C(\mu)$ such that $w\al=v$, $v\al=u$, and $u\al=u$
for some $w,v,u\in\xn$ that satisfy \eqref{ecase1}. (See (2) of Lemma~\ref{l1oru}.)

\begin{lemma}\label{lxy}
Let $\tr=\ga(\mu)$, where $\mu\in N(\cn)$, and $u,v,w\in\xn$ be such that
$\lev(u)=\lev(v)=\lev(w)$ and $u<v<w$. Suppose that $\al\in C(\mu)$ is such
that $w\al=v$, $v\al=u$, and $u\al=u$. Then, there are $a,b,c\in\xn$ such that
$\lev(a)=\lev(b)=\lev(c)$, $a<b<c$, $c\al=b$, $b\al=a$, $a\al=a$, and exactly one of the following conditions is satisfied:
\begin{itemize}
  \item [\rm(a)] there exists $x\in\xn$ such that $x\al=x$ and $a,b,c\in V_1(\tr_x)$;
  \item [\rm(b)] there exist $x,y\in\xn$ such that $\lev(x)=\lev(y)$, $x<y$, $y\al=x$, $a,b\in V_1(\tr_x)$, and $c\in V_1(\tr_y)$.
\end{itemize}
\end{lemma}
\begin{proof}
	Choose $a,b,c\in\xn$ satisfying $a<b<c$, $c\al=b$, $b\al=a$, $a\al=a$, and 
	$\lev(a)=\lev(b)=\lev(c) (=m)$, where $m$ is the smallest level at which these conditions are satisfied. Let  $x=a\mu$, $y=b\mu$, and $z=c\mu$. Then,
	$\lev(x)=\lev(y)=\lev(z)=m-1$, and $x\leq y\leq z$ by Corollary~\ref{cnil} and Definition~\ref{dlevord}.
	Since $\al$ commutes with $\mu$,
	\[
	x\al=(a\mu)\al=(a\al)\mu=a\mu=x,
	\quad
	y\al=(b\mu)\al=(b\al)\mu=a\mu=x,
	\quad
	z\al=(c\mu)\al=(c\al)\mu=b\mu=y.
	\]
The elements $x$, $y$, and $z$ cannot be pairwise distinct since $x<y<z$ would contradict the minimality condition
imposed on $a$, $b$, and $c$. Suppose to the contrary that $x<y=z$. 
Then $x=y\al=z\al=y$, contradicting $x<y$. Therefore, either $x=y=z$ or $x=y<z$. In the former case, $a,b,c\in V_1(\tr_x)$, so (a) holds.
In the latter case, $a,b\in V_1(\tr_x)$, $c\in V_1(\tr_z)$, and $z\al=y=x$, so (b) holds.
\end{proof}

\begin{defi}\label{dlex}
For an integer $m\geq1$, denote by $L_m$ the set of sequences $\lan a_1,\ldots,a_m\ran$ of positive integers.
We define the \emph{lexicographical ordering} $<_l$ on $L_m$ by: $\lan a_1,\ldots,a_m\ran<_l\lan b_1,\ldots,b_m\ran$
if there exists $p\in\{1,\ldots,m\}$
such that $a_p<b_p$ and $a_i=b_i$ for every $1\leq i<p$.
\end{defi}
The relation $<_l$ is a linear ordering on $L_m$ (\cite[Definition~5.9]{HrJe99}).

\begin{defi}\label{dmerge}
Let $\mu\in N(\cn)$ with $\tr=\ga(\mu)$, and let $x,y$ be vertices of $\tr$.
The \emph{merge point} of $x$ and $y$, denoted $\mrp(x,y)$, is the largest element $z\in\xn$
such that there are paths in $\tr$ from $x$ to $z$ and from $y$ to $z$.
\end{defi}

\begin{nota}\label{npre}
 Let $\mu\in N(\cn)$ with $\tr=\ga(\mu)$, and let $x,y\in\xn$ be such that there is no path in $\tr$
from $x$ to $y$ or from $y$ to $x$. We write $x\prec y$ if there are $w,z\in\xn$
such that there are paths from $x$ and $y$ to $w$ and $z$, respectively, $\lev(w)=\lev(z)$, and $w<z$.
\end{nota}

The following lemma describes $\mu\in N(\cn)$ for which there is $\al\in C(\mu)$ such that $w\al=v$, $v\al=u$, and $u\al=u$
for some $w,v,u\in\xn$ that satisfy $w>v>u>1$ and \eqref{ecase1}. (See (2) of Lemma~\ref{l1oru}.)

\begin{lemma}\label{lth4}
	Let $n\geq4$ and $\mu\in N(\cn)$ with $\tr=\ga(\mu)$. 
Then, there exists $\al\in C(\mu)$ with $w\al=v$, $v\al=u$, and $u\al=u$ for some $w,v,u\in\xn$ such that $\lev(w)=\lev(v)=\lev(u)$ and 
$w>v>u>1$ if and only if at least one of the following conditions is satisfied:
\begin{itemize} 
  \item[\rm(1)] there are $x,a,b,c\in\xn$ such that $a,b,c\in V_1(\tr_x)$, $a<b<c$,  $\hei(\tr_b)\geq\hei(\tr_c)$, and $\hei(\tr_a)\geq\hei(\tr_d)$
  for every $a<d<c$;
  \item[\rm(2)] there are $x,y,a,b,c\in\xn$ with $x<y$, $\lev(x)=\lev(y) (=k)$, $a,b\in V_1(\tr_x)$ with $a<b$, and $c\in V_1(\tr_y)$ such that
  the following conditions hold:
  \begin{itemize} 
    \item [\rm(a)] $c=\max(V_1(\tr_y))$ and $\hei(\tr_b)\geq\hei(\tr_c)$;
    \item [\rm(b)] for every $d\in\xn$, $\hei(\tr_a)\geq\hei(\tr_d)$ if $d\in V_1(\tr_x)$ and $d>a$ or $d\in V_1(\tr_y)$ and $d<c$;
    \item [\rm(c)] if $z\in V_k(\tr)$ with $x<z<y$, then $\hei(\tr_a)\geq\hei(\tr_d)$ for every $d\in V_1(\tr_z)$;
    \item [\rm(d)] if $z\in V_k(\tr)$ with $z>y$ and $\mrp(y,z)>\mrp(x,y)$, then $\hei(\tr_z)=0$. 
\end{itemize}
\end{itemize}
\end{lemma}
\begin{proof}
$(\Rightarrow)$
Suppose that there exists $\al\in C(\mu)$ that satisfies the given conditions. We consider two possible cases.
\vskip 1mm
\noindent {\bf Case 1.} There are $x,a,b,c\in\xn$ such that $a,b,c\in V_1(\tr_x)$, $a<b<c$, $c\al=b$, $b\al=a$, and $a\al=a$.
\vskip 1mm
We will prove that in Case~1, condition (1) is satisfied.
Assume that $\lan a,b,c\ran$ is the \emph{minimum} element of $L_3$ with respect to $<_l$
that satisfies these properties. We have $\hei(\tr_b)\geq\hei(\tr_c)$ by $c\al=b$ and Lemma~\ref{lhei1}.
Let $a<d<c$ and $d_1=d\al$. Suppose $a<d\leq b$. Then, $a=a\al\leq d\al\leq b\al=a$. Thus, $d\al=a$, and so
$\hei(\tr_a)\geq\hei(\tr_d)$ by Lemma~\ref{lhei1}. Suppose $b<d<c$. Then, $a=b\al\leq d_1\leq c\al=b$.
If $a<d_1$, then $a<d_1\leq b<d<c$, $d\al=d_1$, $d_1\al=a$ (by the foregoing argument), and $a\al=a$, which contradicts
the $<_l$-minimality of $\lan a,b,c\ran$. Hence, $d\al=d_1=a$, and so $\hei(\tr_a)\geq\hei(\tr_d)$ by Lemma~\ref{lhei1}.
\vskip 1mm
\noindent {\bf Case 2.} There are no $x,a,b,c\in\xn$ such that $a,b,c\in V_1(\tr_x)$, $a<b<c$, $c\al=b$, $b\al=a$, and $a\al=a$.
\vskip 1mm
We will prove that in Case~2, condition (2) is satisfied.
Select $a,b,c\in\xn$ such that $\lev(a)=\lev(b)=\lev(c) (=m)$,
$a<b<c$, $c\al=b$, $b\al=a$, and $a\al=a$, where $m$ is the smallest level at which these conditions are satisfied. (Such $a,b,c$ exist by the hypothesis about $\al$.) 
By Lemma~\ref{lxy}, there exist $x,y\in\xn$ such that $\lev(x)=\lev(y)$, $x<y$, $y\al=x$, $a,b\in V_1(\tr_x)$ with $a<b$, and $c\in V_1(\tr_y)$.
Moreover, since $a\mu=x$, $a\al=a$, and $\al\mu=\mu\al$, we have
$x\al=(a\mu)\al=(a\al)\mu=x$.
Set $k=m-1$ and note that $k\geq1$ and $x,y\in V_k(\tr)$.
Among all quintuples with these properties arising from triples $a,b,c$ at this minimal level $m$, assume that
$\lan x,y,a,b,c\ran$ is the \emph{maximum} element of $L_5$ with respect to $<_l$.
(Such a quintuple exists since $\xn$ is finite.)

To prove (a), suppose to the contrary that there exists $d\in V_1(\tr_y)$ such that $d>c$, and let $d_1=d\al$ and $d_2=d_1\al$. Since $y\al=x$ 
and $d\in V_1(\tr_y)$,
we have $d_1\in V_1(\tr_x)$ by Lemma~\ref{lend}. Since $d_1\in V_1(\tr_x)$, $c\in V_1(\tr_y)$, and $x<y$, we have $d_1<c$
(otherwise $y=c\mu\leq d_1\mu=x$).
Further, $d_1=d\al\geq c\al=b$ (since $d>c$),
$d_2=d_1\al\leq c\al=b$, and $d_2=d_1\al\geq b\al=a$ (since $d_1\geq b$).
If $d_1=b$, then $d>c$, $d\al=d_1=b$, $b\al=a$, and $a\al=a$, which contradicts the $<_l$-maximality of $\lan x,y,a,b,c\ran$.
If $d_1>b$ and $d_2=a$, then $d>c$, $d\al=d_1>b$, $d_1\al=d_2=a$, and $a\al=a$, which again contradicts the $<_l$-maximality of 
$\lan x,y,a,b,c\ran$.
Thus, $d_1>b$ and $d_2>a$. Since $a<d_2\leq b$, we have $a=a\al\leq d_2\al\leq b\al=a$, and so $d_2\al=a$.
Hence, $a,d_2,d_1\in V_1(\tr_x)$, $a<d_2<d_1$, $d_1\al=d_2$, $d_2\al=a$, and $a\al=a$, which contradicts the assumptions about $\al$
in Case~2. We have proved that $c=\max(V_1(\tr_y))$. Since $c\al=b$, we have $\hei(\tr_b)\geq\hei(\tr_c)$ by Lemma~\ref{lhei1}.

To prove (b), suppose that $d\in V_1(\tr_x)$ with $d>a$. If $a<d<b$, then $a=a\al\leq d\al\leq b\al=a$, and so $d\al=a$.
Suppose $d\geq b$, and let $d_1=d\al$ and $d_2=d_1\al$. Since $d\in V_1(\tr_x)$, $c\in V_1(\tr_y)$, and $x<y$, we have $d<c$
(otherwise $y=c\mu\leq d\mu=x$). Thus, since $c\al=b$, we have $a=a\al\leq d\al\leq c\al=b$, and so
$a\leq d_1\leq  b$. Thus, $a=a\al\leq d_1\al\leq b\al=a$, and so $d_2=a$.
Suppose to the contrary that $d_1\ne a$. Then, 
$a,d_1,d\in V_1(\tr_x)$, $a<d_1<d$, and
$d\al=d_1$, $d_1\al=a\al=a$,
which contradicts the assumptions about $\al$ in Case~2. 
Thus, $d\al=d_1=a$. Hence, by Lemma~\ref{lhei1}, $\hei(\tr_a)\geq\hei(\tr_d)$ for every $d>a$.
Now, suppose that $d\in V_1(\tr_y)$ and $d<c$. Let $d_1=d\al$. By $y\al=x$ and Lemma~\ref{lend}, $d_1\in V_1(\tr_x)$.
Since $b<d<c$, we have $a=b\al\leq d\al\leq c\al=b$, and so $a\leq d_1\leq b$. By the foregoing argument, $d_1\al=a$. 
Hence, by Lemma~\ref{lhei1}, $\hei(\tr_a)\geq\hei(\tr_{d_1})\geq\hei(\tr_d)$. 

To prove (c),  let $z\in V_k(\tr)$ with $x<z<y$. Suppose that $d\in V_1(\tr_z)$, and let $d_1=d\al$. Then,
$x=x\al\leq z\al\leq y\al=x$, and so $z\al=x$ and $d_1\in V_1(\tr_x)$. Moreover, $d_1=d\al\geq b\al=a$ (since $d>b$). If $d_1=a$, then $\hei(\tr_a)=\hei(\tr_{d_1})$; otherwise,
(b) gives $\hei(\tr_a)\geq\hei(\tr_{d_1})$. Hence, by Lemma~\ref{lhei1},
$\hei(\tr_a)\geq\hei(\tr_{d_1})\geq\hei(\tr_d)$.

To prove (d), suppose to the contrary that there exists $z\in V_k(\tr)$ such that $z>y$, $\mrp(y,z)>\mrp(x,y)$, and $\hei(\tr_z)>0$.
Select the smallest $z$ with these properties, and let $z'=z\al$.
We claim that $z'=x$. Since $y<z$, we have $x=y\al\leq z\al\leq z$, so $x\leq z'\leq z$.

Suppose to the
contrary that $z'=z$. Let $l\geq1$ be an integer such that there are paths $x=x_0\to\cdots\to x_l=\mrp(x,y)$ and
$y=y_0\to\cdots\to y_l=x_l$ in $\tr$. Since $\mrp(y,z)>\mrp(x,y)$, there is $1\leq q<l$ such that $y_q=\mrp(y,z)$
and there is a path $z=z_0\to\cdots\to z_q=y_q$ in $\tr$. (Note that this forces $l$ to be at least $2$.) By Lemma~\ref{lpath},
$y_q\al=x_q$ (since $y_0\al=y\al=x=x_0$) and $z_q\al=z_q$ (since $z_0\al=z\al=z'=z=z_0$). Now, $x_q=y_q\al=z_q\al=z_q=y_q$. But,
by the definition of $\mrp(x,y)$ and the fact that $q<l$, we have $x_q<y_q$, which is a contradiction.
Hence, $z'<z$. 

Suppose to the contrary that $y<z'<z$. Since $\tr$ is a level-ordered tree (see Proposition~\ref{paincn}), Definition~\ref{dmerge} gives
$\mrp(y,z')\geq \mrp(y,z)>\mrp(x,y)$, and so, by the foregoing argument,
$z'\al<z'$. By  Lemma~\ref{lhei1} and $\hei(\tr_z)>0$, we have $\hei(\tr_{z'})>0$. But this contradicts
the minimality of $z$. Hence,  $z'\leq y$. Thus, since $x<y<z$,
we have $x=x\al\leq z\al=z'\leq y$. This implies $x=x\al\leq z'\al\leq y\al=x$, and so $z'\al=x$. It follows that $z'=x$ since otherwise
we would have $x,z',z\in V_k(\tr)$ with $x<z'<z$, $z\al=z'$, $z'\al=x$, and $x\al=x$, which would contradict the minimality of $m=k+1$.

The claim has been proved; that is, $z\al=z'=x$. Since $\hei(\tr_z)>0$, there exists some $d\in V_1(\tr_z)$.
Let $d_1=d\al$ and $d_2=d_1\al$. Then, $d_1,d_2\in V_1(\tr_x)$ (since $z\al=x$ and $x\al=x$), $d_1=d\al\geq c\al=b$
(since $z > y$, $c \in V_1(\tr_y)$, and $d \in V_1(\tr_z)$, which implies $d > c$),
$d_2=d_1\al\leq c\al=b$ (since $d_1<c$), and $d_2=d_1\al\geq b\al=a$ (since $d_1\geq b$). Suppose to the contrary that $d_2=d_1$. Then,
$d_1\geq b\geq d_2=d_1$, and so $d_1=b$. But then, $a=b\al=d_1\al=d_2=d_1=b$, which contradicts $a<b$. Thus, $d_2<d_1$.

We now have $x<z$, $z\al=x$, $d_2,d_1\in V_1(\tr_x)$ with $d_2<d_1$, $d\in V_1(\tr_z)$, $d\al=d_1$, and $d_1\al=d_2$.
Let $d_3=d_2\al$. Since $a\leq d_2$ and $a\al=a$, we have $a\leq d_3\leq d_2$. Also,
$d_3\in V_1(\tr_x)$ by Lemma~\ref{lend}.
If $d_3<d_2$, then the finite decreasing chain
$d_1>d_2>d_3\geq d_3\al\geq\cdots$ remains in $V_1(\tr_x)$ and eventually reaches a fixed point
(distinct from $1$ since $1\notin V_1(\tr_x)$).
The last three distinct terms of this chain would contradict the assumptions about $\al$ in Case~2.
Hence $d_3=d_2$, and so $d_2=a$ (since $a\leq d_2\leq b$ and
$a=a\al\leq d_2\al=d_2\leq b\al=a$).
Thus, $x<z$, $z\al=x$, $d_2,d_1\in V_1(\tr_x)$ with $d_2<d_1$, $d\in V_1(\tr_z)$, $d\al=d_1$, $d_1\al=d_2$,
and $d_2\al=d_2$. Therefore, $\lan x,z,d_2,d_1,d\ran$ is a quintuple of the same type as
$\lan x,y,a,b,c\ran$, but it is larger with respect to $<_l$ since $y<z$. This contradicts the choice of
$\lan x,y,a,b,c\ran$. We have proved~(d).

$(\Leftarrow)$ Suppose (1) is satisfied. Fix paths $\pth_a$ and $\pth_b$ of maximum length in $\tr_a$ and $\tr_b$, respectively. 
Define $\al:\xn\to\xn$ by:
\vskip -5mm
\begin{align*}
\al|_{V(\tppa)}&=\mbox{$\lam(\tppa,\pth_a)$ (possible by the choice of $\pth_a$)};\\[-1mm]
\al|_{V(\tr_d)}&=\mbox{$\lam(\tr_d,\pth_a)$ for every $a<d<c$  (possible since $\hei(\tr_a)\geq\hei(\tr_d)$  for every $a<d<c$)};\\[-1mm]
\al|_{V(\tr_c)}&=\mbox{$\lam(\tr_c,\pth_b)$  (possible since $\hei(\tr_b)\geq\hei(\tr_c)$)};\\[-1mm]
\end{align*}
\vskip -5mm
\noindent and  $p\al=p$ for every other $p\in\xn$  (see Notations~\ref{nlam1} and~\ref{ntrx}). Then, $\al\in\cn$
 (by Lemma~\ref{lllp}) and $\al\in C(\mu)$ (by Lemma~\ref{lend}). Further, for $u=a$, $v=b$, and $w=c$,
we have $\lev(w)=\lev(v)=\lev(u)$ (since $a,b,c\in V_1(\tr_x)$), $w>v>u$ (since $c>b>a$), $w\al=c\al=b=v$ (by the definition of $\al|_{V(\tr_c)}$),
$v\al=b\al=a=u$ (by the definition of $\al|_{V(\tr_d)}$ for $a<d<c$), and $u\al=a\al=a=u$ (by the definition of $\al|_{V(\tppa)}$).

Suppose (2) is satisfied. For $v\in\xn$, denote by $\pth^v$ the unique path from $v$ to $1$ in $\tr$.
Fix paths $\pth_a$ and $\pth_b$ of maximum length in $\tr_a$ and $\tr_b$, respectively. 
Define $\al:\xn\to\xn$ by:
{\footnotesize
\begin{align*}
\al|_{V(\pth^z)}&=\mbox{$\lam(\pth^z,\pth^x)$ for every $z\in V_k(\tr)$ such that $x<z\leq y$};\\[-1mm]
\al|_{V(\pth^z)}&=\mbox{$\lam(\pth^z,\pth^x)$ for every $z\in V_j(\tr)$ such that $1\leq j<k$, $x\prec z\prec y$, and $\ind(z)=0$};\\[-1mm]
\al|_{V(\pth^z)}&=\mbox{$\lam(\pth^z,\pth^x)$ for every $z\in V_k(\tr)$ such that, $z>y$, and $\mrp(y,z)>\mrp(x,y)$};\\[-1mm]
\al|_{V(\pth^z)}&=\mbox{$\lam(\pth^z,\pth^x)$ for every $z\in V_j(\tr)$ such that $1\leq j<k$, $y\prec z$, $\mrp(y,z)>\mrp(x,y)$, and $\ind(z)=0$};\\[-1mm]
\al|_{V(\tppa)}&=\mbox{$\lam(\tppa,\pth_a)$ (possible by the choice of $\pth_a$)};\\[-1mm]
\al|_{V(\tr_d)}&=\mbox{$\lam(\tr_d,\pth_a)$ for every $d\in V_1(\tr_x)$ with $d>a$ (possible by (b))};\\[-1mm]
\al|_{V(\tr_d)}&=\mbox{$\lam(\tr_d,\pth_a)$ for every $d\in V_1(\tr_z)$, where $z\in V_k(\tr)$ with $x<z<y$ (possible by (c))};\\[-1mm]
\al|_{V(\tr_d)}&=\mbox{$\lam(\tr_d,\pth_a)$ for every $d\in V_1(\tr_y)$ with $d<c$ (possible by (b))};\\[-1mm]
\al|_{V(\tr_c)}&=\mbox{$\lam(\tr_c,\pth_b)$ (possible by (a))};\\[-1mm]
\end{align*}
\vskip -3mm
}
\noindent and  $p\al=p$ for every other $p\in\xn$  (see Notations~\ref{nlam1}, \ref{nlam2}, and~\ref{ntrx}). 
Let $z_1,z_2\in V_k(\tr)$ such that $z_1,z_2>y$, $\mrp(y,z_1)>\mrp(x,y)$, and $\mrp(y,z_2)\leq\mrp(x,y)$. Suppose to the contrary
that $z_2\leq z_1$. Thus,  since $\tr$ is a level-ordered tree, Definition~\ref{dmerge} 
gives $\mrp(y,z_2)\geq\mrp(y,z_1)>\mrp(x,y)$,
which is a contradiction. We have proved that for all $z_1,z_2\in V_k(\tr)$,
\begin{equation}\label{lth4e1}
\mbox{if $z_1,z_2>y$, $\mrp(y,z_1)>\mrp(x,y)$, and $\mrp(y,z_2)\leq\mrp(x,y)$, then $z_1<z_2$}.
\end{equation}
Then, $\al\in\cn$
 (by Lemma~\ref{lllp}, condition~(2)(d), and \eqref{lth4e1}), and $\al\in C(\mu)$ (by Lemma~\ref{lend}). 
Further, for $u=a$, $v=b$, and $w=c$,
we have $\lev(w)=\lev(v)=\lev(u)$ (since $a,b\in V_1(\tr_x)$, $c\in V_1(\tr_y)$, and $\lev(x)=\lev(y)$), $w>v>u$ (since $c>b>a$), 
$w\al=c\al=b=v$ (by the definition of $\al|_{V(\tr_c)}$),
$v\al=b\al=a=u$ (by the definition of $\al|_{V(\tr_d)}$ for $d\in V_1(\tr_x)$ with $d>a$), and $u\al=a\al=a=u$ (by the definition of $\al|_{V(\tppa)}$).
\end{proof}

The following lemma describes $\mu\in N(\cn)$ for which there is $\al\in C(\mu)$ such that $w\al=v$, $v\al=u$, and $u\al=u$
for some $w,v,u\in\xn$ such that $w>v$ and \eqref{ecase2} holds. (See (2) of Lemma~\ref{l1oru}.)

\begin{lemma}\label{lth2}
	Let $n\geq4$ and $\mu\in N(\cn)$ with $\tr=\ga(\mu)$. Then, there exists $\al\in C(\mu)$ with
	$w\al=v$, $v\al=u$, and $u\al=u$ for some $w,v,u\in\xn$ such that $\lev(w)=\lev(v)>\lev(u)>0$ and $w>v>u>1$
	if and only if there are $x,y,z\in V_1(\tr)$ such that $x<y<z$, $\hei(\tr_x)\geq\hei(\tr_y)\geq1$, and the following conditions hold:
\begin{itemize} 
  \item[\rm(a)] for every $s\in V_1(\tr)$, if $s>z$, then $\hei(\tr_s)\leq \hei(\tr_z)$;
  \item[\rm(b)] for every $s\in V_1(\tr)$, if $s<y$, then $\hei(\tr_s)\leq\hei(\tr_z)+1$.
\end{itemize}
\end{lemma}

\begin{proof} 
$(\Rightarrow)$
 Suppose that there exists $\al\in C(\mu)$ that satisfies the given conditions. Consider the unique paths from $w,v$, and $u$ to $1$: $w=w_0\to\cdots\to w_k\to 1$, $v=v_0\to\cdots\to v_k\to 1$, and $u=u_0\to\cdots\to u_m\to 1$ in $\tr$. Note that
$0\leq m<k$ and  $w_k,v_k,u_m\in V_1(\tr)$. Set $x=v_k$, $y=w_k$, and $z=u_m$. 
Since $u\al=u$, $\lev(w)=\lev(v)$, and $\lev(v)>\lev(u)$, Lemma~\ref{lpath} implies that $z\al=z$, $y\al=x$, 
$x\al=1$, and either $v_{k-1}\al=z$ or $v_{k-1}\al=1$. But
$v_{k-1}\al=1$ is impossible
since $v_{k-1}>z$, and so $v_{k-1}\al\geq z\al=z>1$. Hence, $v_{k-1}\al=z$.

We claim that $x<y<z$. If $y\leq x$, then $x=y\al\leq x\al=1$, a contradiction, so $x<y$.
If $z\le x$, then $z=z\al\leq x\al=1$, a contradiction, so $x<z$.
If $z\leq y$, then $z=z\al\leq y\al=x<z$, again a contradiction, so $y<z$.
The claim has been proved.

Since $\lev(x)=\lev(y)=1$, $x<y$, $y\al=x$, and $k\geq1$, we have $\hei(\tr_x)\geq\hei(\tr_y)\geq1$ by Lemma~\ref{lhei1}.
To prove (a), let $s \in V_1(\tr)$ such that $s > z$, and let $s_1 = s\al$. Since $\al$ is both order-preserving and decreasing, 
$z = z\al \leq s\al = s_1 \leq s$. On the other hand, since $s \in V_1(\tr)$ and $v_{k-1} \in V_2(\tr)$, we have $s < v_{k-1}$, and so 
$s_1 = s\al \leq v_{k-1}\al = z$. Hence, $s_1= z$, and so $\hei(\tr_s)\leq\hei(\tr_z)$ by $z=s_1=s\al$ and Lemma~\ref{lhei1}.

To prove (b), let $s\in V_1(\tr)$ with $s<y$. Then, $s\al\leq y\al=x$.
Set $l=\hei(\tr_s)$. Then, there exists a path $s_0\to s_1\to\cdots\to s_l=s\to 1$ in $\tr$. 
If $l=0$, then $\hei(\tr_s)=l=0<\hei(\tr_z)+1$. Suppose $l>0$, and let $r=s_{l-1}\al$. Since $s_{l-1}\in V_2(\tr)$, we have
$s_{l-1}>z$, and so $r=s_{l-1}\al\geq z\al=z$.

Suppose $s\leq x$. Then, $s\al\leq x\al=1$, and so $s\al=1$. 
It then follows by Lemma~\ref{lpath}, $r\geq z$, and $z\in V_1(\tr)$ that $r\in V_1(\tr)$. Hence, by Lemma~\ref{lpath} again,
$s_0\al\to s_1\al\to\cdots\to s_{l-1}\al=r$ is a path in~$\tr_r$. Thus, $l-1\leq\hei(\tr_r)$, and so
$\hei(\tr_s)=l\leq\hei(\tr_r)+1\leq\hei(\tr_z)+1$, where the last inequality follows from~(a) if $r>z$ and is immediate if $r=z$.

Suppose $x<s<y$, and let $s'=s\al$. Then, $s'\leq y\al=x$. Since $s\in V_1(\tr)$ and $\al\in C(\mu)$,
Lemma~\ref{lend} implies that $s'\in V_1(\tr)\cup\{1\}$. If $s'\in V_1(\tr)$, then
$\hei(\tr_{s'})\leq\hei(\tr_z)+1$ (by the foregoing argument), and so
$\hei(\tr_s)\leq\hei(\tr_{s'})\leq\hei(\tr_z)+1$ (by Lemma~\ref{lhei1}). If $s'=1$, then  
$\hei(\tr_s)\leq\hei(\tr_z)+1$ by the same argument as in the case $s\leq x$.
This concludes the proof of (b).

$(\Leftarrow)$ Conversely, suppose there exist $x, y, z \in V_1(\tr)$ such that $x < y < z$,  $\hei(\tr_x)\geq\hei(\tr_y)\geq1$,
and conditions (a) and (b) hold.
Fix paths $\pth_x$ and $\pth_z$ of maximum length in $\tr_x$ and $\tr_z$, respectively. 
Define $\al:\xn\to\xn$ by:

\vskip -5mm
\begin{align*}
\al|_{V(\tr_s)}&=\mbox{$\lamp(\tr_s,\pth_z)$ for every $s\in V_1(\tr)$ with $s<y$ (possible by (b));}\\[-1mm]
\al|_{V(\tr_y)}&=\mbox{$\lam(\tr_y,\pth_x)$ (possible since $\hei(\tr_y)\leq\hei(\tr_x)=\length(\pth_x)$);}\\[-1mm]
\al|_{V(\tppz)}&=\mbox{$\lam(\tppz,\pth_z)$ (possible by the choice of $\pth_z$)};\\[-1mm]
\al|_{V(\tr_s)}&=\mbox{$\lam(\tr_s,\pth_z)$ for every $s\in V_1(\tr)$ with $s>z$ (possible by (a)); }\\[-1mm]
\end{align*}
\vskip -5mm
\noindent and $p\al=p$ for every other $p\in\xn$ (see Notations~\ref{nlam1}, \ref{nlam3}, and~\ref{ntrx}). 
Then, $\al\in\cn$ (by Lemma~\ref{lllp}) and $\al\in C(\mu)$ (by Lemma~\ref{lend}). 

Since $\hei(\tr_y)\geq1$, there exists $w\in V_1(\tr_y)$.
Let $v=w\al$ and $u=v\al$. 
Then, $v\in V_1(\tr_x)$  (by the definition of $\al|_{V(\tr_y)}$ and Lemma~\ref{lpath}) and $u=v\al=z$
(by the definition of $\al|_{V(\tr_s)}$  for $s\in V_1(\tr)$ with $s<y$).
Hence, $\lev(w)=\lev(v)=2>1=\lev(u)$ and $w>v$ (since $\tr$ is level-ordered, $y>x$, $w\in V_1(\tr_y)$, and $v\in V_1(\tr_x)$).
Finally,  $u\al=z\al=z=u$  (by the definition of $\al|_{V(\tppz)}$).
\end{proof}

\begin{defi}\label{dcseq}
Let $\mu\in N(\cn)$ with $\tr=\ga(\mu)$, $x_0,x_1,x_2\in\xn$, and $\al\in C(\mu)$.
We say that $\lan x_0,x_1,x_2\ran$
is a \emph{forbidden $\al$-sequence} in $\tr$ if $\lev(x_0)=\lev(x_1)=\lev(x_2)$,
$x_0>x_1>x_2$, $x_0\al=x_1$, $x_1\al=x_2$, and $x_2\al=x_2$.
\end{defi}

The following lemma describes $\mu\in N(\cn)$ for which there is $\al\in C(\mu)$ such that $w\al=v$, $v\al=u$, and $u\al=u$
for some $w,v,u\in\xn$ such that $v>u$ and \eqref{ecase3} holds,
under the assumption that $\ga(\mu)$ does not have
a forbidden $\al$-sequence.

\begin{lemma}\label{lth3}
	Let $n\geq4$ and $\mu\in N(\cn)$ with $\tr=\ga(\mu)$.
Then, there exists $\al\in C(\mu)$ such that:
\begin{itemize}
\item $\tr$ does not have a forbidden $\al$-sequence; and
\item there are $w,v,u\in\xn$ with $w\al=v$, $v\al=u$, $u\al=u$, $\lev(w)>\lev(v)=\lev(u)>0$, and $w>v>u$
\end{itemize}
if and only if there are $x,y,z\in V_1(\tr)$ 
such that $x<y\leq z$, $\hei(\tr_x)\geq 1$,  and at least one of the following conditions is satisfied:
\begin{itemize} 
  \item[\rm(1)] $y<z$,  there is a path $\pth_x$ in $\tr_x$ with initial vertex of in-degree $0$ and terminal vertex $x$, and 
the following conditions hold:
      \begin{itemize}
                \item[\rm(a)] $\hei(\tzmpx)\leq\hei(\tr_z)+1$ and  $\hei(\tzppx)\leq\hei(\tr_z)+1$;
                \item[\rm(b)] there exists $e\in V_1(\tr)$ such that  $e\geq z$,  and $\hei(\toppx)\leq\hei(\tr_{e})+1$;
		\item[\rm(c)] for every $s\in V_1(\tr)$, if $s>y$, then $\hei(\tr_s)\leq\hei(\tr_y)$;
		\item[\rm(d)] for every $s\in V_1(\tr)$, if $s\leq x$, then $\hei(\tr_s)\leq\hei(\tr_y)+1$;
      \end{itemize} 
  \item[\rm(2)] $y=z$, $\hei(\tr_y)\geq1$, there are paths $\pth_x$ in $\tr_x$ and $\pth_y$, $\qth_y$ in $\tr_y$,
with terminal vertex $x$ for $\pth_x$ and terminal vertex $y$ for $\pth_y$ and $\qth_y$, 
and vertices $f\in V(\pth_x)$, $g\in V(\pth_y)$, and $h\in V(\qth_y)$
such that $\length(\pth_y)\geq\length(\qth_y)$, $\lev(g)=\lev(h)=\lev(\mrp(g,h))+1$, $g<h$, $\dst(f,x)=\dst(h,y)+1$,
$\pth_x$ has maximum length among the paths with terminal vertex $x$ that pass through~$f$,
$\pth_y$ and~$\qth_y$ have maximum length among the paths with terminal vertex $y$ that pass through $g$ and~$h$, respectively,
and the following conditions hold:
     \begin{itemize}
  \item[\rm(i)] $\hei(\tppy)=\length(\pth_y)$;
  \item[\rm(ii)] $\hei(\tppx)\leq\length(\qth_y)+1$;
  \item[\rm(iii)] $\hei(\tmpx)\leq\length(\pth_y)+1$;
  \item[\rm(iv)] $\hei(\tr_f)\leq\length(\qth_h)$, where $\qth_h$ is the path in $\tr$ with initial vertex equal to the initial vertex of $\qth_y$,
and terminal vertex $h$;
  \item[\rm(v)] for every $s\in V_1(\tr)$, if $s>y$, then $\hei(\tr_s)\leq\length(\pth_y)$;  
  \item[\rm(vi)] for every $s\in V_1(\tr)$, if $s<x$, then $\hei(\tr_s)\leq\length(\pth_y)+1$. 
\end{itemize}
\end{itemize}
\end{lemma}

\begin{proof}
$(\Rightarrow)$
Suppose that there exists $\al\in C(\mu)$ that satisfies the given conditions. Consider the unique paths in $\tr$
 from $w$, $v$, and $u$ to $1$: 
$w=w_0\to\cdots\to w_k\to 1$, $v=v_0\to\cdots\to v_m\to 1$, and $u=u_0\to\cdots\to u_m\to 1$. Note that $0\leq m<k$ and 
$w_k,v_m,u_m\in V_1(\tr)$. Set $x=w_k$, $y=u_m$, $z=v_m$.
Let $\pth_x$ be a path
of maximum length in $\tr_x$ that passes through $w$. Note that $w_0=w,\ldots,w_k=x$ lie on $\pth_x$.
Since $\al\mu=\mu\al$, $w\al=v$, and $k\geq m+1=\lev(v)$,
we have $x\al=(w\mu^k)\al=(w\al)\mu^k=v\mu^k=1$.

Moreover, since $\lev(v)=\lev(u)$, $v\al=u$, $u\al=u$, and $w\al=v$,
Lemma~\ref{lpath} implies $z\al=y$, $y\al=y$, and either $w_{k-1}\al=z$ (if $k=m+1$) or $w_{k-1}\al=1$ (if $k>m+1$). However, since 
$w_{k-1}\in V_2(\tr)$ and $z\in V_1(\tr)$, we have $w_{k-1}>z$. This yields $w_{k-1}\al\geq z\al=y>1$,
so $w_{k-1}\al=z$.

If $y\leq x$, then $y=y\al\leq x\al=1$, a contradiction, so $x<y$. Since $\al$ is decreasing, $y=z\al\leq z$. Thus, $x<y\leq z$.
Furthermore, it is clear that $\hei(\tr_x)\geq k\geq 1$. We consider two possible cases.
\vskip 1mm
\noindent {\bf Case 1.} $y<z$. 
\vskip 1mm
In this case, we must prove that (1) is satisfied. To prove (a), first note that  $\hei(\tzmpx)\leq\hei(\tr_{w_{k-1}})+1\leq\hei(\tr_z)+1$, where the last inequality follows from $w_{k-1}\al=z$
and Lemma~\ref{lhei1}. Similarly,  $\hei(\tzppx)\leq\hei(\tr_z)+1$. 

To prove (b), let $e\in V_1(\tr)$ be such that  $e\geq z$ and $\hei(\tr_e)\geq\hei(\tr_t)$ for every $t\in V_1(\tr)$ with $t\geq z$.
Let $l=\hei(\toppx)$. If $l=0$, then $\hei(\toppx)=l=0\leq\hei(\tr_e)+1$.  Suppose $l>0$. Consider a path $r_0\to\cdots\to r_l=x$ in
 $\toppx$. Then,  $r_{l-1}>w_{k-1}$
(see Notation~\ref{ntrx}), and so $r_{l-1}\al\geq w_{k-1}\al=z>1$.
Moreover, $(r_{l-1}\al)\mu=(r_{l-1}\mu)\al=x\al=1$,  and so $r_{l-1}\al>1$ implies $r_{l-1}\al\in V_1(\tr)$.
Since $r_0\to\cdots\to r_{l-1}$ is a path of length $l-1$ in $\tr_{r_{l-1}}$, we have $l-1\leq\hei(\tr_{r_{l-1}})$.
Thus, by Lemma~\ref{lhei1} and the definition of $e$, $\hei(\toppx)=l\leq\hei(\tr_{r_{l-1}})+1\leq\hei(\tr_{r_{l-1}\al})+1\leq\hei(\tr_e)+1$.

To prove (c), let $s\in V_1(\tr)$ with $s>y$. If $y<s\leq z$, then $y=y\al\leq s\al\leq z\al=y$, which forces $s\al=y$. Consequently, $\hei(\tr_s)\leq\hei(\tr_y)$ by Lemma~\ref{lhei1}. Now, suppose that $s>z$. Since $\al$ is order-preserving, $s\al\geq z\al=y$. Moreover, since $s\in V_1(\tr)$ and $w_{k-1}\in V_2(\tr)$, it follows that $s<w_{k-1}$, which implies $s\al\leq w_{k-1}\al=z$. 
Thus, by the foregoing argument and Lemma~\ref{lhei1}, $\hei(\tr_s)\leq\hei(\tr_{s\al})\leq\hei(\tr_y)$.

To prove (d), let $s\in V_1(\tr)$ with $s\leq x$. If $\hei(\tr_s)=0$, then the inequality $\hei(\tr_s)\leq\hei(\tr_y)+1$ holds trivially. 
Assume that $l=\hei(\tr_s)>0$. Then, there exists a path $s_0\to\cdots\to s_{l-1}\to s_l=s\to1$ in $\tr$. Set $r=s_{l-1}\al$. 
Since $s\leq x$ and $x\al=1$, we have $s\al=1$. Since $s_{l-1}\in V_2(\tr)$ and $z\in V_1(\tr)$, we have $s_{l-1}>z$, which gives $r=s_{l-1}\al\geq z\al=y$. Moreover,
$r\mu=(s_{l-1}\al)\mu=(s_{l-1}\mu)\al=s\al=1$. Thus, $r\in V_1(\tr)$, and so, by Lemma~\ref{lpath}, $s_0\al\to\cdots\to s_{l-1}\al=r$ is a path in $\tr_r$, which implies $l-1\leq\hei(\tr_r)$,
and so $l\leq\hei(\tr_r)+1$.

Suppose  $s_{l-1}\leq w_{k-1}$.  Then, $y\leq r=s_{l-1}\al\leq w_{k-1}\al=z$,
which implies $y\leq r\leq z$. Thus, $y=y\al\leq r\al\leq z\al=y$, and so $r\al=y$. Hence, 
$\hei(\tr_s)=l\leq\hei(\tr_r)+1\leq\hei(\tr_y)+1$, where the last inequality follows from Lemma~\ref{lhei1}
and $r\al=y$.

Suppose $s_{l-1}>w_{k-1}$. Then, $r=s_{l-1}\al\geq w_{k-1}\al=z$,
and so $\hei(\tr_s)=l\leq\hei(\tr_r)+1\leq\hei(\tr_y)+1$, where the last inequality follows from (c)
and $r\geq z>y$.
\vskip 1mm
\noindent {\bf Case 2.} $y=z$. 
\vskip 1mm
In this case, we must prove that (2) is satisfied. We may assume that there are no vertices $w^*$, $v^*$, $u^*$, $x^*$, $y^*$,
and $z^*$ of $\tr$
such that
$w^*\al=v^*$, $v^*\al=u^*$, $u^*\al=u^*$, $\lev(w^*)>\lev(v^*)=\lev(u^*)>0$, $w^*>v^*>u^*$, $x^*,y^*,z^*\in V_1(\tr)$,
$x^*<y^*<z^*$, $w^*\in V(\tr_{x^*})$, $u^*\in V(\tr_{y^*})$, and $v^*\in V(\tr_{z^*})$.
Indeed, if such vertices exist, we can apply the foregoing argument to them to prove that (1) is satisfied.

Now, we may further assume that $w$ is the vertex in $\tr_x$ at the smallest level of $\tr_x$ for which vertices $v,u\in V(\tr_y)$ exist
such that $\lev(w)>\lev(v)=\lev(u)>0$, $v>u$, $w\al=v$, $v\al=u$, and $u\al=u$. Let $f$ be the largest vertex
in $\tr_x$ with $\lev(f)=\lev(w)$ for which there exist vertices $g,h\in V(\tr_y)$ such that $\lev(g)=\lev(h)=\lev(v)=\lev(u)$, $h>g$,
$f\al=h$, $h\al=g$, and $g\al=g$.

Let $\pth_x$ be a path of maximum length among the paths in $\tr_x$ that pass through $f$, and let $\pth_y,\qth_y$ be paths of maximum length among the paths in $\tr_y$ that pass through $g$ and $h$, respectively.
By Lemma~\ref{lhei1} and $h\al=g$, we have $\hei(\tr_g)\geq\hei(\tr_h)$. Since $\lev(g)=\lev(h)$, we also have
$\dst(g,y)=\dst(h,y)$. Thus, the maximality of $\pth_y$ and $\qth_y$ implies 
\[
\length(\pth_y)=\hei(\tr_g)+\dst(g,y)\geq \hei(\tr_h)+\dst(h,y)=\length(\qth_y).
\]
The minimality of $\lev(f) (=\lev(w))$ implies $\lev(g)=\lev(h)=\lev(\mrp(g,h))+1$.
 (Indeed, if $g\mu\ne h\mu$, then $f\mu > h\mu > g\mu$, $f\mu\in V(\tr_x)$, $(f\mu)\al=(f\al)\mu = h\mu$,
$ (h\mu)\al = (h\al)\mu = g\mu$, and $(g\mu)\al =(g\al)\mu= g\mu$, contradicting
the choice of $w$.)
Let $\pth^f=f_0\to\cdots\to f_k\to1$ be the unique path in $\tr$ from $f$ to $1$. Note that $k\geq1$ and $f_k=x$.
As in the proof of $w_{k-1}\al=z$ for  $w=w_0\to\ldots w_k\to1$, we obtain  $f_{k-1}\al=y$. Hence,  by Lemma~\ref{lpath},
$\dst(f,x)=\dst(f,f_{k-1}) +1=\dst(f\al,f_{k-1}\al) +1=\dst(h,y)+1$.

To prove (i), select $t\in V(\tppy)$ such that $\lev(t)=\lev(g)$ and $\hei(\tppy)=\hei(\tr_t)+\dst(t,y)$,
and note that $t\geq g$. Let $t_1=t\al$.
Then, $t_1=t\al\geq g\al=g$. Since $\lev(t)=\lev(g)<\lev(f)$, we have $t<f$, and so $t_1=t\al\leq f\al=h$.
It follows that $g\leq  t_1\leq h$, which implies $\lev(t_1)=\lev(t)=\lev(g)$ and $g=g\al\leq  t_1\al\leq h\al=g$. Thus, $t_1\al=g$ and
$\dst(t,y)=\dst(t_1,y)=\dst(g,y)$.   Hence, by Lemma~\ref{lhei1} and the selection of $\pth_y$,
\[
\hei(\tppy)=\hei(\tr_t)+\dst(t,y)\leq\hei(\tr_{t_1})+\dst(t_1,y)\leq\hei(\tr_g)+\dst(g,y)=\length(\pth_y).
\]
On the other hand, $\length(\pth_y)\leq\hei(\tppy)$ by the definition of $\tppy$. Thus, $\hei(\tppy)=\length(\pth_y)$.

To prove (ii), we first show that $\length(\pth_x)\leq\length(\qth_y)+1$. We have $f\al=h$ and $\dst(f,x)=\dst(h,y)+1$, and so,
by Lemma~\ref{lhei1},
\[
\length(\pth_x)=\hei(\tr_f)+\dst(f,x)\leq\hei(\tr_h)+\dst(h,y)+1=\length(\qth_y)+1.
\]
Now, select $t\in V(\tppx)$ such that $\lev(t)=\lev(f)$ and $\hei(\tppx)=\hei(\tr_t)+\dst(t,x)$,
and note that $t\geq f$. If $t=f$, then $\hei(\tppx)=\length(\pth_x)\leq\length(\qth_y)+1$.

Suppose $t>f$, and let $t_1=t\al$ and $t_2=t_1\al$.
Then, since $t>f$, we have $t_1=t\al\geq f\al=h$, and so $t_2=t_1\al\geq h\al=g$.
Since $t\in V(\tppx)$, there is a path from $t$ to $x$. Thus, since $\lev(x\al)=\lev(1)=0$, we have $\lev(t_1)=\lev(t\al)<\lev(t)=\lev(f)$
by Lemma~\ref{lpath}.
Hence, $t_1<f$, and so $t_2=t_1\al\leq f\al=h$. We have proved that $g\leq t_2\leq h$,
which implies $g=g\al\leq t_2\al\leq h\al=g$, so $t_2\al=g$, and
$\lev(t_1\al)=\lev(h)$. Further, $\lev(t_1)\leq\lev(h)$ since $\lev(t_1)<\lev(f)$ and $\lev(f)=\lev(h)+1$. 
But we cannot have $\lev(t_1)<\lev(h)$ since $\lev(t_1\al)=h$ and $\al$ is decreasing. Thus, $\lev(t_1)=\lev(h)=\lev(t_2)$.

Suppose to the contrary that $t_1>h$. If $t_2$ were greater than $g$, then $\lan t_1,t_2,g\ran$ would be a forbidden $\al$-sequence
in $\tr$. Thus, $t_2=g$. We claim that $t_1\in V(\tr_y)$. Suppose not. Let $z^*$ be the (unique) vertex in $V_1(\tr)$ on the path from
$t_1$ to $1$ in $\tr$. Since $t_1>h$ and $t_1\notin V(\tr_y)$, we have $z^*>y$. Then,
$t\al=t_1$, $t_1\al=t_2=g$, $g\al=g$, $\lev(t)>\lev(t_1)=\lev(g)>0$, $t>t_1>g$, $x,y,z^*\in V_1(\tr)$,
$x<y<z^*$, $t\in V(\tr_x)$, $g\in V(\tr_y)$, and $t_1\in V(\tr_{z^*})$. But this contradicts
the assumption made at the beginning of Case~2. The claim has been proved.
We now have $t>f$, $\lev(t)=\lev(f)=\lev(w)$, $g,t_1\in V(\tr_y)$, $\lev(t_2)=\lev(t_1)=\lev(v)=\lev(u)$, $t_1>g$,
$t\al=t_1$,  $t_1\al=g$, and $g\al=g$. This is a contradiction by the maximum property of $f$.
Hence, $t_1=h$, and so
\[
\hei(\tppx)=\hei(\tr_t)+\dst(t,x)\leq\hei(\tr_{t_1})+\dst(t_1,y)+1=\hei(\tr_h)+\dst(h,y)+1=\length(\qth_y)+1.
\]

To prove (iii), select $t\in V(\tmpx)$ such that $\lev(t)=\lev(f)$ and $\hei(\tmpx)=\hei(\tr_t)+\dst(t,x)$,
and note that $t\leq f$. Let $t_1=t\al$. Then, $t_1=t\al\leq f\al=h$ and $t_1=t\al\geq g\al=g$ (since $\lev(t)=\lev(f)>\lev(g)$, and so
$t>g$). Hence, $g\leq t_1\leq h$, and so $t_1\in V(\tppy)$.
By (i), $\hei(\tppy)\leq\length(\pth_y)$, and so, by Lemma~\ref{lhei1} and $t_1=t\al$,
\[
\hei(\tmpx)=\hei(\tr_t)+\dst(t,x)\leq\hei(\tr_{t_1})+\dst(t_1,y)+1\leq\hei(\tppy)+1\leq\length(\pth_y)+1.
\]

To prove (iv), let $\qth_h$ be as described. Since $f\al=h$, Lemma~\ref{lhei1} and the choice of $\qth_y$ give
$\hei(\tr_f)\leq\hei(\tr_h)=\length(\qth_h)$.

To prove (v), let $s\in V_1(\tr)$ with $s>y$. If $\hei(\tr_s)<\length(\pth_y)$, then there is nothing to prove.
Otherwise,  there is $t\in V(\tr_s)$
such that $\lev(t)=\lev(g)$ and $\hei(\tr_s)=\hei(\tr_t)+\dst(t,s)$. Let $t_1=t\al$.
Since $s>y$, we have $t>g$. As in the proof of (i), we obtain $g\leq t_1\leq h$. Hence, $\lev(t_1)=\lev(t)=\lev(g)$ and $t_1\in V(\tppy)$, so $\dst(t,s)=\dst(t_1,y)$. Thus, by
Lemma~\ref{lhei1} and~(i),
\[
\hei(\tr_s)=\hei(\tr_t)+\dst(t,s)\leq\hei(\tr_{t_1})+\dst(t_1,y)\leq\hei(\tppy)\leq\length(\pth_y).
\]

To prove (vi), let $s\in V_1(\tr)$ with $s<x$. If $\hei(\tr_s)<\length(\pth_x)$, then
\[
\hei(\tr_s)<\length(\pth_x)\leq\length(\qth_y)+1\leq\length(\pth_y)+1,
\]
where the second inequality was proved in (ii). Otherwise,
select $t\in V(\tr_s)$ such that $\lev(t)=\lev(f)$ and $\hei(\tr_s)=\hei(\tr_t)+\dst(t,s)$,
and note that $t<f$ (since $s<x$). The rest of the proof
goes exactly as the corresponding part of the proof of (iii) except that we replace ``$\hei(\tmpx)$'' with 
``$\hei(\tr_s)$.''

$(\Leftarrow)$ Conversely, suppose there exist $x, y, z \in V_1(\tr)$  
such that $x < y\leq z$ and  $\hei(\tr_x)\geq1$.

Suppose (1) is satisfied.
Fix paths $\pth_y$, $\pth_z$, and $\pth_e$ of maximum length in $\tr_y$, $\tr_z$, and $\tr_e$, respectively.
We assume that  $\pth_z=\pth_e$ if $z=e$.
Define $\al:\xn\to\xn$ by: 
\vskip -5mm
\begin{align*}
\al|_{V(\tr_s)}&=\mbox{$\lamp(\tr_s,\pth_y)$ for every $s\in V_1(\tr)$ with $s<x$ (possible by (d));}\\[-1mm]
\al|_{V(\tompx)}&=\mbox{$\lamp(\tompx,\pth_y)$ (possible by (d));}\\[-1mm]
\al|_{V(\tzmpx)}&=\mbox{$\lamp(\tzmpx,\pth_z)$ (possible by (a));}\\[-1mm]
\al|_{V(\tzppx)}&=\mbox{$\lamp(\tzppx,\pth_z)$ (possible by (a));}\\[-1mm]
\al|_{V(\toppx)}&=\mbox{$\lamp(\toppx,\pth_e)$ (possible by (b));}\\[-1mm]
\al|_{V(\tppy)}&=\lam(\tppy,\pth_y);\\[-1mm]
\al|_{V(\tr_s)}&=\mbox{$\lam(\tr_s,\pth_y)$ for every $s\in V_1(\tr)$ with $s>y$ (possible by (c));}\\[-1mm]
\end{align*}
\vskip -5mm
\noindent and $p\al=p$ for every other $p\in\xn$ (see Notations~\ref{nlam1},  \ref{nlam3}, and~\ref{ntrx}). 
Then, $\al\in\cn$ (by Lemma~\ref{lllp}) and $\al\in C(\mu)$ (by Lemma~\ref{lend}). 

Since $\hei(\tr_x)\geq1$ and $\pth_x$ has initial vertex of in-degree $0$, $\length(\pth_x)\geq1$, and so
there exists $w\in V_1(\tr_x)$ such that $w$ lies on $\pth_x$.
Let $v=w\al$ and $u=v\al$. Then, $v=z$ (by the definition of $\al|_{V(\tzmpx)}$) and $u=v\al=z\al=y$ (by the definition of 
$\al|_{V(\tr_s)}$ for $s\in V_1(\tr)$ with $s>y$). 
Hence, $\lev(w)=2>\lev(v)=1=\lev(u)$ and $v=z>y=u$. Finally, $u\al=y\al=y=u$ (by the definition of $\al|_{V(\tppy)}$).

Suppose $x_0,x_1\in\xn$ are such that $\lev(x_0)=\lev(x_1)$, $x_0>x_1$, and $x_0\al=x_1$. By the definition of $\al$,
$x_0\in V(\tppy)$ or $x_0\in V(\tr_s)$ for some $s\in V_1(\tr)$ with $s>y$, and $x_1\in V(\pth_y)$.
Since every vertex lying on $\pth_y$ is fixed by $\al$, $x_1\al=x_1$. It follows that $\tr$ does not have a forbidden $\al$-sequence
(see Definition~\ref{dcseq}).

Suppose (2) is satisfied. Let $f'=f\mu$ and note that $f'$ lies on $\pth_x$.
Split the path $\pth_x$ into two paths: $\pth_{\!f}$ and $\pth^{(f'\!\!,x)}$, where $\pth_{\!f}$ is the path in 
$\tr$ with initial vertex equal to the initial vertex of $\pth_x$,
and terminal vertex $f$, and $\pth^{(f'\!\!,x)}$ is the path from $f'$ to $x$.
Since $f'$ need not have the in-degree $0$, we do not apply Notation~\ref{ntrx} to $\pth^{(f'\!\!,x)}$ directly; instead, we set
$\tmpfpr=\tmpx\setminus\tr_{\!f}$ and $\tppfpr=\tppx\setminus\tr_{\!f}$.
Define $\al:\xn\to\xn$ by: 
\vskip -5mm
\begin{align*}
\al|_{V(\tr_s)}&=\lamp(\tr_s,\pth_y) \mbox{ for every $s\in V_1(\tr)$ with $s<x$ (possible by (vi));}\\[-1mm]
\al|_{V(\tmpf)}&=\mbox{$\lam(\tmpf,\qth_h)$ (possible by (iv));}\\[-1mm]
\al|_{V(\tmpfpr)}&=\mbox{$\lamp(\tmpfpr,\pth_y)$ (possible by (iii));}\\[-1mm]
\al|_{V(\tppf)}&=\mbox{$\lam(\tppf,\qth_h)$ (possible by (iv));}\\[-1mm]
\al|_{V(\tppfpr)}&=\mbox{$\lamp(\tppfpr,\qth_y)$ (possible by (ii));}\\[-1mm]
\al|_{V(\tppy)}&=\mbox{$\lam(\tppy,\pth_y)$ (possible by (i));}\\[-1mm]
\al|_{V(\tr_s)}&=\mbox{$\lam(\tr_s,\pth_y)$ for every $s\in V_1(\tr)$ with $s>y$ (possible by (v));}\\[-1mm]
\end{align*}
\vskip -5mm
\noindent and $p\al=p$ for every other $p\in\xn$ (see Notations~\ref{nlam1},  \ref{nlam3}, and~\ref{ntrx}). Then, $\al\in\cn$ (by Lemma~\ref{lllp}) 
and $\al\in C(\mu)$ (by Lemma~\ref{lend} and $\lev(g)=\lev(h)=\lev(\mrp(g,h))+1$). We note that the hypothesis $\lev(g)=\lev(h)=\lev(\mrp(g,h))+1$
guarantees that $\al$ is a digraph homomorphism from $\tr$ to $\tr$  splitting $\pth_x$. Indeed, let $r=\mrp(g,h)$.
The edge in $\tr$ that links $\pth_{\!f}$ and $\pth^{(f'\!\!,x)}$ is $f\to f'$. We have $f\al=h$ (by the definition of $\al|_{V(\tppf)}$)
and $f'\al=r$ (by the definition of $\al|_{V(\tmpfpr)}$).
Since $\lev(h)=\lev(r)+1$ and $r=\mrp(g,h)$, we have $h\to r$ (see Definition~\ref{dmerge}), that is, $f\al\to f'\al$, as required.
Also, since  $\lev(g)=\lev(h)=\lev(\mrp(g,h))+1$, the definition of $\al|_{V(\tmpfpr)}$ is consistent with the definition of $\al|_{V(\tppfpr)}$.

Now, let $w=f$, $v=h$, and $u=g$. Then, $\lev(v)=\lev(h)=\lev(g)=\lev(u)$, $\lev(w)=\lev(f)>\lev(h)=\lev(v)$ 
(since $\dst(w,x)=\dst(f,x)=\dst(h,y)+1=\dst(v,y)+1$), and $v=h>g=u$. Further,
$w\al=f\al=h=v$  (by the definition of $\al|_{V(\tmpf)}$), $v\al=h\al=g=u$  (by the definition of $\al|_{V(\tppy)}$),
and $u\al=g\al=g=u$  (by the definition of $\al|_{V(\tppy)}$).

Finally, the proof that $\tr$ does not have a forbidden $\al$-sequence
is identical to the one we have in the case when (1) is satisfied.
\end{proof}

We are now ready to state and prove the main theorem of this subsection.

\begin{theorem}\label{tiv}
Let $n\geq4$ and $\mu\in N(\cn)$ with $\tr=\ga(\mu)$.  Then, $C(\mu)\ne C_{nil}(\mu)\cup C_{id}(\mu)$
if and only if at least one of the following conditions is satisfied:
\begin{itemize}
\item[\rm(1)] there are $x,y,z\in V_1(\tr)$ such that $x<y<z$, and both {\rm(a)} and {\rm(b)} of {\rm Lemma~\ref{lwv1u}} hold;
\item[\rm(2)] {\rm(1)} or {\rm(2)} of {\rm Lemma~\ref{lth4}} holds;
\item[\rm(3)] there are $x,y,z\in V_1(\tr)$ such that $x<y<z$, $\hei(\tr_x)\geq\hei(\tr_y)\geq1$, and both {\rm(a)} and {\rm(b)} of {\rm Lemma~\ref{lth2}} hold;
\item[\rm(4)] there are $x,y,z\in V_1(\tr)$ such that $x<y\leq z$, $\hei(\tr_x)\geq1$, and either {\rm(1)} or {\rm(2)} of {\rm Lemma~\ref{lth3}} holds.
\end{itemize}
\end{theorem}
\begin{proof}
($\imp$)
Suppose $C(\mu)\ne C_{nil}(\mu)\cup C_{id}(\mu)$. Then, there exists $\al\in C(\mu)$
that satisfies (1) or (2) of Lemma~\ref{l1oru}. If $\al$ satisfies (1) of Lemma~\ref{l1oru}, then (1) holds by Lemma~\ref{lwv1u}.

Suppose $\al$ satisfies (2) of Lemma~\ref{l1oru}, that is, $w\al=v$, $v\al=u$, and $u\al=u$ for some $w,v,u\in\xn$ with $w>v>u>1$.
Then, \eqref{ecase1}, \eqref{ecase2}, or \eqref{ecase3} holds for $w,v,u$. (By Lemma~\ref{lth1}, \eqref{ecase4} is impossible.)
If \eqref{ecase1} holds, then (2) holds by Lemma~\ref{lth4}. If \eqref{ecase2} holds, then (3) holds by Lemma~\ref{lth2}.

Assume \eqref{ecase3} holds. Suppose $\tr$ has a forbidden $\al$-sequence, say $\lan x_0,x_1,x_2\ran$. Set $w=x_0$, $v=x_1$, and $u=x_2$. Then, by Definition~\ref{dcseq},
$\lev(w)=\lev(v)=\lev(u)$, $w>v>u$, $w\al=v$, $v\al=u$, and $u\al=u$. 
Hence, (2) holds by Lemma~\ref{lth4}. If $\tr$ does not have a forbidden $\al$-sequence, then (4) holds by
Lemma~\ref{lth3}.

$(\Leftarrow)$ Suppose that (1), (2), (3), or (4) is satisfied. Then, by Lemmas~\ref{lwv1u}, \ref{lth4}, \ref{lth2}, and~\ref{lth3}, respectively,
there exists $\al\in C(\mu)$ that satisfies (1) or (2) of Lemma~\ref{l1oru}. Hence, $C(\mu)\ne C_{nil}(\mu)\cup C_{id}(\mu)$ by Lemma~\ref{l1oru}.
\end{proof}

From Theorem~\ref{tiv} we can deduce various corollaries that  give sufficient conditions for a
nilpotent $\mu\in\cn$ to satisfy $C(\mu)=C_{nil}(\mu)\cup C_{id}(\mu)$.
We will present one such corollary.

\begin{cor}\label{civ}
Let $n\geq4$ and $\mu\in N(\cn)$ with $\tr=\ga(\mu)$.  Suppose the following conditions are satisfied:
\begin{itemize}
  \item[\rm(a)] for every $x\in\xn$, $|V_1(\tr_x)|\leq2$;
  \item[\rm(b)] for every $x\in\xn$, if $V_1(\tr_x)=\{a,b\}$ with $a<b$, then $\hei(\tr_a)<\hei(\tr_b)$.
\end{itemize}
Then, $C(\mu)=C_{nil}(\mu)\cup C_{id}(\mu)$.
\end{cor}
\begin{proof}
First, by (a), neither (1) nor (3) of Theorem~\ref{tiv} is satisfied. Second,
(1) of Lemma~\ref{lth4} does not hold by (a), and (2) of Lemma~\ref{lth4} does not hold by (a) and (b). Hence,
(2) of Theorem~\ref{tiv} is not satisfied.

It remains to exclude (4) of Theorem~\ref{tiv}.
Condition~(1) of Lemma~\ref{lth3} is impossible by (a).
Suppose that condition~(2) of that lemma holds, and let
$\pth_y$, $\qth_y$, $g$, and $h$ be as in that condition.
Set $p=\mrp(g,h)$. Since
$\lev(g)=\lev(h)=\lev(p)+1$, condition (a) gives
$V_1(\tr_p)=\{g,h\}$. As $g<h$, condition (b) yields
$\hei(\tr_g)<\hei(\tr_h)$.
By the maximality of $\pth_y$ and $\qth_y$ and the equality
$\dst(g,y)=\dst(h,y)$, we obtain
\[
\length(\pth_y)
=\hei(\tr_g)+\dst(g,y)
<\hei(\tr_h)+\dst(h,y)
=\length(\qth_y),
\]
contradicting $\length(\pth_y)\geq\length(\qth_y)$.
Thus, (4) of Theorem~\ref{tiv} is also impossible.
Therefore, $C(\mu)=C_{nil}(\mu)\cup C_{id}(\mu)$ by Theorem~\ref{tiv}.
\end{proof}

\begin{example}\label{e1civ}
Consider  
\begin{equation*}
\mu=\left(\begin{smallmatrix}
	1 & 2 & 3 & 4 & 5 & 6 & 7 & 8 & 9 & 10 & 11 & 12 & 13 & 14 & 15 & 16 & 17 & 18 & 19 & 20 & 21 & 22 & 23 & 24\\
	1 & 1 & 1 & 2 & 2 & 3 & 3 & 4 & 4 &  5 &5 &  6 &   6 &   7 &   7 &   9 &   11 & 13 & 14 & 15 & 17 & 20 & 20 & 23
	\end{smallmatrix}\right) \in N(\mathcal{C}_{24}).
\end{equation*}
By looking at $\tr=\ga(\mu)$ (see Figure~\ref{fciv}), we can see immediately that it satisfies (a) and (b)
of Corollary~\ref{civ}. Hence, $C(\mu)=C_{nil}(\mu)\cup C_{id}(\mu)$.
\end{example}

\begin{figure}[htbp]
	\centering
	\begin{tikzpicture}[
	level distance=1.0cm,
	level 1/.style={sibling distance=6.0cm},
	level 2/.style={sibling distance=3.0cm},
	level 3/.style={sibling distance=1.5cm},
	level 4/.style={sibling distance=1.2cm},
	level 5/.style={sibling distance=1.2cm},
	level 6/.style={sibling distance=1.2cm},
	every node/.style={circle, fill=black, inner sep=1.2pt}
	]
	\node[label=above:1] {}
	child { node[label=above left:2] {}
		child { node[label=above left:4] {}
			child { node[label=below:8] {} }
			child { node[label=right:9] {}
				child { node[label=left:16] {} }
			}
		}
		child { node[label=above right:5] {}
			child { node[label=below:10] {} }
			child { node[label=right:11] {}
				child { node[label=right:17] {}
					child { node[label=below:21] {} }
				}
			}
		}
	}
	child { node[label=above right:3] {}
		child { node[label=above left:6] {}
			child { node[label=below:12] {} }
			child { node[label=right:13] {}
				child { node[label=right:18] {} }
			}
		}
		child { node[label=above right:7] {}
			child { node[label=left:14] {}
				child { node[label=right:19] {} }
			}
			child { node[label=right:15] {}
				child { node[label=right:20] {}
					child { node[label=left:22] {} }
					child { node[label=right:23] {}
						child { node[label=right:24] {} }
					}
				}
			}
		}
	};
	\end{tikzpicture}
	\caption{The tree $\tr=\ga(\mu)$ for the nilpotent $\mu$ in Example~\ref{e1civ}.}
	\label{fciv}
\end{figure}

\begin{example}\label{e2civ}
	Consider  $\mu=\left( \begin{smallmatrix}
	1 & 2 & 3 & 4 & 5 & 6 & 7 & 8 & 9 \\
	1 & 1 & 1 & 2 & 2 & 3 & 5 & 6 & 6 
	\end{smallmatrix}\right) \in N(\mathcal{C}_{9})$ (see  Figure~\ref{fcentmu} for
	the tree $\tr=\ga(\mu)$). We already know that $C(\mu)=C_{nil}(\mu)\cup C_{id}(\mu)$ (see Examples~\ref{ecentmu}, \ref{ecentnilmu}, and~\ref{ecentidmu}). We can also establish this equality using Theorem~\ref{tiv},
without enumerating the elements of $C(\mu)$. Corollary~\ref{civ} does not apply because
$\hei(\tr_2)=\hei(\tr_3)=2$, contrary to condition (b).

Since $|V_1(\tr_x)|\leq2$ for every $x\in\xn$, conditions
(1) and (3) of Theorem~\ref{tiv}, condition (1) of Lemma~\ref{lth4}, and condition (1) of Lemma~\ref{lth3}
are all impossible.
	
	Suppose to the contrary that (2) of Lemma~\ref{lth4} holds, and let $x,y,a,b,c$ be as in (2) of Lemma~\ref{lth4}.
Then, $x=2$, $y=3$, $a=4$, $b=5$, and $c=6$. By (2)(b) of Lemma~\ref{lth4},
$\hei(\tr_a)\geq\hei(\tr_b)$, contradicting $\hei(\tr_4)=0<1=\hei(\tr_5)$.

	Suppose to the contrary that (2) of Lemma~\ref{lth3} holds, and let $x,y,f,g,h$ be as in (2) of Lemma~\ref{lth3}.
Necessarily, $x=2$, $y=3$, $g=8$, and $h=9$. Condition~(2) of Lemma~\ref{lth3} therefore gives
$
\dst(f,2)=\dst(9,3)+1=3.$
However, $f\in V(\tr_2)$ and $\hei(\tr_2)=2$,
so $\dst(f,2)\leq2$, a contradiction.

	We have proved that none of (1), (2), (3), and (4) of Theorem~\ref{tiv} is satisfied. Therefore, $C(\mu)=C_{nil}(\mu)\cup C_{id}(\mu)$.
\end{example}

\section{Centralizers of general elements of $\cn$}\label{sgen}
\setcounter{equation}{0}
In this section, we study the centralizers of arbitrary
elements $\bt\in\cn$ in terms of the ordered tree components of
$\ga(\bt)$, extending and using results for
idempotents in Section~\ref{side} and nilpotents in Section~\ref{snil}.

The following theorem characterizes the centralizer 
$C(\bt)=\{\al\in\cn:\al\bt=\bt\al\}$ of an arbitrary element $\bt\in\cn$. 
It generalizes Theorems~\ref{tcide} and~\ref{tcnilg}.

\begin{theorem}\label{tcgen}
Let $\bt\in\cn$ with $\fix(\bt)=\{e_1<\cdots<e_m\}$ and $\ga(\bt)=\{\tr_1<\cdots<\tr_m\}$
as in {\rm Notation~\ref{ntrcn}}, and let $\al\in\cn$.
Then, $\al\in C(\bt)$ if and only if the following conditions are satisfied:
\begin{itemize} 
  \item [\rm(1)] $(\im(\bt))\al\subseteq \im(\bt)$ and $(\fix(\bt))\al\subseteq\fix(\bt)$;
  \item [\rm(2)] if $\{e_{i_1}<\cdots<e_{i_t}\}=(\fix(\bt))\al$, $\{j_1<\cdots<j_t\}$ is such that $j_r=\min\{j\in[m]:e_j\al=e_{i_r}\}$,
and $j_{t+1}=m+1$, then for all $j\in[m]$ and $r\in[t]$,
\begin{itemize} 
    \item [\rm(a)]  if $j=j_{r+1}-1$, then for all $y\in V(\tr_j)\cap\im(\bt)$ and $x\in\xn$, 
if $x=y\al$, then $x\in V(\tr_{i_r})\cap\im(\bt)$ and $\al|_{J}\in\cnst(J,I)$, where $J=y\bt^{-1}$ and $I=x\bt^{-1}$
 (see {\rm Notation~\ref{ncnstij}});
  \item [\rm(b)] if $j_r\leq j<j_{r+1}-1$, then $(V(\tr_j))\al=\{e_{i_r}\}$.
\end{itemize}
\end{itemize} 
\end{theorem}
\begin{proof}
Let $\al\in C(\bt)$. Suppose $x \in \im(\bt)$. Then, $x=y\bt$ for some $y\in\xn$. Thus, since $\al\bt=\bt\al$,
$x\al=(y\bt)\al=(y\al)\bt\in\im(\bt)$.
Suppose $e \in \fix(\bt)$. Then, $(e\al)\bt=(e\bt)\al=e\al$, and so, $e\al\in\fix(\bt)$. We have proved~(1).

To prove (2), let $\{e_{i_1}<\cdots<e_{i_t}\}$, $\{j_1<\cdots<j_t\}$, and $j_{t+1}$ be as in~(2), 
and let $j\in[m]$ and $r\in[t]$. By the definition of $j_r$ and the fact that $\al$
is order-preserving,
\begin{equation}\label{tcgene1}
e_j\al=e_{i_r}\miff j_r\leq j\leq j_{r+1}-1.
\end{equation}
Suppose that $j=j_{r+1}-1$, $y\in V(\tr_j)\cap\im(\bt)$, and $x=y\al$.
Then, $e_{j}\al = e_{i_r}$ by \eqref{tcgene1}.  Since $y\in V(\tr_j)\cap\im(\bt)$, there exist an integer $k\geq0$
and a vertex $z\in V(\tr_j)$ such that $y\bt^k=e_j$ and $z\bt=y$. Thus, $x\bt^k=(y\al)\bt^k=(y\bt^k)\al=e_j\al=e_{i_r}$,
which implies $x\in V(\tr_{i_r})$. Further, $x=y\al=(z\bt)\al=(z\al)\bt$, and so $x\in V(\tr_{i_r})\cap\im(\bt)$.
Let $J=y\bt^{-1}$, $I=x\bt^{-1}$, and $w\in J$. Then, $w\to y$ in $\ga(\bt)$, and so $w\al\to y\al=x$ in $\ga(\bt)$
by Lemma~\ref{lend}. Hence, $(w\al)\bt=x$, and so $w\al\in I$. We have proved that
$\al$ maps $J$ to $I$, that is, $\al|_J: J \to I$. The sets $I$ and $J$ are nonempty intervals by
Proposition~\ref{paincn}. Further, $I\leq J$ because $x=y\al\leq y$ and $\bt$ is
order-preserving.
Since $\al$ is decreasing and order-preserving, $\al|_J \in \cnst(J,I)$.
This concludes the proof of~2(a).

To prove 2(b), suppose that $j_r\leq j<j_{r+1}-1$, and let
$z\in V(\tr_j)$. By the ordering of the components,
$e_{j_r}\leq z<e_{j+1}\leq e_{j_{r+1}-1}$.
Equation~\eqref{tcgene1} and order preservation of $\al$ give
$e_{i_r}=e_{j_r}\al\leq z\al
\leq e_{j_{r+1}-1}\al=e_{i_r}$.
Hence, $(V(\tr_j))\al=\{e_{i_r}\}$.

Conversely, suppose that conditions (1) and (2) hold.
Let $z\in\xn$, and choose $j\in[m]$ such that $z\in V(\tr_j)$.
Since $j_1=1$ and $j_{t+1}=m+1$, there exists $r\in[t]$ such that
$j_r\leq j<j_{r+1}$.
If $j<j_{r+1}-1$, then $z(\al\bt)=(z\al)\bt=e_{i_r}=(z\bt)\al=z(\bt\al)$
by 2(b) and the fact that $z\bt\in V(\tr_j)$. Suppose $j= j_{r+1}-1$, and
let $y=z\bt$ and $x=y\al$. Then, $y\in V(\tr_j)\cap\im(\bt)$,  and so, by 2(a), $x\in V(\tr_{i_r})\cap\im(\bt)$
and $z\al\in x\bt^{-1}$. Hence, $z(\al\bt)=(z\al)\bt=x=y\al=(z\bt)\al=z(\bt\al)$. Therefore,
$\al\in C(\bt)$.
\end{proof}

In Subsection~\ref{snil5}, we observed that for a nilpotent $\mu\in\cn$,
$C(\mu)=C_{id}(\mu)$ if and only if $n\leq2$. The corresponding result for a general
$\bt\in\cn$, with $\ga(\bt)=\{\tr_1<\cdots<\tr_m\}$, is almost equally obvious:
$C(\bt)=C_{id}(\bt)$ if and only if $n\leq2$ or $n=3$ and $m=2$.

We now describe $\bt\in\cn$ for which $C(\bt)$ consists of nilpotents together with the identity
or of nilpotents and idempotents. (See Theorems~\ref{tiii} and~\ref{tiv} obtained for nilpotents.)

The following theorem shows that only nilpotents $\bt$ can have a centralizer without nontrivial idempotents
or be such that $C(\bt)=\cnl(\bt)\cup\{\vep_n\}$.

\begin{theorem}\label{tgeniii}
	Let $n\geq3$ and $\bt\in\cn$. Then, the following conditions are equivalent:
	\begin{itemize} 
		\item [\rm(1)] $C(\bt)$ contains no idempotents except $\vep_n$ and $\thet_n$;
		\item [\rm(2)] $C(\bt)=\cnl(\bt)\cup\{\vep_n\}$;
		\item [\rm(3)] $\bt$ is a nilpotent, and $\tr=\ga(\bt)$ satisfies $|V_1(\tr)|=1$ and, for all $x,y,z\in\xn$ with $x\ne1$, if $y,z\in V_1(\tr_x)$ and $y<z$, then $\hei(\tr_y)<\hei(\tr_z)$.
	\end{itemize}
\end{theorem}
\begin{proof}
The proof of ($(1)\Rightarrow (2)$)  is exactly the same as that of ($(1)\Rightarrow (2)$) in Theorem~\ref{tiii}, with $\mu$ replaced by $\bt$.
	
	($(2)\Rightarrow (3)$) Assume (2) holds. Since $\bt\in C(\bt)$, either $\bt$ is a nilpotent or $\bt=\vep_n$.
Suppose to the contrary that $\bt=\vep_n$. Then, $C(\bt)=\cn$ and $m=n\geq3$. Define $\al:\xn\to\xn$ by $2\al=1$ and 
$p\al=p$ for every $p\ne2$. Then, $\al\in C(\bt)$, $\al\notin \cnl(\bt)$ (since $3\al=3$),
and $\al\ne\vep_n$ (since $2\al=1$). This contradicts~(2). Hence $\bt$ is a nilpotent,
and so (3) holds by Theorem~\ref{tiii}.
	
Finally, ($(3)\Rightarrow (1)$) follows from Theorem~\ref{tiii}.
\end{proof}

The next theorem complements Theorem~\ref{tiv} by treating non-nilpotent transformations.

\begin{theorem}\label{tgeniv}
	Let $n\geq2$ and $\bt\in\cn$ with $\ga(\bt)=\{\tr_1<\cdots<\tr_m\}$, where $m\geq2$. Then, $C(\bt)=\cnl(\bt)\cup\cid(\bt)$
	if and only if $m\in\{2,3\}$ and the following conditions are satisfied:
	\begin{itemize}
		\item[\rm(1)] $|V(\tr_j)|\leq2$ for every $j\in[m]$;
		\item[\rm(2)] if $m=3$, then either $|V(\tr_2)|=1$ or $|V(\tr_2)|=2$ and $|V(\tr_1)|=|V(\tr_3)|=1$.
	\end{itemize}
\end{theorem}
\begin{proof}
	$(\Rightarrow)$ Suppose $C(\bt)=\cnl(\bt)\cup\cid(\bt)$. Then, $\bt$ is not a nilpotent (since $m\geq2$), and so
$\bt$ is an idempotent since $\bt\in C(\bt)$. Let $e_i$ be the root of $\tr_i$ for every $i\in[m]$.

Suppose to the contrary that $m\geq4$. 
Define $\al:\xn\to\xn$ by: $p\al=1$ if $p\in V(\tr_1)\cup V(\tr_2)$, $p\al=e_2$ if $p\in V(\tr_3)$, and $p\al=p$ for every other $p$.
Then, $\al\in\cn$ (since it is decreasing and order-preserving), $\al\in C(\bt)$ (by Theorem~\ref{tcide}), 
$\al$ is not a nilpotent (since $e_4\al=e_4\ne1$),
and $\al$ is not an idempotent (since $e_3\al=e_2$ and $e_2\al=1\ne e_2$). This contradicts $C(\bt)=\cnl(\bt)\cup\cid(\bt)$.
Hence, $m\in\{2,3\}$.

Suppose to the contrary that for some $j\in[m]$, $|V(\tr_j)|\geq3$, that is, $V(\tr_j)=\{e_j<x_1<\cdots<x_k\}$ for some $k\geq2$.
Note that $x_1,\ldots,x_k\in V_1(\tr_j)$ (since $\bt$ is an idempotent). Define $\al:\xn\to\xn$ by $e_j\al=x_1\al=e_j$, $x_2\al=x_1$,
and $p\al=p$ for every other $p$. Then, $\al\in\cn$, $\al\in C(\bt)$, $\al$~is not a nilpotent (since $e_2\al=e_2\ne1$),
and $\al$ is not an idempotent (since $x_2\al=x_1$ and $x_1\al=e_j\ne x_1$). This contradicts $C(\bt)=\cnl(\bt)\cup\cid(\bt)$,
which proves (1).

Suppose $m=3$ and $|V(\tr_2)|\ne1$. By (1), $|V(\tr_2)|=2$ and $|V(\tr_1)|,|V(\tr_3)|\leq2$.
Let $V(\tr_2)=\{e_2<y\}$.  Suppose to the contrary that $|V(\tr_1)|=2$, and let $V(\tr_1)=\{1<x\}$.
Define $\al:\xn\to\xn$ by: $1\al=x\al=e_2\al=1$, $y\al=x$, and $p\al=p$ for every $p\in V(\tr_3)$.
Then, $\al\in\cn$, $\al\in C(\bt)$, $\al$~is not a nilpotent (since $e_3\al=e_3\ne1$),
and $\al$ is not an idempotent (since $y\al=x$ and $x\al=1\ne x$). This contradicts $C(\bt)=\cnl(\bt)\cup\cid(\bt)$.

Suppose to the contrary that $|V(\tr_3)|=2$, and let $V(\tr_3)=\{e_3<z\}$.
Define $\al:\xn\to\xn$ by: $e_2\al=y\al=e_3\al=e_2$, $z\al=y$, and $p\al=p$ for every $p\in V(\tr_1)$.
Then, $\al\in\cn$, $\al\in C(\bt)$, $\al$~is not a nilpotent (since $e_2\al=e_2\ne1$),
and $\al$ is not an idempotent (since $z\al=y$ and $y\al=e_2\ne y$). This again contradicts $C(\bt)=\cnl(\bt)\cup\cid(\bt)$.
We have proved~(2).
	
	$(\Leftarrow)$ Suppose that $m\in\{2,3\}$ and (1) and (2) are satisfied. By (1), $\bt$ is an idempotent.
Suppose $\al\in C(\bt)$.

Let $m=2$. By Theorem~\ref{tcide}, $(V(\tr_1))\al\subseteq V(\tr_1)$ and $e_2\al\in\{1,e_2\}$.
If $e_2\al=1$, we have $(V(\tr_1))\al=\{1\}$ (since $\al$ is order-preserving and $\tr_1<\tr_2$), and $(V(\tr_2))\al\subseteq V(\tr_1)$.
Hence, $\al$ is a nilpotent. If $e_2\al=e_2$, then $(V(\tr_2))\al\subseteq V(\tr_2)$, and so $\al$ is an idempotent.

Let $m=3$, and set $E_j=V(\tr_j)$ for $j\in[3]$. By Theorem~\ref{tcide}, the possible pairs $(e_2\al,e_3\al)$ are
$(1,1)$, $(1,e_2)$, $(1,e_3)$, $(e_2,e_2)$, and $(e_2,e_3)$. 
Suppose $e_2\al=1$ and $e_3\al<e_3$. Then $E_1\al=\{1\}$, $E_2\al\subseteq E_1$, and $E_3\al\subseteq E_1\cup E_2$.
Thus, $\fix(\al)=\{1\}$, and so $\al$ is a nilpotent. 
Suppose $e_2\al=1$ and $e_3\al=e_3$.  Then
$E_2\al =\{1\}$ (if $|E_2| = 2$, then $|E_1| = 1$ by (2)), and so $(E_1\cup E_2)\al = \{1\}$.
Therefore, $\al$ is an idempotent. 
Suppose $e_2\al=e_3\al=e_2$. Since $\al$ is order-preserving and $E_2<E_3$, we have $E_2\al=\{e_2\}$. By (2), either $|E_2|=1$ or $|E_3|=1$,
so $E_3\al=\{e_2\}$. It follows that $\al$ is an idempotent. 
Finally, if $e_2\al=e_2$ and $e_3\al=e_3$, then $\al$ is an idempotent.

The foregoing argument, together with (1) and (2), implies that $C(\bt)=\cnl(\bt)\cup\cid(\bt)$.
\end{proof}

\begin{example}\label{exgeniv1}
By Theorem~\ref{tgeniv}, the non-nilpotent transformations $\bt\in\mathcal{C}_4$ satisfying
$C(\bt)=\cnl(\bt)\cup\cid(\bt)$ are $1133$, $1134$, $1224$, and $1233$ (see Figure~\ref{fig4}).
In $\mathcal{C}_5$, the only such transformation is $11344$ (see Figure~\ref{fig5}). For $n\geq6$, there are no such transformations.
\end{example}

\begin{figure}[htbp]
	\centering
	\begin{tikzpicture}[
	level distance=1.0cm,
	level 1/.style={sibling distance=0.9cm},
	every node/.style={circle, fill=black, inner sep=1.2pt}
	]
	\node[label=above:1] at (0,0) {}
	child { node[label=below:2] {} };
	\node[label=above:3] at (0.9,0) {}
	child { node[label=below:4] {} };
	\node[label=above:1] at (3.4,0) {}
	child { node[label=below:2] {} };
	\node[label=above:3] at (4.3,0) {};
	\node[label=above:4] at (5.2,0) {};
	\node[label=above:1] at (6.8,0) {};
	\node[label=above:2] at (7.7,0) {}
	child { node[label=below:3] {} };
	\node[label=above:4] at (8.6,0) {};
	\node[label=above:1] at (10.2,0) {};
	\node[label=above:2] at (11.1,0) {};
	\node[label=above:3] at (12.0,0) {}
	child { node[label=below:4] {} };
	\end{tikzpicture}
	\caption{The non-nilpotent $\bt\in\mathcal{C}_4$ such that $C(\bt)=\cnl(\bt)\cup\cid(\bt)$.}
	\label{fig4}
\end{figure}
\begin{figure}[htbp]
	\centering
	\begin{tikzpicture}[
	level distance=1.0cm,
	level 1/.style={sibling distance=0.9cm},
	every node/.style={circle, fill=black, inner sep=1.2pt}
	]
	\node[label=above:1] at (0,0) {}
	child { node[label=below:2] {} };
	\node[label=above:3] at (0.9,0) {};
	\node[label=above:4] at (1.8,0) {}
	child { node[label=below:5] {} };
	\end{tikzpicture}
	\caption{The only non-nilpotent $\bt\in\mathcal{C}_5$ such that $C(\bt)=\cnl(\bt)\cup\cid(\bt)$.}
	\label{fig5}
\end{figure}

\begin{example}\label{exgeniv}
	Consider $\bt=1222\in\mathcal{C}_4$ (see Figure~\ref{fig7}). Since (1) of Theorem~\ref{tgeniv} fails, 
$C(\bt)\ne\cnl(\bt)\cup\cid(\bt)$. Here, $\al=1223\in C(\bt)$ is neither a nilpotent nor an idempotent.
	
	Consider $\gam=11335\in\mathcal{C}_5$ (see Figure~\ref{fig7}). Since (2) of Theorem~\ref{tgeniv} fails,
$C(\gam)\ne\cnl(\gam)\cup\cid(\gam)$. Here, 
$\al=11125\in C(\gam)$ is neither a nilpotent nor an idempotent.
\end{example}
\begin{figure}[htbp]
	\centering
	\begin{tikzpicture}[
	level distance=1.0cm,
	level 1/.style={sibling distance=0.9cm},
	every node/.style={circle, fill=black, inner sep=1.2pt}
	]
	\node[label=above:1] at (0,0) {};
	\node[label=above:2] at (1.35,0) {}
	child { node[label=below:3] {} }
	child { node[label=below:4] {} };
	\node[label=above:1] at (4.5,0) {}
	child { node[label=below:2] {} };
	\node[label=above:3] at (5.4,0) {}
	child { node[label=below:4] {} };
	\node[label=above:5] at (6.3,0) {};
	\end{tikzpicture}
	\caption{The transformations $\bt$ and $\gam$ in Example~\ref{exgeniv}.}
	\label{fig7}
\end{figure}

\end{document}